\documentclass[12pt,a4paper]{article}

\usepackage[latin2]{inputenc}
\usepackage[T1]{fontenc}
\usepackage[english]{babel}
\usepackage{amsmath}
\usepackage{amssymb}
\usepackage{graphicx}
\usepackage{tikz}
\usepackage{caption}
\usepackage{subcaption}
\usepackage{float}
\usepackage{color}
\usepackage{enumitem}
\usetikzlibrary{patterns}
\newtheorem{thr}{Theorem}
\newtheorem{dfn}{Definition}
\newtheorem{lmm}{Lemma}
\newtheorem{prp}{Proposition}
\newtheorem{crl}{Corollary}
\newtheorem{clm}{Claim}

\newtheorem{obs}{Observation}
\newtheorem{cnj}{Conjecture}
\newtheorem{prb}{Problem}
\newenvironment{prf}{\noindent\textit{Proof.}\relax}{\vspace{3mm}}%
\newenvironment{case}[1]{\vspace{3mm}\noindent\textsc{Case} #1.\relax}{}%
\newenvironment{prfof}[1]{\vspace{3mm}\noindent\textit{Proof of #1.}\relax}{\vspace{3mm}}%

\numberwithin{equation}{section}
\makeatletter
\newcommand*{\storecounter}[2]{%
  \edef\@currentlabel{\the\value{#1}}% Store current counter value in \@currentlabel
  \label{#2}% Store label
}
\makeatother

\title{Betweenness Structures of Small Linear Co-Size}
\date{\today}
\author{P\'eter G.N. Szab\'o\thanks{This work was supported by the National Research, Development and Innovation Office -- NKFIH,  No. 108947.}\\
Email: szape@cs.bme.hu\\
\\
Alfr\'ed R\'enyi Institute of Mathematics\\
Hungarian Academy of Sciences\\
Budapest, H-1364, Hungary\\
\textit{and}\\
Department of Computer Science and Information Theory\\
Budapest University of Technology and Economics\\
Budapest, H-1111, Hungary}

\begin{document}

\maketitle

%%% utilize notations tau and theta???

%%% final checks:
%% proof validation - check
%% wording check - check
%% spell check - check
%% correct comma usage - check
%% format check - check
%% ask Gyula - check
%% bibtex - check
%% export figures - check
%% elsevier format

%%% shortening plan:
%% put appendix into a separate paper; put more lemmas in the appendix - check
%% reduce size of figures/delete unnecessary figures - \includegraphics[scale=0.5]{...} - check
%% investigate if changing order of theorems helps - check
%% investigate if stating the theorems only for large n helps
%% shorten proofs where possible

\begin{abstract}
One way to study the combinatorics of finite metric spaces is to study the betweenness relation associated with the metric space.
%The notion of almost metrizable (pseudometric) betweenness was introduced as a combinatorial generalization of metric betweenness.
%The co-size of a betweenness is the number of its (non-degenerate) triangles.
In the hypergraph metrization problem, one has to find and characterize metric betweennesses whose collinear triples (or alternatively, non-degenerate triangles) coincide with the edges of a given $3$-uniform hypergraph.
Metrizability of different kinds of hypergraphs was investigated in the last decades. Chen showed that steiner triple systems are not metrizable, while Richmond and Richmond characterized linear betweennesses, i.e. metric betweennesses that realize the complete $3$-uniform hypergraph. The latter result was also generalized to almost-metric betweennesses by Beaudou et al. In this paper, we further extend this theory by characterizing the largest nonlinear almost-metric betweennesses that satisfy certain hereditary properties, as well as the ones that contain a small linear number of non-degenerate triangles.
\end{abstract}


%---------------------------------------------------------------------------------------------------

\section{Introduction}\label{Sint}

Metric space is one of the most successful concepts of mathematics, with a wide range of applications in many fields including computer science, quantitative geometry, topology, molecular chemistry and phylogenetics.
Although finite metric spaces are trivial objects from a topological point of view, they have surprisingly complex and intriguing combinatorial properties, which were investigated from different angles over the last fifty years \cite{aboulker2015chen, bandelt1992canonical, beaudou2015debruijn, buneman1974note, chen2008problems, chvatal2014debruijn, dress1987parsimonious, mascioni2004equilateral}.

A well known approach to the combinatorics of finite metric spaces is to study the betweenness relation associated with the metric space.
We say that point $y$ is between points $x$ and $z$ in a metric space $M = (X, d)$ if $$d(x, y) + d(y, z) = d(x, z).$$ We also say that the triple $\{x, y, z\}$ is collinear. A (non-degenerate) triangle is a triple that is not collinear. The collinear triples/triangles form a $3$-uniform hypergraph called the collinearity/triangle hypergraph of the metric space. A $3$-uniform hypergraph is \emph{metrizable} if it is the collinearity hypergraph of some finite metric space. We can equivalently talk about metrizability of a hypergraph as triangle hypergraph, since it is the complement of the collinearity hypergraph. There are two general types of hypergraph metrization problems.

\begin{prb}
Decide whether a given $3$-uniform hypergraph is metrizable.
\end{prb}

\begin{prb}
Characterize metric betweennesses that realize a giv\-en $3$\--\-uni\-form hypergraph.
\end{prb}

Only a couple of partial results are known to these problems.
% In \cite{chvatal2004sylvester}, Chv\'atal proved results about the uniqueness of Steiner graphs, which imply that no Steiner triple system with more than 3 points is metrizable. A more general result was proved in \cite{chen2003sylvester} by Chen, who showed that no $(v, k, 1)$ design with $k\geq 3$ and $v > k is metrizable. In partucular, no projective plane of order higher than 1, nor any affine plane of order higher than $2$ is metrizable. Finally, Beaudou showed that no complement of a steiner triple system with more than 3 points is metrizable \cite{beaudou2013lines}.
In \cite{chen2003sylvester}, Chen showed as a consequence of the Sylvester-Chv\'atal Theorem that no $(v, k, 1)$ design with $k\geq 3$ and $v > k$ is metrizable. In particular, no Steiner triple system with more than 3 points, finite projective plane of order higher than 1 or finite affine plane of order higher than $2$ is metrizable.
Further, Beaudou et al. proved that no complement of a Steiner triple system with more than 3 points is metrizable \cite{beaudou2013lines}.
%It is also worth mentioning here that metric betweennesses are recognizable in polynomial time (Chv\'atal, \cite{chvatal2004sylvester}).

%In \cite{chen2003sylvester}, Chen showed as a consequence of the Sylvester-Chv\'atal Theorem that Steiner-triple systems and $3$-shadows of finite projective planes are not metrizable.

In \cite{richmond1997metric}, Richmond and Richmond characterized metrizable betweennesses with a complete collinearity hypergraph i.e. with the maximum number of collinear triples. %that contain no (non-degenerate) triangles.
That result was generalized to almost-metrizable (pseudometric) betweennesses by Beaudou et al. \cite{beaudou2013lines}.
The aim of this paper is to further extend these extremal results by characterizing the largest almost-metrizable betweennesses that have at least one triangle. We also characterize infinite families of almost-metrizable betweennesses that have a linear number of triangles.

First, we introduce the base definitions in Section \ref{Sdefs}. Then, we state our main results in Section \ref{Smain}, and prove them in Section \ref{Sprf} after some preparations in Section \ref{Sgen}.
We conclude the paper with some interesting remarks in Section \ref{Sconc}.
 %and with an Appendix that contains some technical parts of the proof.

%---------------------------------------------------------------------------------------------------

\section{Definitions}\label{Sdefs}

A metric space $M = (X, d)$ is \emph{finite} if $|X| <\infty$. In this paper, all metric spaces will be assumed to be finite. Further, a \emph{triple} will always mean an unordered triple if not stated otherwise.
A \emph{betweenness structure} is a pair $\mathcal{B} = (X,\beta)$ where $X$ is a nonempty finite set and $\beta\subseteq X^3$ is a ternary relation called the \emph{betweenness relation} of $\mathcal{B}$.
The \emph{order} of $\mathcal{B}$ is $n(\mathcal{B}) = |X|$.
The relation $(x, y, z)\in\beta$ will be denoted by $(x\ y\ z)_\mathcal{B}$ or simply by $(x\ y\ z)$ if $\mathcal{B}$ is clear from the context and we say that $y$ is \emph{between} $x$ and $z$. We also say that the triple $\{x, y, z\}$ is \emph{collinear} in $\mathcal{B}$. The \emph{size} of $\mathcal{B}$ is the number of collinear triples in $\mathcal{B}$.

%%% ezt talán lejjebb, a linearity-nál kellene leírni???
A non-collinear triple of $\mathcal{B}$ is called a \emph{triangle}. We denote the set of triangles in $\mathcal{B}$ by $\Delta(\mathcal{B})$, and we define the \emph{co-size} of $\mathcal{B}$ to be $|\mathcal{B}|_\Delta = |\Delta(\mathcal{B})|$. We can associate two complementary $3$-uniform hypergraphs to a betweenness structure: the hypergraph of triangles and the hypergraph of collinear triples. In this paper we prefer to use the \emph{triangle hypergraph} of $\mathcal{B}$, denoted by $\mathcal{H}(\mathcal{B})$, as we will study betweenness structures of linear co-size.
Accordingly, when we speak about metrizability of a hypergraph, we mean metrizability as a triangle hypergraph (which is the complement of the hypergraph in the usual definition).
The degree of a point $x\in X$ in $\mathcal{H}(\mathcal{B})$ will be denoted by $d_\mathcal{B}(x)$. % or simply by $d(x)$ if $\mathcal{B}$ is clear from the context.
%%% kell az utolsó megjegyzés???

The \emph{substructure} of $\mathcal{B}$ induced by a nonempty subset $Y\subseteq X$ is the betweenness structure $\mathcal{B}\vert_Y = (Y,\beta\cap Y^3)$.
%We write $\mathcal{A}\leq\mathcal{B}$ if $\mathcal{A}$ is a substructure of $\mathcal{B}$.
The substructure induced by $X\backslash\{x\}$ will also be denoted by $\mathcal{B} - x$.
We say that the betweenness structure $\mathcal{B}_1 = (X,\beta_1)$ is an \emph{extension} of the betweenness structure $\mathcal{B}_2 = (X,\beta_2)$ (in notation $\mathcal{B}_1\preccurlyeq\mathcal{B}_2$) if $\beta_1\supseteq\beta_2$ (here, the reversed direction of ``$\succcurlyeq$'' is intentional, as we want the betweenness structure induced by the constant zero pseudometric to be the smallest element in this ordering).

%One way to study the combinatorics of finite metric spaces is to study the betweenness relation associated with the metric space.
The \emph{betweenness structure induced by a finite metric space} $M = (X, d)$ is $\mathcal{B}(M) = (X,\beta_M)$ where $$\beta_M =\{(x, y, z)\in X^3: d(x, z) = d(x, y) + d(y, z)\}$$ is the \emph{betweenness relation} of $M$. Note that $\{x, y, z\}$ is a triangle in $\mathcal{B}(M)$ if and only if the triangle inequality holds strictly for $x$, $y$ and $z$ in any combination.

The betweenness structure $\mathcal{B}$ is \emph{metrizable} if it is induced by some  metric space $M = (X, d)$. We say that a 3-uniform hypergraph $\mathcal{H}$ is \emph{metrizable} if there exists a metrizable betweennes structure $\mathcal{B}$ such that $\mathcal{H} =\mathcal{H}(\mathcal{B})$.
The betweenness relation of a metrizable betweenness structure satisfies the following elementary properties for all $x, y, z\in X$:
\begin{enumerate}[label={(P\arabic*)}, ref={(P\arabic*)}]
\item $(x\ x\ z)$\label{Ecoll1};
\item $(x\ y\ z)\Rightarrow (z\ y\ x)$\label{Ecoll2};
\item $(x\ y\ z)\wedge (y\ x\ z)\Rightarrow x = y$\label{Ecoll3};
\end{enumerate}
and additionally, for all $x, y, z, w\in X$,
\begin{enumerate}[label={(P\arabic*)}, ref={(P\arabic*)}]
\setcounter{enumi}{3}
\item $(x\ y\ z)\wedge (x\ w\ y)\Rightarrow (x\ w\ z)\wedge (w\ y\ z)$.\label{Ecoll4}
\end{enumerate}
The \emph{trichotomy} of betweenness follows straight from properties \ref{Ecoll1}--\ref{Ecoll3}: for any three distinct points $x, y, z\in X$, at most one of the relations $(x\ y\ z)$, $(y\ z\ x)$, $(z\ x\ y)$ can hold.
Property \ref{Ecoll4}, that we call the \emph{four relations property} or \emph{f.r.p.} in short, is the simplest non-trivial property of metric betweennesses.

It is easy to see that these elementary properties are not sufficient to guarantee the metrizability of a betweenness structure (think about the Fano plane). We call a betweenness structure \emph{almost-metrizable} if it satisfies properties \ref{Ecoll1}--\ref{Ecoll4}. These betweennesses are usually called ``pseudometric'' in the related literature, however, we want to avoid confusion with a different meaning of the term, a betweenness structure induced by a pseudometric, i.e. a generalized metric where zero distances are allowed.
Quite interestingly, several properties of finite metric spaces can be seamlessly generalized to almost-metrizable betweenness structures (Proposition \ref{Plin} is a good example). Our main results will be stated for almost-metrizable betweennesses, and every betweenness structure will be assumed to be almost-metrizable in the rest of the paper if not stated otherwise.

The \emph{adjacency graph of a betweenness structure} $\mathcal{B} = (X,\beta)$ is the simple graph $G(\mathcal{B}) = (X, E)$ where the edges are such pairs of points for which no third point lies between them, %(see Figure \ref{Fadj})
or more formally,
$$E(\mathcal{B}) =\left\{\{x,z\}\in\binom{X}{2}:\nexists\,y\in X\backslash\{x, z\},\,(x\ y\ z)_\mathcal{B}\right\}.$$
%These edges are also called primitive pairs by some authors.
Further, the \emph{adjacency graph of a finite metric space} $M$ is defined to be $G(M) = G(\mathcal{B}(M))$. We can make the following observations about the adjacency graph.

%\begin{figure}[t]
%\centering
%\includegraphics[scale=0.75]{fms_qus_v3-figure0}
%%\begin{tikzpicture}[scale=1, pont/.style={circle, fill=black, inner sep=0.5mm}]
%%
%%\draw[->] (-0.75,0)--(3,0) node[anchor=north] {{\small $x$}} ;
%%\draw[->] (0,-0.75)--(0,3) node[anchor=east] {{\small $y$}};
%%\draw[dashed] (-0.75,1.5)--(3,1.5);
%%\draw[dashed] (1.5,-0.75)--(1.5,3);
%%
%%\node (A) at (0.75,1.5) [pont] {};
%%\node (B) at (1.5,1.5) [pont] {};
%%\node (C) at (2.25,1.5) [pont] {};
%%\node (D) at (1.5,0.75) [pont] {};
%%\node (E) at (1.5,2.25) [pont] {};
%%\draw (2.625,2.625) node {$\mathcal{B}$};
%%
%%\node (O) at (3.75,-0.75) {};
%%\node (F) at (5.25,1.5) [pont] {};
%%\node (G) at (6,1.5) [pont] {};
%%\node (H) at (6.75,1.5) [pont] {};
%%\node (I) at (6,0.75) [pont] {};
%%\node (J) at (6,2.25) [pont] {};
%%\draw (7.125,2.625) node {$G(\mathcal{B})$};
%%
%%\draw[thick] (F)--(H);
%%\draw[thick] (I)--(J);
%%\draw[thick] (F)--(I);
%%\draw[thick] (I)--(H);
%%\draw[thick] (H)--(J);
%%\draw[thick] (J)--(F);
%%\draw[->,decorate,decoration={snake, amplitude = 0.4mm, segment length = 2.25mm, post length = 0.5mm}](3.75,1.5)--(4.5,1.5);
%%
%%\end{tikzpicture}
%\caption{The adjacency graph of a betweenness structure induced by five points in the Euclidean plane}\label{Fadj}
%\end{figure}

\begin{obs}%\label{Oconn}
The adjacency graph of a betweenness structure is con\-nec\-ted.
\end{obs}

\begin{obs}\label{Osubmetr}
Let $\mathcal{B}$ be a betweenness structure and let $Y$ be a nonempty set of points in $\mathcal{B}$. Then $G(\mathcal{B})[Y]\leq G({\mathcal{B}\vert_Y})$.
\end{obs}

%------

A \emph{weighted graph} is a triple $W = (V, E,\omega)$ where $G = (V, E)$ is a simple graph and $\omega$ is a positive real-valued function on the set of edges, also called the \emph{edge weighting} of $W$. We note that every simple graph $G = (V, E)$ can be regarded as a weighted graph with $\omega\equiv 1$ as edge weighting. We will freely move between these interpretations as convenient.

By \emph{graph} we will mean a connected weighted graph in the rest of the paper if not stated otherwise. We will write ``\emph{simple graph}'' if we want to emphasize that all of the edge weights are equal to one.
%Further, every graph will be assumed to be connected.
We use notations $P_n$, $C_n$, $K_n$ and $K_{n_1, n_2}$ in the usual sense for the (non-weighted) path, cycle, complete graph of order $n$ and for the complete bipartitie graph with parts of size $n_1$ and $n_2$, respectively.

Let $W$ be a graph on vertex set $V$. The length of a path in $W$ is the sum of the weights on its edges.
The \emph{metric space induced by $W$} is $M(W) = (V, d_W)$ where $d_W$ is the usual \emph{graph metric} of $W$, i.e. $d_W(u, v)$ is the length of the shortest path between $u$ and $v$ in $W$.
Note that every finite metric space $M = (X, d)$ is induced by some graph $W$. For example, take $d$ as the edge weighting on a complete simple graph on vertex set $X$.
It can be also proved that the adjacency graph is the smallest simple graph that can induce the metric space with an appropriate edge weighting.

The \emph{betweenness structure induced by $W$} is the betweenness structure induced by $M(W)$, denoted by $\mathcal{B}(W)$. We also say that $W$ is the spanner graph of the betweenness structure $\mathcal{B}(W)$.
 %If betweenness structure $\mathcal{B}$ is induced by graph $W$, then we say that $W$ is a \emph{spanner graph} of $\mathcal{B}$.
In order to simplify notations, we will write $(x\ y\ z)_W$ instead of $(x\ y\ z)_{\mathcal{B}(W)}$.
Note that $\mathcal{B}(W)$ is always metrizable, and $(x\ y\ z)_{W}$ holds if and only if $y$ is on a shortest path connecting $x$ and $z$ in $W$.
A betweenness structure (or finite metric space) is
\begin{itemize}
\item\emph{graphic} if it is induced by a simple graph;
\item\emph{ordered} if it is induced by a path;
\item\emph{orderable} if it has an ordered extension.
\end{itemize}
We remark that betweenness structures are typically not graphic.

%------

%Let $G = (V, E)$ be a simple graph.
%The \emph{metric space induced by $G$} is $M(G) = (V, d_G)$ where $d_G$ is the usual \emph{graph metric} of $G$, i.e. $d_G(u, v)$ is the length of the shortest path between $u$ and $v$ in $G$.
%The \emph{betweenness structure induced by $G$} is the betweenness structure induced by $M(G)$, also denoted by $\mathcal{B}(G)$. In order to simplify notations, we will write $(x\ y\ z)_G$ instead of $(x\ y\ z)_{\mathcal{B}(G)}$.
%Note that $\mathcal{B}(G)$ is always metrizable and $(x\ y\ z)_{\mathcal{B}(G)}$ holds if and only if $y$ is on a shortest path connecting $x$ and $z$ in $G$.
%A betweenness structure (or finite metric space) is
%\begin{itemize}
%\item\emph{graphic} if it is induced by a simple graph;
%\item\emph{ordered} if it is induced by a path;
%\item\emph{orderable} if it has an ordered extension.
%\end{itemize}
%If the graphic betweenness structure $\mathcal{B}$ is induced by the simple graph $G$, then we say that $G$ is the \emph{spanner graph} of $\mathcal{B}$. We remark that betweenness structures are typically not graphic.

We will denote the ordered betweenness structure induced by the path $P = x_1x_2\ldots x_n$ by $[x_1, x_2,\ldots, x_n]$.
Let $\mathcal{B}$ be a betweenness structure on ground set $X$ and let $Y =\{y_1, y_2,\ldots, y_\ell\}$ be a subset of $X$. We write $(y_1\ y_2\ \ldots\ y_\ell)_\mathcal{B}$ if $\mathcal{B}\vert_Y = [y_1, y_2,\ldots, y_\ell]$. Notice that for three points, this gives back the usual notation of betweenness.
%The statement $\mathcal{B}\vert_Y = [y_1, y_2,\ldots, y_\ell]$ will be denoted by $(y_1\ y_2\ \ldots\ y_\ell)_\mathcal{B}$ (which, for three points, gives back the usual notation for betweennesses).
We close this section with an observation that is an easy consequence of f.r.p.
\begin{obs}\label{Oblowup}
Let $\mathcal{B}$ be a betweenness structure on ground set $X =\allowbreak\{x_1,\allowbreak x_2,\allowbreak\ldots,\allowbreak x_n\}$ and let $1\leq i\leq j\leq n$ be integers such that for $Y =\allowbreak\{x_i,\allowbreak x_{i + 1},\allowbreak\ldots,\allowbreak x_j\}$ and $Z =\allowbreak\{x_1,\allowbreak\ldots,\allowbreak x_i,\allowbreak x_j,\allowbreak\ldots,\allowbreak x_n\}$, $\mathcal{B}\vert_Y =\allowbreak [x_i,\allowbreak x_{i + 1},\allowbreak\ldots,\allowbreak x_j]$ and $\mathcal{B}\vert_Z =\allowbreak [x_1,\allowbreak x_2,\allowbreak\ldots,\allowbreak x_i,\allowbreak x_j,\allowbreak\ldots,\allowbreak x_n]$. Then
$$\mathcal{B} =[x_1, x_2,\ldots, x_n].$$
%Let $\mathcal{B}$ be betweenness structure and let $y$ be a point of $\mathcal{B}$ such that
%$\mathcal{B} - y = [x_1, x_2,\ldots, x_{n - 1}]$ and $(x_i\ y\ x_{i + 1})$ holds. Then
%$$\mathcal{B} =\allowbreak [x_1,\allowbreak x_2,\allowbreak\ldots,\allowbreak x_i,\allowbreak y,\allowbreak x_{i + 1},\allowbreak\ldots,\allowbreak x_{n - 1}].$$
\end{obs}

Two important subcases of Observation \ref{Oblowup} that we will extensively use later are $|Y| = 3$ and $|Z| = 3$. Note that we get back f.r.p. by setting $|Y| = |Z| = 3$.

%A \emph{weighted graph} is a triple $W = (V, E,\omega)$ where $G = (V, E)$ is a simple graph and $\omega$ is a positive real-valued function on the set of edges, also called the \emph{edge weighting} of $W$. We note that every simple graph $G = (V, E)$ can be regarded as a weighted graph $W(G)$, with edge weighting $\omega_G =\mathbf{1}_{E}$ (where $\mathbf{1}_A$ denotes the indicator function of the set $A\subseteq E$). We will freely move between interpretations as convenient.
%
%Similarly to the definitions above, we can define the metric space $M(W)$ and the betweenness structure $\mathcal{B}(W)$ induced by a weighted graph $W$ where the ``length'' of a path is the sum of the weights on its edges. We remark that these definitions are compatible with the corresponding definitions for simple graphs. Also note that every finite metric space $M = (X, d)$ is induced by some weighted graph $W$. For example, take $d$ as the edge weighting on a complete graph on ground set $X$.
%It can also be proved that the adjacency graph is the smallest simple graph that can induce the metric space with an appropriate edge weighting.

%---------------------------------------------------------------------------------------------------

\section{Main Results}\label{Smain}
%We use the usual notations $P_n$ and $C_n$ for the path and the cycle of length $n$, respectively.
Let $\mathcal{P}_n$ and $\mathcal{C}_n$ denote the graphic betweenness structures induced by $P_n$ and $C_n$, respectively.
\begin{dfn}
A betweenness structure $\mathcal{B} = (X,\beta)$ is \textit{linear} if any triple $T\in\binom{X}{3}$ is collinear, or equivalently, if $\mathcal{B}$ is of co-size $0$.
\end{dfn}

Observe that all ordered betweenness structures are linear and any substructures of a linear betweenness structure are linear as well. On the other hand, however, orderedness does not follow from linearity, as $\mathcal{C}_4$ shows.
%It is a natural question to ask, what are the linear betweenness structures up to isomorphism, and whether linearity and orderedness are equivalent.
In \cite{richmond1997metric}, Richmond and Richmond gave a full characterization linear metrizable betweenness structures, which was extended to almost-metrizable betweenness structures by Beaudou et al. in \cite{beaudou2013lines}, Lemma 1. We reformulate the latter result in our notations.
%This question was answered for metrizable betweennesses by Richmond and Richmond in \cite{richmond1997metric},
% and also by Dovgoshei and Dordovskii \cite{dovgoshei2009betweenness}.
%That result was extended to almost-metrizable betweenness structures by Beaudou et al. in \cite{beaudou2013lines}, Lemma 1.

\begin{prp}[Beaudou et al. \cite{beaudou2013lines}]\label{Plin}
Up to isomorphism, the linear be\-tween\-ness structures are $\mathcal{P}_n$ ($n\geq 1$) and $\mathcal{C}_4$.
%A betweenness structure $\mathcal{B}$ is linear if and only if $\mathcal{B}\simeq\mathcal{P}_n$ for some $n\geq 1$ or $\mathcal{B}\simeq\mathcal{C}_4$.
\end{prp}

Because of Proposition \ref{Plin}, we can use linearity and orderedness interchangeably when $n\neq 4$.
A \emph{line} in a betweenness structure $\mathcal{B} = (X,\beta)$ is a set of points $Y\subseteq X$ that induces a linear substructure.
Lines inducing an ordered substructure are called \emph{ordered line}s, while the ones inducing a $\mathcal{C}_4$ are called \emph{cyclic line}s. A betweenness structure is \emph{regular} if it does not contain any cyclic lines.
This seems to be an important distinction as a lot of questions are much easier to answer for regular betweenness structures than for irregular ones.

% Further, the set of betweenness structures of linear co-size $kn - c$ will be denoted by $\Gamma(n, k, c)$.
Let $\Phi$ be a hereditary property of betweenness structures. We will denote the set of betweenness structures of order $n$ and co-size $m$ by $B(n, m)$, and among those, the set of betweenness structures that satisfy $\Phi$ by $B_\Phi(n, m)$.
Our quantities of interest are the following:
\begin{itemize}
\item $\tau(n, k) =\min\{m > k : B(n, m)\neq\emptyset\}$; 
\item $\tau_\Phi(n, k) =\min\{m > k : B_\Phi(n, m)\neq\emptyset\}$
\item $\gamma(k, c) =\max\{n\in\mathbb{Z} : B(n, kn - c)\neq\emptyset\}$;
\item $\sigma(k, c) =\max\{n\in\mathbb{Z} : B(n, kn - c) =\emptyset\}$;
\item $\vartheta_{\min}(k) =\min\{c\in\mathbb{Z} :\gamma(k, c) =\infty\}$;
\item $\vartheta_{\max}(k) =\max\{c\in\mathbb{Z} :\gamma(k, c) =\infty\}$.
\end{itemize}

%\begin{dfn}
%A betweenness structure $\mathcal{B}$ of order $n$ is \textit{quasilinear} if it is a nonlinear betweenness structure of minimum co-size.
%\end{dfn}
%\begin{dfn}
%Let $\tau_n$ denote the minimum co-size of a non-linear betweenness structure of order $n$. A betweenness structure $\mathcal{B}$ of order $n$ is \textit{quasilinear} if $|\mathcal{B}|_\Delta =\tau_n$.
%A betweenness structure $\mathcal{B}$ of order $n$ is \textit{quasilinear} if it is a nonlinear betweenness structure of minimum co-size.
%\end{dfn}
\begin{dfn}
A betweenness structure $\mathcal{B}$ of order $n$ is \textit{quasilinear} if $|\mathcal{B}|_\Delta =\tau(n, 0)$, i.e. it is a nonlinear betweenness structure of minimum co-size.
\end{dfn}
As part of our main results, we will extend Proposition \ref{Plin} by characterizing quasilinear betweenness structures.

\begin{figure}[t]
\centering
\includegraphics{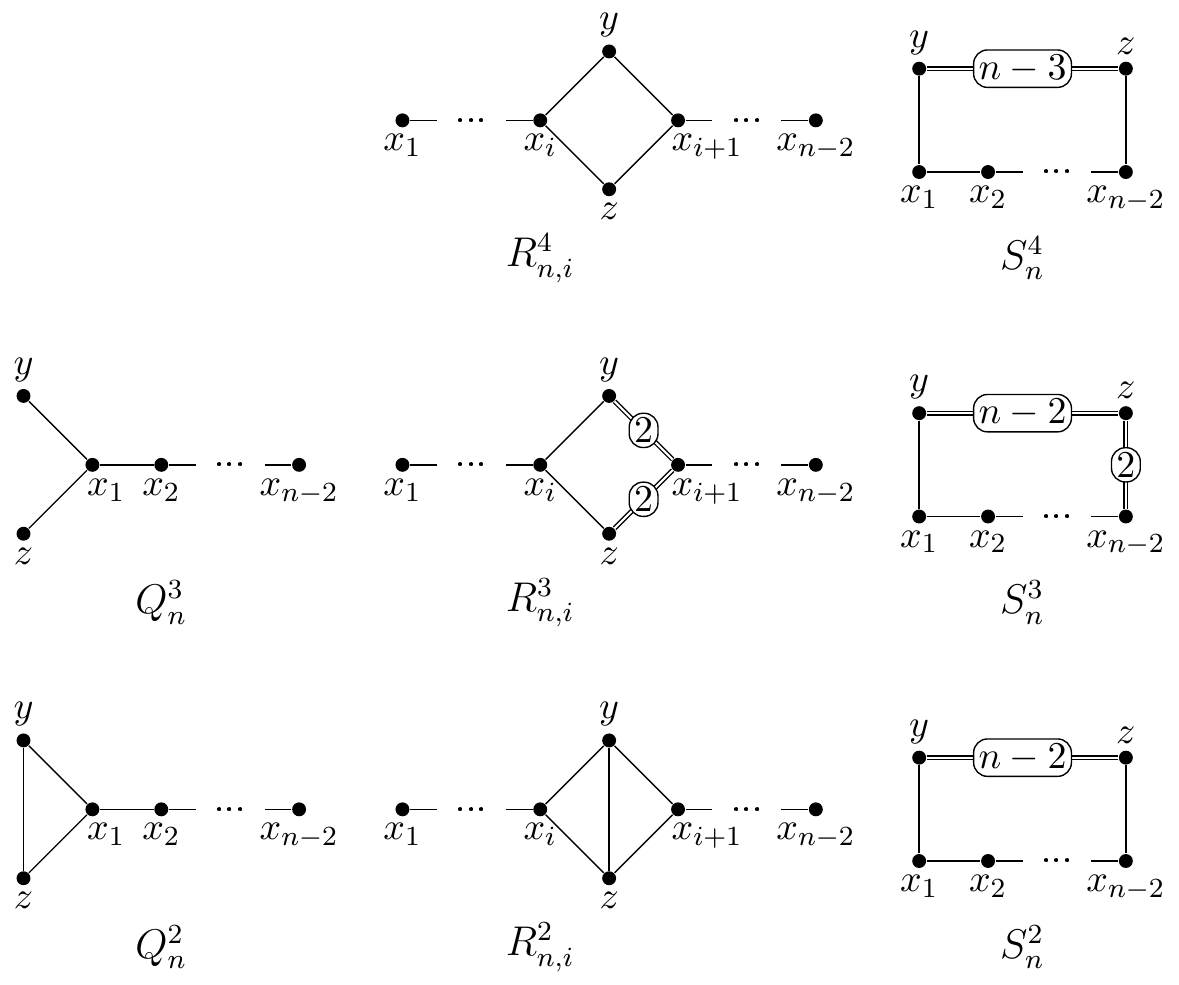}
\caption{Graphs $Q_n^c$ ($2\leq c\leq 3$), $R_{n, i}^c$ ($2\leq c\leq 4$) and $S_n^c$ ($2\leq c\leq 4$). Edges of weight different from $1$ are indicated by double-lines and labeled with the corresponding edge weight.}\label{Fsni}
\end{figure}

Below, we introduce the most important graph classes that will appear in the theorems below (see Figure \ref{Fsni}). All of these graphs are defined on vertex set $X =\{x_1, x_2,\ldots, x_{n - 2}, y, z\}$ relative to the path $P = x_1x_2\ldots x_{n - 2}$. We indicate the admissible values of parameters between parentheses. The range of parameter $i$ will be chosen such that the obtained graphs are pairwise non-isomorphic. These ranges will be $I_n^2 = I_n^4 =\{i\in\mathbb{N} : 1\leq i\leq\left\lceil\frac{n - 3}{2}\right\rceil\}$ and $I_n^3 = \{i\in\mathbb{N} : 1\leq i\leq n - 3\}$.
\begin{itemize}
\item $R_{n, i}^4$ ($n\geq 5$, $i\in I_n^4$): delete edge $\{x_i, x_{i + 1}\}$ from $P$ and add edges $\{y, x_i\}$, $\{z, x_i\}$, $\{y, x_{i + 1}\}$ and $\{z, x_{i + 1}\}$;
\item $S_n^4$ ($n\geq 5$): add edges $\{x_1, y\}$, $\{x_{n - 2}, z\}$, and add edge $\{y, z\}$ of weight $n - 3$;
\item $Q_n^3$ ($n\geq 4$): add edges $\{y, x_1\}$ and $\{z, x_1\}$;
\item $R_{n, i}^3$ ($n\geq 4$, $i\in I_n^3$): delete edge $\{x_i, x_{i + 1}\}$, add edges $\{y, x_i\}$, $\{z, x_i\}$, and add edges $\{y, x_{i + 1}\}$, $\{z, x_{i + 1}\}$ of weight $2$;
\item $S_n^3$ ($n\geq 4$): add edge $\{x_1, y\}$, and add edge $\{x_{n - 2}, z\}$ of weight $2$ and edge $\{y, z\}$ of weight $n - 2$;
\item $Q_n^2$ ($n\geq 3$): add edges $\{y, x_1\}$, $\{z, x_1\}$ and $\{y, z\}$;
\item $R_{n, i}^2$ ($n\geq 4$, $i\in I_n^2$): delete edge $\{x_i, x_{i + 1}\}$, and add edges $\{y, x_i\}$, $\{z, x_i\}$, $\{y, x_{i + 1}\}$, $\{z, x_{i + 1}\}$ and $\{y, z\}$;
\item $S_n^2$ ($n\geq 3$): add edges $\{x_1, y\}$, $\{x_{n - 2}, z\}$, and add edge $\{y, z\}$ of weight $n - 2$.
\end{itemize}

\begin{figure}[t]
\centering
\includegraphics{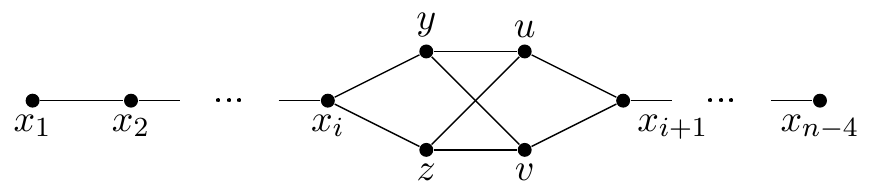}
%\begin{tikzpicture}[pont/.style={circle, fill=black, inner sep=0.5mm},
%mybox/.style={rectangle, rounded corners, draw=black, fill=white, inner sep=0.5mm}]
%
%\node (X1) at (0,0) [pont] {};
%\node (X2) at (1,0) [pont] {};
%\node (XI) at (3,0) [pont] {};
%\node (XI+1) at (6,0) [pont] {};
%\node (XN-4) at (8,0) [pont] {};
%\node (Z) at (4,-0.5) [pont] {};
%\node (Y) at (4,0.5) [pont] {};
%\node (V) at (5,-0.5) [pont] {};
%\node (U) at (5,0.5) [pont] {};
%
%\draw (X1) node[anchor=north] {\small$x_1$};
%\draw (X2) node[anchor=north] {\small$x_2$};
%\draw (XI) node[anchor=north] {\small$x_i$};
%\draw (XI+1) node[anchor=north west] {\small$x_{i+1}$};
%\draw (XN-4) node[anchor=north] {\small$x_{n-4}$};
%\draw (Y) node[anchor=south] {\small$y$};
%\draw (Z) node[anchor=north] {\small$z$};
%\draw (U) node[anchor=south] {\small$u$};
%\draw (V) node[anchor=north] {\small$v$};
%
%\draw (X1)--(X2) node {};
%\draw (X2)--(1.5,0) node {};
%\draw (2,0) node {\small ...};
%\draw (2.5,0)--(XI) node {};
%\draw (XI)--(Y) node {};
%\draw (XI)--(Z) node {};
%\draw (Y)--(U) node {};
%\draw (Y)--(V) node {};
%\draw (Z)--(U) node {};
%\draw (Z)--(V) node {};
%\draw (U)--(XI+1) node {};
%\draw (V)--(XI+1) node {};
%\draw (XI+1)--(6.5,0) node {};
%\draw (7,0) node {\small ...};
%\draw (7.5,0)--(XN-4) node {};
%
%\end{tikzpicture}
\caption{Graph $T_{n, i}$}\label{Ftni}
\end{figure}

Further, we define the simple graph $T_{n, i}$ ($n\geq 6$, $1\leq i\leq\left\lceil\frac{n - 5}{2}\right\rceil$) with vertices $\{x_1, x_2,\ldots x_{n - 4}, y, z, u, v\}$ and edges $\{\{x_j, x_{j + 1}\} : 1\leq j\leq n - 5, j\neq i\}$, $\{x_i, y\},\{x_i, z\},\{y, u\},\{y, v\},\{z, u\},\{z, v\},\{u, x_{i + 1}\}$ and $\{v, x_{i + 1}\}$ (see Figure \ref{Ftni}).
We will denote the betweenness structures induced by graphs $Q_n^c$, $R_{n, i}^c$, $S_n^c$ and $T_{n, i}$ by $\mathcal{Q}_n^c$, $\mathcal{R}_{n, i}^c$, $\mathcal{S}_n^c$ and $\mathcal{T}_{n, i}$, respectively.

We divide our main results into two groups. The starting point of the first three theorems is the characterization of quasilinear betweenness structures. That result can be then easily extended to other extremal problems of similar form: characterize nonlinear betweenness structures of minimum co-size that satisfy certain hereditary properties. We consider two of the most important hereditary properties: \emph{regularity} and \emph{orderability}.

\begin{figure}[t]
\centering
\includegraphics{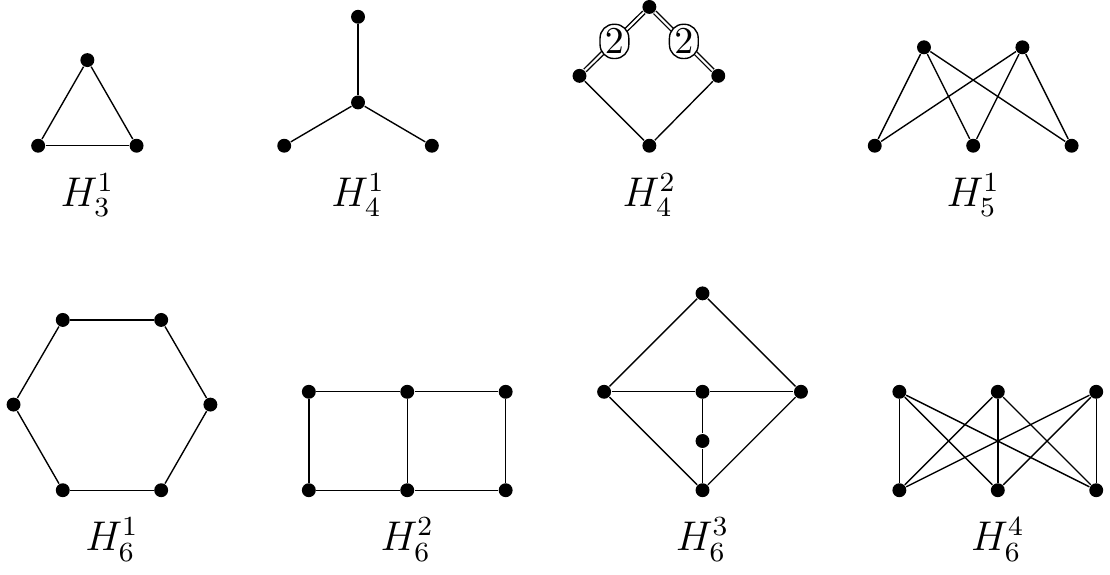}
\caption{Spanner graphs of small exceptional quasilinear betweenness structures. Edges of weight different from $1$ are indicated by double-lines.}\label{Fsmallgr}
\end{figure}

\begin{thr}\label{Tqus1}$ $
\begin{enumerate}
\item For all $n\geq 3$, $\tau(n, 0) =\max\{1, n - 4\}$;\label{Equs11}
\item up to isomorphism, the quasilinear betweenness structures are the following:
\begin{itemize}
\item $\mathcal{R}_{n, i}^4$  for $n\geq 5$, $i\in I_n^4$;
\item $\mathcal{S}_n^4$ for $n\geq 5$;
\item $\mathcal{B}(G)$ where $G$ is one of the graphs in Figure \ref{Fsmallgr}.
\end{itemize}\label{Equs12}
\end{enumerate}
\end{thr}

\begin{thr}\label{Tqus2}
Let $\Phi$ be the ``regular'' property.
Then for all $n\geq 3$, $$\tau_\Phi(n, 0) =\max\{1, n - 3\}.$$ Further, up to isomorphism, the nonlinear regular betweenness structures of minimum co-size are the following:
\begin{itemize}
\item $\mathcal{Q}_n^3$ for $n\geq 4$;
\item $\mathcal{R}_{n, i}^3$ for $n\geq 4$, $i\in I_n^3$;
\item $\mathcal{S}_n^3$ for $n\geq 4$;
\item $\mathcal{B}(K_3)$.
\end{itemize}
\end{thr}

\begin{thr}\label{Tqus3}
Let $\Phi$ be the ``orderable'' property.
Then for all $n\geq 3$, $$\tau_\Phi(n, 0) = n - 2.$$
Further, up to isomorphism, the nonlinear orderable betweenness structures of minimum co-size are the following:
\begin{itemize}
\item $\mathcal{Q}_n^2$ for $n\geq 3$;
\item $\mathcal{R}_{n, i}^2$ for $n\geq 4$, $i\in I_n^2$;
\item $\mathcal{S}_n^2$ for $n\geq 3$.
\end{itemize}
%for all $n\geq 3$, $$\tau_\Phi(n, 0) = n - 2.$$ Further, 
\end{thr}

A $k$-uniform hypergraph $\mathcal{H} = (V,\mathcal{E})$ is a \emph{$\Delta$-star} (or $\Delta$-system) if there exists a set of points $K\subseteq V$ such that for any edges $E, F\in\mathcal{E}$, $E\cap F = K$.
% that is contained in every edge and the edges do not intersect outside of $K$.
% sets $\{E\backslash K : E\in\mathcal{E}\}$ are pairwise disjoint
We call $K$ the \emph{kernel} of the $\Delta$-star and denote it by $\ker(\mathcal{H})$. $\mathcal{H}$ is a \emph{tight star} if it is a $\Delta$-star with kernel of size $k - 1$. Note that the hypergraph with only one edge can be regarded as a tight star. We will apply these definitions to $3$-uniform hypergraphs only.

With the second group of theorems, we focus our attention on betweenness structures of linear co-size $kn - c$. In particular, we fully characterize case $k = 1$ in Theorem \ref{Tlinsize1}, give a sharp upper bound on $c$ in case of $k = 2$ in Theorem \ref{Tlinsize3} and characterize the corresponding extremal betweenness structures in Theorem \ref{Tlinsize2}. Interestingly, there is a series of gaps in the sequence of possible co-sizes that separate the cases with different leading coefficient $k$ if $n$ is large enough. These gaps should be subject to future research.

\begin{thr}\label{Tlinsize1}
Let $c$ be an integer and $N_c = 11 - c$. Then the following hold.
\begin{enumerate}
\item If $c > 4$, then $B(n, n - c) =\emptyset$, except for $n = c$, in which case $B(n, 0)$ consists of all the linear betweenness structures of order $n$ (see Proposition \ref{Plin}).\label{Elinsize11} % \vartheta_{\max}(1)\leq 4
\item If $2\leq c\leq 4$, then $B(n, n - c)\neq\emptyset$ if and only if $n\geq c$. Further, if $n\geq N_c$, then the betweenness structures $\mathcal{B}\in B(n, n - c)$ can be characterized as follows:
\begin{itemize}
\item if $c = 4$, then $\mathcal{B}$ is isomorphic to either $\mathcal{R}_{n, i}^4$ ($i\in I_n^4$) or $\mathcal{S}_n^4$;
\item if $c = 3$, then $\mathcal{B}$ is isomorphic to either $\mathcal{Q}_n^3$, $\mathcal{R}_{n, i}^3$ ($i\in I_n^3$) or $\mathcal{S}_n^3$;
\item if $c = 2$, then $\mathcal{B}$ is isomorphic to either $\mathcal{Q}_n^2$, $\mathcal{R}_{n, i}^2$ ($i\in I_n^2$) or $\mathcal{S}_n^2$.
\end{itemize}\label{Elinsize12}
\item If $c < 2$ and $n\geq N_c$, then $B(n, n - c) =\emptyset$.\label{Elinsize13} % \vartheta_{\min}(1)\geq 2
\item If $c\leq 4$ and $n = N_c - 1$, then $B(n, n - c)\neq\emptyset$. Moreover, there exists a betweenness structure $\mathcal{B}\in B(n, n - c)$ such that $\mathcal{H}(\mathcal{B})$ is not a tight star.\label{Elinsize14}
\end{enumerate}
\end{thr}

\begin{thr}\label{Tlinsize3}
$$\tau(n, n - 2)=\left\{\begin{array}{ll}2n - 10 &\text{if }n\geq 9\\
n - 1 &\text{if }4\leq n\leq 8\end{array}\right.$$
%For all $n\geq 9$, $\tau(n, n - 2) = 2n - 10$.
%For $4\leq n\leq 8$, $\tau(n, n - 2) = n - 1$.
\end{thr}

\begin{thr}\label{Tlinsize2}
$ $
\begin{enumerate}
\item Let $c\geq 11$ be an integer. Then $B(n, 2n - c)\neq\emptyset$ if and only if $n = c / 2$ or $c - 4\leq n\leq c - 2$.\label{Tlinsize21}
\item $B(n, 2n - 10)\neq\emptyset$ if and only if $n\geq 5$.
Further, if $n\geq 9$ and $\mathcal{B}\in B(n, 2n - 10)$, then $\mathcal{B}\simeq\mathcal{T}_{n, i}$ for some $1\leq i\leq\left\lceil\frac{n - 5}{2}\right\rceil$.\label{Tlinsize22}
\item For all $6\leq n < 9$, there exists a betweenness structure $\mathcal{B}\in B(n, 2n - 10)$ such that $\mathcal{B}\not\simeq\mathcal{T}_{n, i}$ for any $1\leq i\leq\left\lceil\frac{n - 5}{2}\right\rceil$.\label{Tlinsize23}
\end{enumerate}
\end{thr}

Now, we can easily determine some of the quantities defined above on the basis of Theorem \ref{Tlinsize1} and Theorem \ref{Tlinsize2}.

\begin{crl}
$ $
\begin{itemize}
\item $\vartheta_{\min}(1) = 2$;
\item $\vartheta_{\max}(1) = 4$;
\item $\sigma(1, c) =\left\{\begin{array}{ll}
c - 1 & \text{if } 2\leq c\leq 4\\
\infty & \text{otherwise}
\end{array}\right.$;
\item $\gamma(1, c) =\left\{\begin{array}{ll}
10 - c & \text{if } c < 2\\
\infty & \text{if } 2\leq c\leq 4\\
c &  \text{if } c > 4\\
\end{array}\right.$.
\end{itemize}
\end{crl}

\begin{crl}
$ $
\begin{itemize}
\item $\vartheta_{\max}(2) = 10$;
\item $\sigma(2, 10) = 4$;
\item $\gamma(2, c) =\left\{\begin{array}{ll}
c - 2 & \text{if } c > 10\\
\infty & \text{if } c = 10
\end{array}\right.$.
\end{itemize}
\end{crl}
% that define the $8$ infinite families of nonlinear betweenness structures of linear co-size.

%--------------------------------------------------------------------------------------------------

\section{General Lemmas}\label{Sgen}

In this section, we state and prove most of the lemmas that we will use in the proof of the main results.
 %serve as building blocks to the main results.

\begin{lmm}\label{Lcyc}
$ $
\begin{itemize}
\item For all $n\geq 6$, $\tau(n, 0)\geq 2$.
\item Additionally, if $\Phi$ is the ``regular'' property, then $\tau_\Phi(5, 0)\geq 2$.
\end{itemize}
\end{lmm}
\begin{prf}
Suppose to the contrary that there exists a nonlinear betweenness structure $\mathcal{B}$ of order $n$ and co-size $1$ such that either $n\geq 6$, or $n = 5$ and $\mathcal{B}$ is regular. Let $X$ be the ground set and $T =\{x, y, z\}$ be the sole triangle of $\mathcal{B}$. Further, let $G$ be the adjacency graph of $\mathcal{B}$ and set $G_p = G(\mathcal{B} - p)$ for all points $p\in T$.

Observe first that for all $p\in T$, \begin{equation}\label{EQ1}
G_p = G - p.
\end{equation}
Since $\mathcal{B}$ is either regular or $n\geq 6$ holds, Proposition \ref{Plin} implies that $\mathcal{B} - p$ is ordered, hence, \begin{equation}\label{EQ1a}
G_p\simeq P_{n - 1}.
\end{equation} Further, $G - p\leq G_p$ by Observation \ref{Osubmetr}, and $G_p\leq G - p$ is also true, otherwise $(u\ p\ v)_\mathcal{B}$ would be true for an edge $\{u, v\}$ of $G_p$, and $\mathcal{B}$ would be linear by Observation \ref{Oblowup}. This completes the proof of (\ref{EQ1}).

Also observe that for all $p\in T$, $\mathcal{B} - p$ is ordered, hence, $\mathcal{B} - p =\mathcal{B}(G_p)$ and we obtain by (\ref{EQ1}) that
\begin{equation}\label{EQ1d}
\mathcal{B} - p =\mathcal{B}(G - p)\simeq\mathcal{P}_{n - 1},
\end{equation}

Our next goal is to show that $G\simeq C_n$. Let $N^{+}_G(w)$ denote the closed neighborhood of a point $w$ in $G$ (i.e. $w\in N^{+}_G(w)$).
%Let $w\neq x$ be a point and let $N^{+}_G(w)$ denote the closed neighborhood of $w$ in $G$.
%As $G_x\simeq P_{n - 1}$, $d_{G_x}(w)\leq 2$ and we obtain by (\ref{EQ1}) that $x$ can be the only extra neighbor of $w$ in $G$. Hence, $d_G(w)\leq 3$ and if $d_G(w) = 3$, then $x$ is a neighbor of $w$. The same argument holds for $w\neq y$ and $w\neq z$, from which it follows that
%for all $w\in X$, $d_G(w)\leq 3$ and if $d_G(w) = 3$, then $T\subseteq N^{+}_G(w)$.

\begin{clm}\label{Cgcic1}
For all $w\in X$, $d_G(w)\leq 3$ and if $d_G(w) = 3$, then $T\subseteq N^{+}_G(w)$.
\end{clm}
\begin{prf}
The degree of a point $w\neq x$ in $G_x\simeq P_{n - 1}$ is at most $2$ and because of (\ref{EQ1}), $x$ can be the only extra neighbor of $w$ in $G$. Hence, $d_G(w)\leq 3$ and if $d_G(w) = 3$, then $x$ is a neighbor of $w$. The same argument holds for $w\neq y$ and $w\neq z$, from which the claim follows. $\square$
\end{prf}

\begin{clm}\label{Cgcic2}
For all $p\in T$, $d_G(p) = 2$.
\end{clm}
\begin{prf}
We can suppose without loss of generality that $p = x$. We obtain from Claim \ref{Cgcic1} that $d_G(x)\leq 3$, thus, it is enough to show that $d_G(x)\neq 3$ and $d_G(x)\neq 1$ (obviously, $d_G(x) > 0$ since $G$ is connected).

Suppose first that $d_G(x) = 3$ and let $u_1, u_2, u_3$ be the neighbors of $x$ such that $(u_1\ u_2\ u_3)_{G_x}$ holds. Then $y = u_2$, otherwise $G_y$ would contain a cycle by (\ref{EQ1}), in contradiction with (\ref{EQ1a}). Similarly, we obtain that $z = u_2$, which contradicts $y\neq z$.

Next, suppose that $d_{G}(x) = 1$. It follows from (\ref{EQ1}) and (\ref{EQ1a}) that $G$ is a tree and \begin{equation}\label{EQ1b}
d_G(y) = d_G(z) = 1,
\end{equation} because otherwise $G_y$ or $G_z$ would be disconnected. Let $w$ be the sole neighbor of $x$. Clearly, $w$ is not a leaf, hence, $w\notin T$. Further, $w$ cannot be an end-vertex of $G_x$, otherwise $G$ would have only two leaves. Therefore, $d_G(w) = 3$ and so $w$ must be adjacent to both $y$ and $z$ by Claim \ref{Cgcic1}. This is, however, impossible since $x, y, z$ and $w$ would form a connected component of $G$ by (\ref{EQ1b}), contradicting $n\geq 5$.
$\square$
\end{prf}

Now, we prove that $$G\simeq C_n.$$
Because of (\ref{EQ1}), it is enough to show that the two neighbors of $x$ guaranteed by Claim \ref{Cgcic2} are the end vertices of $G_x$. Let $u$ and $v$ be the two neighbors of $x$ and suppose to the contrary that $u$ is an inner vertex of $G_x$. Let $P_{uv}$ denote the subpath of $G_x$ that connects $u$ and $v$. Now, $d_G(u) = 3$, so Claim \ref{Cgcic1} implies that $u$ is adjacent to both $y$ and $z$, one of which, e.g. $y$, is not in $P_{uv}$. However, this contradicts (\ref{EQ1a}) as $G_y$ would contain the cycle formed by $P_{uv}$ and the edges $\{x, u\}$ and $\{x, v\}$ by (\ref{EQ1}).
%\end{prf}

%Now, for all $p\in T$, $\mathcal{B} - p$ is ordered, hence, $\mathcal{B} - p =\mathcal{B}(G_p)$ and we obtain by (\ref{EQ1}) that
%\begin{equation}\label{EQ1c}
%\mathcal{B} - p =\mathcal{B}(G - p).
%\end{equation}

%\begin{clm}\label{Cgcic4}
%For all $p\in T$, $\mathcal{B} - p =\mathcal{B}(G - p)$.
%\end{clm}
%\begin{prf}
%Since $\mathcal{B} - p$ is ordered, $\mathcal{B} - p =\mathcal{B}(G_p)$. Observation (\ref{EQ1}) completes the proof.
%%$\mathcal{B} - p$ is a linear betweenness structure because $T$ is the only triangle in $\mathcal{B}$. Every linear betweenness structure is graphic by Proposition \ref{Plin}, hence, $\mathcal{B} - p$ is induced by its adjacency graph $G_p = G - p$ (see (\ref{EQ1})).
%$\square$
%\end{prf}

%Claim \ref{Cgcic3} yields that $G\simeq C_n$.
Let $A_{xy}, A_{yz}$ and $A_{zx}$ denote the three arcs that we obtain by deleting $x, y$ and $z$ from $G\simeq C_n$. Since $n\geq 5$, there exist two distinct points $u$ and $v$ different from $x, y$ an $z$. There are two cases depending on whether $u$ and $v$ are on the same arc.

If $u$ and $v$ are on the same arc, say $A_{xz}$, then we can assume without loss of generality that $(x\ u\ v\ z)_G$ holds. %(see Figure \ref{SFcyc1}).
However, this implies $(u\ v\ y)_{G - x}$ and $(v\ u\ y)_{G - z}$, which lead to a contradiction 
%$(u\ v\ y)_\mathcal{B}$ and $(v\ u\ y)_\mathcal{B}$
by (\ref{EQ1d}).

In the second case when $u$ and $v$ are on distinct arcs, for example $u\in V(A_{xy})$ and $v\in V(A_{yz})$, $(u\ v\ z)_{G - x}$ and $(u\ z\ v)_{G - y}$ hold, %(see Figure \ref{SFcyc2})
which lead to a contradiction again by (\ref{EQ1d}). This completes the proof of Lemma \ref{Lcyc}. $\square$
\end{prf}

\begin{lmm}\label{Lgen1}
Let $\Phi$ be a hereditary property of betweenness structures and suppose that there exists an integer $c\geq 2$ such that $\tau_\Phi(c + 2, 0)\geq 2$. Then for all $n\geq 3$, $\tau_\Phi(n, 0)\geq n - c$.
\end{lmm}
\begin{prf}
We prove $\tau_\Phi(n, 0)\geq n - c$ by induction on $n$.
%First, notice that condition $\tau_\Phi(c + 2, 0)\geq 2$ implies $c\leq 2$.
$\tau_\Phi(n, 0)\geq n - c$ is obvious for $3\leq n\leq c + 1$, and it is also true for $n = c + 2$ by the lemma's assumption.

Next, suppose that $n > c + 2$ and for all $n' < n$, $\tau_\Phi(n', 0)\geq n' - c$.
% for all $n' < n$.
Let $\mathcal{B}$ be a nonlinear betweenness structure of order $n$ that satisfies $\Phi$ and suppose to the contrary that $$|\mathcal{B}|_\Delta < n - c.$$
Notice that $|\mathcal{B}|_\Delta > 1$ because otherwise any substructure of order $c + 2$ of $\mathcal{B}$ that contains the single triangle of $\mathcal{B}$ would violate the assumption of the lemma (here we relied on the assumption that $\Phi$ is hereditary).
Now, let $T_1$ and $T_2$ be two distinct triangles of $\mathcal{B}$, $x$ be a point in $T_1\backslash T_2$ and set $\mathcal{B}' =\mathcal{B} - x$. Since $x\notin T_2$ and $\Phi$ was hereditary, $\mathcal{B}'$ is a nonlinear betweenness structure on $n - 1$ points that satisfies $\Phi$. However, \begin{equation*}\begin{split}
|\mathcal{B}'|_\Delta & \leq |\mathcal{B}|_\Delta - 1\\
& < (n - 1) - c\end{split}\end{equation*} in contradiction with the induction hypothesis.
$\square$
\end{prf}

\begin{obs}\label{Otstar}
Let $c$ and $n$ be integers such that $n > 2c - 1$ and let $\mathcal{B}\in B(n, n - c)$ be a nonlinear betweenness structure on ground set $X$ such that for all points $x\in X$, $d_\mathcal{B}(x) = n - c$ or $d_\mathcal{B}(x)\leq 1$. Then $\mathcal{H}(\mathcal{B})$ is a tight star.
\end{obs}

\begin{prf}
%If $n = c + 1$, then there is only one triangle in $\mathcal{B}$ and $\mathcal{H}(\mathcal{B})$ is clearly a tight star. So, we can assume that $n > c + 1$.
%
%Let $k$ be the number of points $x\in X$ for which $d_\mathcal{B}(x) = n - c$. Clearly, $0\leq k\leq 3$, and it is enough to show that $$k\geq 2.$$
%Counting the number of pairs $\{(x, T) : T\in\Delta(\mathcal{B}), x\in T\}$ in two ways, we obtain $$3(n - c)\leq k(n - c) + (n - k),$$ from which $$k\geq 2 -\frac{c - 2}{n - c - 1}$$ follows by the assumption $n > c + 1$. We complete the proof by showing that $\frac{c - 2}{n - c - 1} < 1$: since the denominator is strictly positive, $\frac{c - 2}{n - c - 1} < 1$ if and only if $c - 2 < n - c - 1$, which is equivalent to the assumption $n > 2c - 1$ of the lemma. $\square$
Let $k$ be the number of points $x\in X$ for which $d_\mathcal{B}(x) = n - c$. It is obvious that $\mathcal{H}(\mathcal{B})$ is a $\Delta$-star and $0\leq k\leq 3$, so it is enough to show that $k\geq 2$.
If $k = 0$, then $|\mathcal{B}|_\Delta\leq n/3$ and if $k = 1$, then $|\mathcal{B}|_\Delta\leq (n - 1)/2$. Since $\mathcal{B}$ is nonlinear, $n\geq 3$ and hence $n/3\leq (n - 1)/2$. Further, $n > 2c - 1$ implies $$(n - 1)/2 < n - c = |\mathcal{B}|_\Delta,$$ thus, $k\geq 2$ and $\mathcal{H}(\mathcal{B})$ is a tight star.
%Let $k$ be the number of points $x\in X$ for which $d_\mathcal{B}(x) = n - c$. It is obvious that $\mathcal{H}(\mathcal{B})$ is a $\Delta$-star and $0\leq k\leq 3$, so it is enough to show that $k\geq 2$.
%Suppose to the contrary that $k\leq 1$.
%Counting the number of pairs $\{(x, T) : x\in T\in\Delta(\mathcal{B})\}$ in two ways, we obtain
%$$3(n - c)\leq k(n - c) + (n - k)$$ and so
%$$3(n - c)\leq k(n - c - 1) + n,$$
%from which
%$$3(n - c)\leq (n - c - 1) + n$$
%follows by the assumptions $k\leq 1$ and $n > c$.
%This, however, contradicts $n > 2c - 1$.
$\square$
\end{prf}

\begin{lmm}\label{Lgen2a}
Let $c$ be an integer and $\Phi$ be a hereditary property of betweenness structures such that $\tau_\Phi(n', 0)\geq n' - c$ for all $n'\geq 3$. Further, let $n > c$ be an integer and $\mathcal{B}\in B_\Phi(n, n - c)$ be a nonlinear betweenness structure. Then $\mathcal{H}(\mathcal{B})$ is a $\Delta$-star. Further, if $n > 2c - 1$, then $\mathcal{H}(\mathcal{B})$ is a tight star.
\end{lmm}
\begin{prf}
First, notice that condition $\tau_\Phi(n', 0)\geq n' - c$ for $n' = 3$ implies that $c\geq 2$. Further, we can assume that $n > c + 1$: if not, then $|\mathcal{B}|_\Delta = n - c\leq 1$ and $\mathcal{H}(\mathcal{B})$ is clearly a $\Delta$-star.

Now, we can prove that $\mathcal{H}(\mathcal{B})$ is a $\Delta$-star by showing that for all points $x\in X$ \begin{equation}\label{EQ1c}
d_\mathcal{B}(x) = n - c\emph{ or } d_\mathcal{B}(x)\leq 1.
\end{equation}
In addition to (\ref{EQ1c}), if $n > 2c - 1$, then $\mathcal{B}$ satisfies the conditions of Observation \ref{Otstar}, thus, $\mathcal{H}(\mathcal{B})$ is also a tight star.

In order to prove (\ref{EQ1c}), suppose to the opposite that there exits a point $x$ such that $1 < d_\mathcal{B}(x) < n - c$. Then $\mathcal{B} - x$ is clearly a nonlinear betweenness structure of order $n - 1\geq c + 1\geq 3$ that satisfies property $\Phi$. Further, $$|\mathcal{B} - x|_\Delta < n - c - 1,$$ which contradicts $\tau_\Phi(n - 1, 0)\geq n - c - 1$. $\square$
\end{prf}

\begin{lmm}\label{Lgen3}
Let $\mathcal{B}\in B(n, n - c)$ be a nonlinear betweenness structure on ground set $X$ such that $\mathcal{H}(\mathcal{B})$ is a tight star. Further, in case of $n = 5$ suppose that $\mathcal{B}$ is regular. Then $2\leq c\leq 4$ and the following hold:
\begin{itemize}
\item if $c = 2$, then $\mathcal{B}$ is isomorphic to either $\mathcal{Q}_n^2$, $\mathcal{R}_{n, i}^2$ ($i\in I_n^2$) or $\mathcal{S}_n^2$ ($\mathcal{R}_{n, i}^2$ is only possible if $n\geq 4$);
\item if $c = 3$, then $\mathcal{B}$ is isomorphic to either $\mathcal{Q}_n^3$, $\mathcal{R}_{n, i}^3$ ($i\in I_n^3$) or $\mathcal{S}_n^3$;
\item if $c = 4$, then $\mathcal{B}$ is isomorphic to either $\mathcal{R}_{n, i}^4$ ($i\in I_n^4$) or $\mathcal{S}_n^4$.
\end{itemize}
\end{lmm}
\begin{prf}
Let $\{y, z\}$ be the kernel of $\mathcal{H}(\mathcal{B})$. First, we prove that the points of $X\backslash\{y, z\}$ can be ordered as $x_1, x_2,\ldots, x_{n - 2}$ such that one of the following cases hold:
\begin{itemize}
\item\textsc{Case} 1: $(y\ x_1\ x_2\ \ldots\ x_{n - 2})_\mathcal{B}$ and $(z\ x_1\ x_2\ \ldots\ x_{n - 2})_\mathcal{B}$;%\label{Egen1}
\item\textsc{Case} 2: there exists an index $1\leq i < n - 2$ such that $(x_1\ \allowbreak x_2\ \allowbreak\ldots\ \allowbreak x_i\ \allowbreak y\ \allowbreak x_{i + 1}\ \allowbreak\ldots\ \allowbreak x_{n - 2})_\mathcal{B}$ and $(x_1\ \allowbreak x_2\ \allowbreak\ldots\ \allowbreak x_i\ \allowbreak z\ \allowbreak x_{i + 1}\ \allowbreak\ldots\ \allowbreak x_{n - 2})_\mathcal{B}$;%\label{Egen2}
\item\textsc{Case} 3: $(y\ x_1\ x_2\ \ldots\ x_{n - 2})_\mathcal{B}$ and $(x_1\ x_2\ \ldots\ x_{n - 2}\ z)_\mathcal{B}$.%\label{Egen3}
\end{itemize}

The substructures $\mathcal{B} - y$ and $\mathcal{B} - z$ are ordered by Proposition \ref{Plin}: they are clearly linear, and neither one is isomorphic to $\mathcal{C}_4$ because otherwise $n = 5$ and $\mathcal{B}$ is irregular. It follows that $\mathcal{B} - y - z$ is ordered too, hence, with an appropriate ordering of the points of $X\backslash\{y, z\}$, $$\mathcal{B} - y - z = [x_1, x_2,\ldots, x_{n - 2}].$$
Points $x_1, x_2,\ldots, x_{n - 2}$ must be in the same order in both $\mathcal{B} - y$ and $\mathcal{B} - z$. Let $1\leq j, k\leq n - 1$ be the positions of $y$ and $z$ in $\mathcal{B} - z$ and $\mathcal{B} - y$, respectively.

We show that if $j\neq k$, then \begin{equation}\label{EQ6}
j = 1\emph{ and }k = n - 1,\emph{ or }j = n - 1\emph{ and }k = 1.
\end{equation} If (\ref{EQ6}) is false, then we can assume without loss of generality that $1 < j < n - 1$. Now, $(x_{j - 1}\ y\ x_j)_\mathcal{B}$ is true but $(x_{j - 1}\ z\ x_j)_\mathcal{B}$ is false, hence, we can apply Observation \ref{Oblowup} to $\mathcal{B} - y$ to obtain that $\mathcal{B}$ is ordered in contradiction with its nonlinearity.
Hence, either $j = k$ or (\ref{EQ6}) holds. Reversing the ordering of points $x_1, x_2,\ldots, x_{n - 2}$ if necessary, $\mathcal{B}$ satisfies one of the cases listed above. We complete the proof by showing that
\begin{itemize}
%[label={\Alph*.}, ref={\Alph*}]
%\renewcommand{\theenumi}{\Alph{enumi}}
\item if Case 1 holds, then $2\leq c\leq 3$ and $\mathcal{B}\simeq\mathcal{Q}_n^c$;
\item if Case 2 holds, then $2\leq c\leq 4$ and $\mathcal{B}\simeq\mathcal{R}_{n, i}^c$ for some $i\in I_n^c$;
\item if Case 3 holds, then $2\leq c\leq 4$ and $\mathcal{B}\simeq\mathcal{S}_n^c$.
\end{itemize}

It is clear that $c\geq 2$ in all three cases since $\mathcal{H}(\mathcal{B})$ is a tight star.

\begin{case}{1}
We show that if $(y\ z\ x_j)_\mathcal{B}$ holds for some $1\leq j\leq n - 2$, then $\{y, z, x_k\}$ is a collinear triple for all $1\leq k\leq n - 2$, $k\neq j$, which would violate the nonlinearity of $\mathcal{B}$.
There are two possibilities:
\begin{itemize}
\item if $k < j$, then $(y\ z\ x_j)_\mathcal{B}$ and $(z\ x_k\ x_j)_\mathcal{B}$ implies $(y\ z\ x_k)_\mathcal{B}$ by f.r.p.;
\item if $j < k$, then $(y\ z\ x_j)_\mathcal{B}$ and $(y\ x_j\ x_k)_\mathcal{B}$ implies $(y\ z\ x_k)_\mathcal{B}$ by f.r.p.
\end{itemize}
Similarly, $(z\ y\ x_j)_\mathcal{B}$ cannot hold for any $1\leq j\leq n - 2$.

Next, we show that if $(y\ x_j\ z)_\mathcal{B}$ holds, then $j = 1$.
Suppose that there exists an integer $k$ such that $1\leq k < j$. Now, $(y\ x_j\ z)_\mathcal{B}$ and $(y\ x_k\ x_j)_\mathcal{B}$ implies $(x_k\ x_j\ z)_\mathcal{B}$ by f.r.p. that would contradict $(z\ x_k\ x_j)_\mathcal{B}$ from the case's assumptions. So, there is no such $k$ and consequently, $j = 1$.

This also means that only one extra betweennesses, $(y\ x_1\ z)_\mathcal{B}$, can hold in $\mathcal{B}$. If it does hold, then $c = 3$ and $\mathcal{B}\simeq\mathcal{Q}_n^3$. If it does not hold, then $c = 2$ and $\mathcal{B}\simeq\mathcal{Q}_n^2$. Note that $\mathcal{Q}_n^c$ is defined because $n\geq c + 1$ by the non-linearity of $\mathcal{B}$.
%If it does not hold, then $c = 2$ and $\mathcal{B}\simeq\mathcal{Q}_n^2$. If $(y\ x_1\ z)_\mathcal{B}$ is true, then $c = 3$ and $\mathcal{B}\simeq\mathcal{Q}_n^3$.
\end{case}

\begin{case}{2}
Similarly to the previous case, we show that if $(x_j\ y\ z)_\mathcal{B}$ holds for some $1\leq j\leq n - 2$, then $\{y, z, x_k\}$ is a collinear triple for all $1\leq k\leq n - 2$, $k\neq j$, which would violate the nonlinearity of $\mathcal{B}$. We can assume by symmetry that $j\leq i$. There are three possibilities:
\begin{itemize}
\item if $k < j$, then $(x_j\ y\ z)_\mathcal{B}$ and $(x_k\ x_j\ z)_\mathcal{B}$ implies $(x_k\ y\ z)_\mathcal{B}$ by f.r.p.;
\item if $j < k\leq i$, then $(x_j\ y\ z)_\mathcal{B}$ and $(x_j\ x_k\ y)_\mathcal{B}$ implies $(x_k\ y\ z)_\mathcal{B}$ by f.r.p.;
\item if $i < k$, then $(x_j\ y\ z)_\mathcal{B}$ and $(x_j\ z\ x_k)_\mathcal{B}$ implies $(y\ z\ x_k)_\mathcal{B}$ by f.r.p.
\end{itemize}
Similarly, $(y\ z\ x_j)_\mathcal{B}$ cannot hold for any $1\leq j\leq n - 2$.

Next, we show that if $(y\ x_j\ z)_\mathcal{B}$ holds, then $j = i$ or $j = i + 1$.
Assume again that $j\leq i$. If there exists an integer $k$ such that $j < k\leq i$, then $(y\ x_j\ z)_\mathcal{B}$ and $(x_j\ x_k\ y)_\mathcal{B}$ implies $(x_k\ x_j\ z)_\mathcal{B}$ by f.r.p., in contradiction with $(x_j\ x_k\ z)_\mathcal{B}$. So there is no such $k$, thus, $j = i$. Similarly, $j = i + 1$ follows from $j\geq i$ by symmetry.

This means that only two extra betweennesses, $(y\ x_i\ z)_\mathcal{B}$ and $(y\ x_{i + 1}\ z)_\mathcal{B}$, can hold in $\mathcal{B}$.
If both of them hold, then $c = 4$ and $\mathcal{B}\simeq\mathcal{R}_{n, i}^4$.
If only one of them holds, then $c = 3$ and $\mathcal{B}\simeq\mathcal{R}_{n, i}^3$ or $\mathcal{B}\simeq\mathcal{R}_{n, n - 2 - i}^3$. 
Finally, if there are no extra betweennesses, then $c = 2$ and $\mathcal{B}\simeq\mathcal{R}_{n, i}^2$. Note that $\mathcal{R}_{n, i}^c$ is defined because $n\geq c + 1$ by the non-linearity of $\mathcal{B}$ and also $n\geq 4$ by the assumption of Case 2.
%If both are false, then $c = 2$ and $\mathcal{B}\simeq\mathcal{R}_{n, i}^2$. If only one of them holds, then $c = 3$ and $\mathcal{B}\simeq\mathcal{R}_{n, i}^3$ or $\mathcal{B}\simeq\mathcal{R}_{n, n - 2 - i}^3$. Finally, if both are true, then $c = 4$ and $\mathcal{B}\simeq\mathcal{R}_{n, i}^4$.
\end{case}

\begin{case}{3}
First, we show that if $(y\ x_j\ z)_\mathcal{B}$ holds for some $1\leq j\leq n - 2$, then $\{y, z, x_k\}$ is a collinear triple for all $1\leq k\leq n - 2$, $k\neq j$, which would violate the nonlinearity of $\mathcal{B}$.
We can assume by symmetry that $k < j$. Now, relations $(y\ x_j\ z)_\mathcal{B}$ and $(y\ x_k\ x_j)_\mathcal{B}$ imply $(y\ x_k\ z)_\mathcal{B}$ by f.r.p., which is exactly what we wanted to show.

To finish this case, we show that $(x_j\ y\ z)_\mathcal{B}$ implies $j = 1$ (and similarly, $(y\ z\ x_j)_\mathcal{B}$ implies $j = n - 2$).
Suppose that there exists an integer $k$ such that $1\leq k < j$.
Now, $(x_j\ x_k\ z)_\mathcal{B}$ follows from $(x_j\ y\ z)_\mathcal{B}$ and $(y\ x_k\ x_j)_\mathcal{B}$ by f.r.p., in contradiction with $(x_k\ x_j\ z)_\mathcal{B}$. This means that there is no such $k$, hence, $j = 1$.

So, only two extra betweennesses, $(x_1\ y\ z)_\mathcal{B}$ and $(y\ z\ x_{n - 2})_\mathcal{B}$, can hold in $\mathcal{B}$.
If both of them holds, then $c = 4$ and $\mathcal{B}\simeq\mathcal{S}_n^4$. If exactly one of them holds, then $c = 3$ and $\mathcal{B}\simeq\mathcal{S}_n^3$. Finally, if there are no extra betweennesses, then $c = 2$ and $\mathcal{B}\simeq\mathcal{S}_n^2$. Note again that $\mathcal{S}_n^c$ is defined because $n\geq c + 1$ by the non-linearity of $\mathcal{B}$.
\end{case}
$\square$
%If both are false, then $c = 2$ and $\mathcal{B}\simeq\mathcal{S}_n^2$. If exactly one of them holds, then $c = 3$ and $\mathcal{B}\simeq\mathcal{S}_n^3$. Finally, if both are true, then $c = 4$ and $\mathcal{B}\simeq\mathcal{S}_n^4$.
%\end{case}
%$\square$
\end{prf}

%%% Include or not???
%\begin{crl}\label{CRphi}
%Let $\Phi$ be a hereditary property of betweenness structures and suppose
%that there exists an integer $c\geq 2$ such that $\tau_\Phi(c + 2, 0)\geq 2$. Further, let $n > 2c - 1$ be an integer, $\mathcal{B}\in B_\Phi(n, n - c)$ be a betweenness structure and suppose that either $n\geq 6$ or $\mathcal{B}$ is regular. Then $c\leq 4$, and
%\begin{itemize}
%\item if $c = 2$ or $c = 3$, then $\mathcal{B}$ is isomorphic to either $\mathcal{Q}_n^c$, $\mathcal{R}_{n, i}^c$ ($i\in I_n^c$) or $\mathcal{S}_n^c$.
%\item if $c = 4$, then $\mathcal{B}$ is isomorphic to either $\mathcal{R}_{n, i}^4$ ($i\in I_n^4$) or $\mathcal{S}_n^4$;
%\end{itemize}
%\end{crl}
%
%\begin{prf}
%Chain together Lemmas \ref{Lgen1}, \ref{Lgen2a} and \ref{Lgen3}.
%$\square$
%\end{prf}
%We note that Corollary \ref{CRphi} does not rule out the existence of betweenness structures $\mathcal{B}\in B_\Phi(n, n - c)$ for $c < 2$. It will be a consequence of Lemma \ref{Lgen2b}.

%--------------------------------------------------------------------------------------------------

\section{Proof of the Main Results}\label{Sprf}

\begin{prfof}{Theorem \ref{Tqus1}}

\ref{Equs11}. It is obvious that $R_{n, i}^4$ ($n\geq 5$, $i\in I_n^4$), $S_n^4$ ($n\geq 5$) and all the graphs in Figure \ref{Fsmallgr} induce quasilinear betweenness structures.
Let $\Phi$ be the \emph{trivial property}, i.e. $\Phi$ is true for all betweenness structures. Lemma \ref{Lcyc} shows that the condition of Lemma \ref{Lgen1} holds for $\Phi$ and $c = 4$, thus, we obtain $$\tau(n, 0)\geq n - 4$$ for all $n\geq 3$. For $n\geq 5$, $\mathcal{S}_n^4$ proves the sharpness of this bound, and it is also easy to see that $\tau(3, 0) =\tau(4, 0) = 1$.

\ref{Equs12}.
%Quasilinear betweenness structures of order $n\leq 6$ are characterized by Lemma \ref{Lsmallgr}. %, the proof of which can be found in the Appendix (Section \ref{Sapp}).
The lemma below characterizes quasilinear betweenness structures of order $n\leq 6$. We skip the proof of this statement as it is a long but straightforward case analysis. See \cite{gnszabo2018qusapxprep} for a fully detailed proof.
%The proof is somewhat long and technical, so we present it in the Appendix at the end of the paper.
\begin{lmm}\label{Lsmallgr}
Up to isomorphism, the quasilinear betweenness structures of order at most $7$ are the following:
%A betweenness structure of order $1\leq n\leq 7$ is quasilinear if and only if it is isomorphic to one of the following betweenness structures:
\begin{itemize}
\item $\mathcal{R}_{n, i}^4$ for $5\leq n\leq 7$, $i\in I_n^4$;
\item $\mathcal{S}_n^4$ for $5\leq n\leq 7$;
\item $\mathcal{B}(G)$ where $G$ is one of the graphs in Figure \ref{Fsmallgr}.
\end{itemize}
\end{lmm}

Next, let $\mathcal{B}\in B(n, n - 4)$ be a quasilinear betweenness structure of order $n\geq 7$.
The following claim is an easy consequence of Lemma \ref{Lsmallgr}.

\begin{crl}\label{CRn7}
Let $\mathcal{A}\in B(7, 3)$ be a betweenness structure. Then $\mathcal{H}(\mathcal{A})$ is a tight star.
\end{crl}

Now, we can apply either Claim \ref{CRn7} ($n = 7$) or Lemma \ref{Lgen2a} with $c = 4$ ($n\geq 8$) to obtain that $\mathcal{H}(\mathcal{B})$ is a tight star. Finally, Lemma \ref{Lgen3} for $c = 4$ proves that $\mathcal{B}\simeq\mathcal{R}_{n, i}^4$ for some $i\in I_n^4$ or $\mathcal{B}\simeq\mathcal{S}_n^4$.
$\square$
\end{prfof}

\begin{prfof}{Theorem \ref{Tqus2}}
First, notice that every betweenness structure $\mathcal{B}\in B(5, 1)$ contains a cyclic line by Theorem \ref{Tqus1}, hence, $\tau_\Phi(5, 0)\geq 2$. We can now apply Lemma \ref{Lgen1} to $\Phi$ and $c = 2$ to obtain that $$\tau_\Phi(n, 0) \geq n - 3$$ for all $n\geq 3$.

$\tau_\Phi(3, 0) = 1$ is obvious. For $n\geq 4$, $\mathcal{Q}_n^3$, $\mathcal{R}_{n, i}^3$ $(i\in I_n^3)$ and $\mathcal{S}_n^3$ are nonlinear regular betweenness structures of co-size $n - 3$, hence, they are of minimum co-size, too, and for all $n\geq 4$, $$\tau_\Phi(n, 0) = n - 3.$$

Next, we characterize the extremal cases. Let $\mathcal{B}$ be a nonlinear regular betweenness structures of minimum co-size. Case $n = 3$ is trivial, so we can assume that $n\geq 4$ and $|\mathcal{B}|_\Delta = n - 3$.
It is enough to show that $\mathcal{H}(\mathcal{B})$ is a tight star. Then $\mathcal{B}$ is characterized by Lemma \ref{Lgen3} applied with $c = 3$.

If $n\geq 6$, then Lemma \ref{Lgen2a} applied to $\Phi$ and $c = 3$ proves that $\mathcal{H}(\mathcal{B})$ is a tight star.
If $n = 4$, then $|\mathcal{B}|_\Delta = 1$ and $\mathcal{H}(\mathcal{B})$ is a tight star again. Finally, Claim \ref{Cn5} below completes the proof with case $n = 5$.

\begin{clm}\label{Cn5}
Let $\mathcal{B}\in B(5, 2)$ be a regular betweenness structure. Then $\mathcal{H}(\mathcal{B})$ is a tight star.
\end{clm}
%\end{prfof}
%We prove Claim \ref{Cn5} in the Appendix attached to the end of the paper. %(Section \ref{Sapp}).
We skip the proof of Claim \ref{Cn5} here, as it is a straightforward case analysis. We refer the interested reader to \cite{gnszabo2018qusapxprep} where a complete proof can be found.
%The proof of Claim \ref{Cn5} can be found in \cite{gnszabo2018qusapxprep}.

\noindent$\square$
\end{prfof}

\begin{prfof}{Theorem \ref{Tqus3}}
First, notice that the betweenness structures induced by the graphs $H_4^1$ and $H_4^2$ are not orderable, hence, we obtain from Theorem \ref{Tqus1} that there are no orderable quasilinear betweenness structures on $4$ points, i.e. $\tau_\Phi(4, 2)\geq 2$.
%First, we show that every nonlinear orderable betweenness structure $\mathcal{B}$ of order $4$ is of co-size at least $2$.
%Assume to the contrary that $|\mathcal{B}|_\Delta = 1$, and let $[y_1, y_2, y_3, y_4]$ be an ordered extension of $\mathcal{B}$. We can suppose by symmetry that the sole triangle of $\mathcal{B}$ is either $\{y_1, y_2, y_3\}$ or $\{y_1, y_2, y_4\}$. In the former case $(y_1\ y_2\ y_4)_\mathcal{B}$ and $(y_2\ y_3\ y_4)_\mathcal{B}$, while in the latter case $(y_1\ y_2\ y_3)_\mathcal{B}$ and $(y_1\ y_3\ y_4)_\mathcal{B}$ leads to a contradiction by f.r.p.

Now, we can apply Lemma \ref{Lgen1} with $\Phi$ and $c = 2$ to obtain that for all $n\geq 3$, $$\tau_\Phi(n, 0)\geq n - 2.$$
It is obvious that $\mathcal{Q}_n^2$ ($n\geq 3$), $\mathcal{R}_{n, i}^2$ ($n\geq 4$, $i\in I_n^2$) and $\mathcal{S}_n^2$ ($n\geq 3$) are nonlinear orderable betweenness structures of co-size $n - 2$, hence, they are of minimum co-size as well, and for all $n\geq 3$, $$\tau_\Phi(n, 0) = n - 2.$$

Next, we characterize the extremal cases. Let $\mathcal{B}$ be a nonlinear orderable betweenness structure of co-size $n - 2$.
It is clear that $\mathcal{H}(\mathcal{B})$ is a tight star if $n = 3$.
The same holds for $n\geq 4$ as shown by Lemma \ref{Lgen2a} with parameter $c = 2$. Finally, since orderable betweenness structures are also regular, we can apply Lemma \ref{Lgen3} to obtain the desired characterization. $\square$
%For $n = 4$, notice that the betwenness structures induced by graphs $H_4^1$ and $H_4^2$ are not orderable, hence, we obtain from Theorem \ref{Tqus1} that there are no orderable quasilinear betweenness structures on $4$ points, i.e. $\tau_\Phi(4, 2)\geq 2$.
%As $\Phi(\mathcal{C}_4)$ is false,
%Now, it follows from Lemma \ref{Lgen1} and Lemma \ref{Lgen2a} applied with $c = 2$ that $\mathcal{H}(\mathcal{B})$ is a tight star for $n\geq 4$. Finally, we apply Lemma \ref{Lgen3} to obtain the desired characterization. $\square$
\end{prfof}

\begin{lmm}\label{L5-c}
Let $c$ be an integer and $\mathcal{B}\in B(n, n - c)$ be a betweenness structure. Then for all points $x$ of $\mathcal{B}$, $$d_\mathcal{B}(x) = n - c\emph{ or }d_\mathcal{B}(x)\leq 5 - c.$$
\end{lmm}
\begin{prf}
Suppose that $d_\mathcal{B}(x) < n - c$ and set $\mathcal{B}' =\mathcal{B} - x$. $\mathcal{B}'$ is a nonlinear betweenness structure with $n - 1$ points and $n - c - d_\mathcal{B}(x) > 0$ triangles, hence, Theorem \ref{Tqus1} yields $$n - c - d_\mathcal{B}(x)\geq\tau(n - 1, 0)\geq n - 5,$$ from which $d_\mathcal{B}(x)\leq 5 - c$ follows. $\square$
\end{prf}

\begin{lmm}\label{Lgen2b}
Let $c\leq 4$ and $n\geq N_c = 11 - c$ be integers, and let $\mathcal{B}\in B(n, n - c)$ be a betweenness structure. Then $\mathcal{H}(\mathcal{B})$ is a tight star.
\end{lmm}

\begin{prf}
We prove the lemma by descending induction on $c$.
If $c = 4$ and $n\geq N_4 = 7$, then Theorem \ref{Tqus1} shows that  $\mathcal{H}(\mathcal{B})$ is a tight star.

Next, suppose that $c < 4$ and the lemma is true for all $4\geq c' > c$. First, observe that $$n >\max\{c, 2c - 1\}$$ since $N_c > c$ and $N_c > 2c - 1$ for $c\leq 3$.
If for all points $x\in X$, $d_\mathcal{B}(x) = n - c$ or $d_\mathcal{B}(x)\leq 1$, then $\mathcal{B}$ satisfies the conditions of Observation \ref{Otstar} and we are done. 
So, assume that $x$ is a point of $\mathcal{B}$ such that \begin{equation}\label{EQ7}
1 < d_\mathcal{B}(x) < n - c.
\end{equation} First, notice that Lemma \ref{L5-c} gives \begin{equation}\label{EQ8}
d_\mathcal{B}(x)\leq 5 - c.
\end{equation}
Set $\mathcal{B}' =\mathcal{B} - x$ and $$c' = c + d_\mathcal{B}(x) - 1.$$ $\mathcal{B}'$ has $n - 1$ points and $n - c - d_\mathcal{B}(x) = (n - 1) - c'$ triangles. It is obvious from (\ref{EQ7}) that $c' < n - 1$, hence, $c'\leq 4$ by Theorem \ref{Tqus1}.
Also notice that $c < c'$ and \begin{equation*}\begin{split}
n - 1 & \geq N_c - 1\\ & \geq 10 - c\\ & \geq 10 - c - (d_\mathcal{B}(x) - 2)\\ & \geq N_{c'},
\end{split}\end{equation*}
thus, $\mathcal{H}(\mathcal{B}')$ must be a tight star by the induction hypothesis. Let $\{y, z\}$ be the kernel of $\mathcal{H}(\mathcal{B}')$. If $d_\mathcal{B}(y) = d_\mathcal{B}(z) = n - c$, then $\mathcal{H}(\mathcal{B})$ is clearly a tight star and the proof is complete. We show below that this is indeed the case. Assume to the contrary that, for example,
$$d_\mathcal{B}(z) < n - c.$$ Then, Lemma \ref{L5-c} yields \begin{equation}\label{EQ9}
d_\mathcal{B}(z)\leq 5 - c
\end{equation} on one hand, and (\ref{EQ8}) gives \begin{equation}\label{EQ10}\begin{split}
d_{\mathcal{B}'}(z) & = n - c - d_\mathcal{B}(x)\\ & \geq n - c - (5 - c)\\ & \geq n - 5,
\end{split}\end{equation} on the other hand.
However, by combining (\ref{EQ9}) and (\ref{EQ10}) we obtain \begin{equation*}\begin{split}
n - 5 & \leq d_{\mathcal{B}'}(z)\\ & \leq d_{\mathcal{B}}(z)\\ & \leq 5 - c,
\end{split}\end{equation*} which contradicts $n\geq N_c$.
$\square$
\end{prf}

\begin{prfof}{Theorem \ref{Tlinsize1}}

\ref{Elinsize11}. Let $c > 4$ and $1\leq n$ be integers and $\mathcal{B}\in B(n, n - c)$ be a betweenness structure. It is obvious that $B(n, n - c) =\emptyset$ if $n < c$, and $B(n, n - c)\neq\emptyset$ if $n = c$.
Further, if $c < n$, then $$|\mathcal{B}|_\Delta = n - c\geq n - 4$$ by Theorem \ref{Tqus1},  in contradiction with $c > 4$.
% Hence, $B(n, n - c)\neq\emptyset$ if and only if $n = c$.

\ref{Elinsize12}. Let $2\leq c\leq 4$ and $1\leq n$ be integers and $\mathcal{B}\in B(n, n - c)$ be a betweenness structure. First, we show that $$B(n, n - c)\neq\emptyset\Leftrightarrow n\geq c.$$ It is obvious that $B(n, n - c) =\emptyset$ if $n < c$, and $B(n, n - c)$ consists of the linear betweenness structures of order $n$ if $n = c$. Finally, if $n > c$, then $\mathcal{S}_n^c\in B(n, n - c)$.

Now, suppose that $n\geq N_c$. Lemma \ref{Lgen2b} guarantees that $\mathcal{H}(\mathcal{B})$ is a tight star. Further, since $5 < N_c\leq n$, we can apply Lemma \ref{Lgen3} to obtain the desired characterization.

\ref{Elinsize13}. Let $c < 2$ and $n\geq N_c$ be integers and $\mathcal{B}\in B(n, n - c)$ be a betweenness structure. We can apply Lemma \ref{Lgen2b} again to obtain that $\mathcal{H}(\mathcal{B})$ is a tight star. However, a tight star can have at most $n - 2$ edges, hence, $$n - c = |\mathcal{B}|_\Delta\leq n - 2,$$ contradicting $c < 2$.

\ref{Elinsize14}. Finally, let $n = N_c - 1$. Note that $n\geq 6$ for $c\leq 4$.
It is easy to see that for $1\leq i\leq\left\lceil\frac{n - 5}{2}\right\rceil$, $\mathcal{T}_{n, i}\in B(n, n - c)$ and $\mathcal{H}(\mathcal{T}_{n, i})$ is not a tight star (see Figure \ref{Ftni}).
$\square$
\end{prfof}

\begin{figure}
\centering
\includegraphics{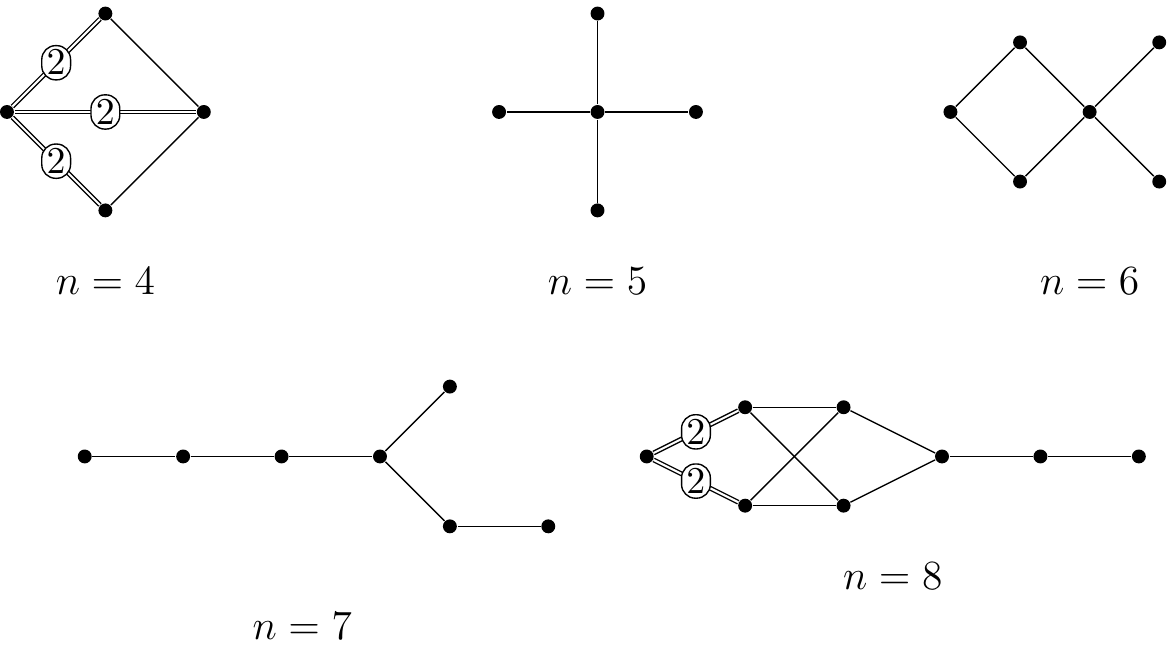}
\caption{Examples for spanner graphs of small graphic betweenness structures of co-size $n - 1$.}\label{Flinsize3}
\end{figure}

\begin{prfof}{Theorem \ref{Tlinsize3}}
Let $n\geq 9$ be an integer. Since $n - 2 < 2n - 10$ and $\mathcal{T}_{n, 1}\in B(n, 2n - 10)$, $\tau(n, n - 2)\leq 2n - 10$. Next, we show that $$\tau(n, n - 2)\geq 2n - 10.$$

Let $\mathcal{B}$ be a betweenness structure of order $n$ such that $|\mathcal{B}|_\Delta > n - 2$ and set $c = n - |\mathcal{B}|_\Delta$. Then $c < 2$ and $\mathcal{B}\in B(n, n - c)$, from which $$n\leq N_c - 1 = 10 + |\mathcal{B}|_\Delta - n$$ follows by Part \ref{Elinsize13} of Theorem \ref{Tlinsize1}, so $|\mathcal{B}|_\Delta\geq 2n - 10$.

Finally, the betweenness structures induced by the graphs in Figure \ref{Flinsize3} proves that $\tau(n, n - 2) = n - 1$ for $4\leq n\leq 8$.
$\square$
\end{prfof}

\begin{prfof}{Theorem \ref{Tlinsize2}}

\ref{Tlinsize21}. Suppose that $\mathcal{B}\in B(n, 2n - c)$ and set $c' = c - n$. Now, \begin{equation}\label{EQ11}
\mathcal{B}\in B(n, n - c').
\end{equation} On one hand, if $n < c - 4$, i.e. $c' > 4$, then $n = c'$ by Part \ref{Elinsize11} of Theorem \ref{Tlinsize1}, which gives $n = c / 2$. On the other hand, if $n > c - 2$ i.e. $c' < 2$, then $c\geq 11$ implies $n\geq N_{c'}$, and Part \ref{Elinsize13} of Theorem \ref{Tlinsize1} yields $B(n, n - c') =\emptyset$ in contradiction with (\ref{EQ11}). In summary, either $n = c / 2$ or $c - 4\leq n\leq c - 2$. In the latter case, $\mathcal{S}_n^c$ is the evidence for $B(n, 2n - c)\neq\emptyset$.

\ref{Tlinsize22}. It is easy to see that $$B(n, 2n - 10)\neq\emptyset\Leftrightarrow n\geq 5:$$ if $n = 5$, then take the ordered betweenness structure on $5$ points; otherwise take a betweenness structure that satisfies Part \ref{Elinsize14} of Theorem \ref{Tlinsize1} with $c = 10 - n\leq 4$, for example, $\mathcal{T}_{n, 1}$.

Next, suppose that $n\geq 9$. Set $c = 10 - n$ and let $\mathcal{B}\in B(n, n - c) = B(n, 2n - 10)$ be a betweenness structure on ground set $X$. Below, we prove that $\mathcal{B}\simeq\mathcal{T}_{n, i}$ for some $1\leq i\leq\left\lceil\frac{n - 5}{2}\right\rceil$.

\begin{clm}\label{Cstruct}
%Let $n\geq 9$ be an integer and $\mathcal{B}\in B(n, 2n - 10)$ be a betweenness structure. Then o

\begin{figure}[t]
\centering
\begin{subfigure}[b]{0.33\textwidth}
\centering
\includegraphics{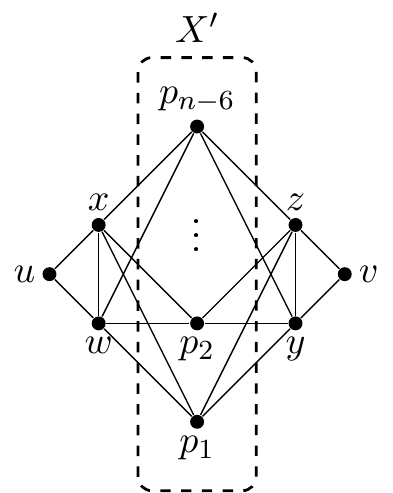}
%\begin{tikzpicture}[pont/.style={circle, fill=black, inner sep=0.5mm}]
%
%\node (W) at (1,-0.5) [pont] {};
%\node (X) at (1,0.5) [pont] {};
%\node (Y) at (3,-0.5) [pont] {};
%\node (Z) at (3,0.5) [pont] {};
%\node (U) at (0.5,0) [pont] {};
%\node (V) at (3.5,0) [pont] {};
%\node (P1) at (2,-1.5) [pont] {};
%\node (P2) at (2,-0.5) [pont] {};
%\node (PN-6) at (2,1.5) [pont] {};
%
%\draw (W) node[anchor=north] {\small$w$};
%\draw (X) node[anchor=south] {\small$x$};
%\draw (Y) node[anchor=north] {\small$y$};
%\draw (Z) node[anchor=south] {\small$z$};
%\draw (U) node[anchor=east] {\small$u$};
%\draw (V) node[anchor=west] {\small$v$};
%\draw (P1) node[anchor=north] {\small$p_1$};
%\draw (P2) node[anchor=north] {\small$p_2$};
%\draw (2,0.5) node {\small$\vdots$};
%\draw (PN-6) node[anchor=south] {\small$p_{n - 6}$};
%
%\draw (U)--(W) node {};
%\draw (U)--(X) node {};
%\draw (W)--(X) node {};
%\draw (W)--(P1) node {};
%\draw (X)--(P1) node {};
%\draw (W)--(P2) node {};
%\draw (X)--(P2) node {};
%\draw (W)--(PN-6) node {};
%\draw (X)--(PN-6) node {};
%\draw (Y)--(Z) node {};
%\draw (Y)--(V) node {};
%\draw (Z)--(V) node {};
%\draw (P1)--(Y) node {};
%\draw (P1)--(Z) node {};
%\draw (P2)--(Y) node {};
%\draw (P2)--(Z) node {};
%\draw (PN-6)--(Y) node {};
%\draw (PN-6)--(Z) node {};
%
%\draw[thick, dashed, rounded corners] (1.4,-2.2) rectangle (2.6,2.2);
%\draw (2,2.5) node {$X'$};
%\end{tikzpicture}
\caption{Case A}\label{SFstruct1}
\end{subfigure}%
\begin{subfigure}[b]{0.33\textwidth}
\centering
\includegraphics{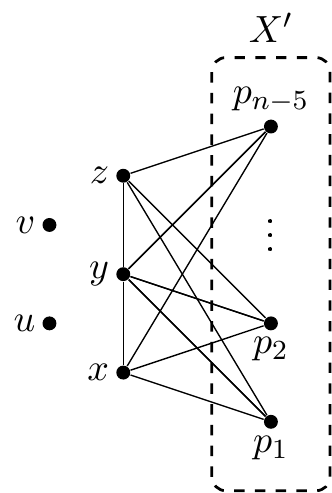}
\caption{Case B}\label{SFstruct2}
\end{subfigure}%
\begin{subfigure}[b]{0.33\textwidth}
\centering
\includegraphics{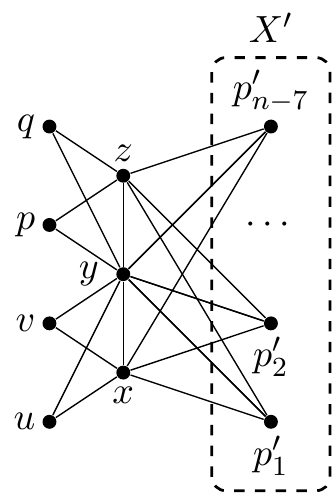}
\caption{Case C}\label{SFstruct3}
\end{subfigure}%
\caption{Triangle hypergraphs in Claim \ref{Cstruct}. The hyperedges are represented by triangles.}\label{Fstruct}
\end{figure}

One of the following three cases hold:
\begin{enumerate}[label={\Alph*.}, ref={\Alph*}]
\item there exist distinct points $u, v, w, x, y, z\in X$ such that with $X' =\allowbreak X\backslash\{u,\allowbreak v,\allowbreak w,\allowbreak x,\allowbreak y,\allowbreak z\}$, $$\Delta(\mathcal{B}) =\{\{p, w, x\} : p\in X'\}\cup\{\{p, y, z\} : p\in X'\}\cup\{\{u, w, x\}, \{v, y, z\}\};$$\label{Estruct1}
\item there exist distinct points $u, v, x, y, z\in X$ such that with $X'\allowbreak = X\backslash\{u,\allowbreak v,\allowbreak x,\allowbreak y,\allowbreak z\}$, $$\Delta(\mathcal{B}) =\{\{p, x, y\} : p\in X'\}\cup\{\{p, y, z\} : p\in X'\};$$\label{Estruct2}
\item there exist distinct points $p, q, u, v, x, y, z\in X$ such that with $X' = X\backslash\{p, q, u, v, x, y, z\}$, $$\Delta(\mathcal{B}) =\{\{p', x, y\} : p'\in X'\}\cup\{\{p', y, z\} : p'\in X'\}\cup$$ $$\{\{u, x, y\}, \{v, x, y\}, \{p, y, z\}, \{q, y, z\}\}.$$\label{Estruct3}
\end{enumerate}
\end{clm}
\begin{prf}
%Set $$c = 10 - n.$$ Then $\mathcal{B}\in B(n, n - c)$. 
Observe that \begin{equation}\label{EQ11a}
\mathcal{H}(\mathcal{B})\emph{ is not a tight star},
\end{equation} because $|\mathcal{B}|_\Delta = 2n - 10 > n - 2$, the latter being the maximum number of edges in a tight star.

%Now, we obtain from Observation \ref{Otstar} that there exists a point $x\in X$ such that $$0 < d_\mathcal{B}(x) < n - c.$$
%Additionally, we assume that $x$ is a point of maximum degree such that $d_\mathcal{B}(x) < n - c$. Because of Lemma \ref{L5-c}, \begin{equation}\label{EQ12}
%d_\mathcal{B}(x)\leq 5 - c = n - 5.
%\end{equation}

\begin{clm}\label{Cdx5-c}
Let $p\in X$ be a point such that $0 < d_\mathcal{B}(p) < n - c$ and suppose that $d_\mathcal{B}(p)$ is maximal with this property. Then
\begin{enumerate} 
\item $\mathcal{H}(\mathcal{B} - p)$ is a tight star;\label{Edx5-c1}
\item $d_{\mathcal{B}}(p) = n - 5$.\label{Edx5-c2}
\end{enumerate}
Further, if $q$ is a point in the kernel of $\mathcal{H}(\mathcal{B} - p)$ such that $d_\mathcal{B}(q) < n - c$, then
%$q\in\ker(\mathcal{H}(\mathcal{B} - p))$ such that $d_\mathcal{B}(q) < n - c$, then
\begin{enumerate}\setcounter{enumi}{2}
\item $d_\mathcal{B}(q) = n - 5$ and no triangle $T\in\Delta(\mathcal{B})$ contains both $p$ and $q$.\label{Edx5-c3}
%\item no triangle $T\in\Delta(\mathcal{B})$ contains both $p$ and $q$.\label{Edx5-c4}
\end{enumerate}
\end{clm}
\begin{prf}

\ref{Edx5-c1}. Set $\mathcal{B}' =\mathcal{B} - p$. 
First of all, observe that Lemma \ref{L5-c} implies \begin{equation}\label{EQ12}
d_\mathcal{B}(p)\leq 5 - c = n - 5.
\end{equation}
Further, notice that $$\mathcal{B}'\in B(n - 1, n - 1 - c')$$ where $$c' = c + d_{\mathcal{B}}(p) - 1.$$ We show that $\mathcal{H}(\mathcal{B}')$ is a tight star by applying Lemma \ref{Lgen2b}. Two conditions must be met:
\begin{enumerate}[label={(\roman*)}, ref={(\roman*)}]
\item $c'\leq 4$;\label{Econd1}
\item $n - 1\geq N_{c'}$.\label{Econd2}
\end{enumerate}
Condition \ref{Econd1} holds because of (\ref{EQ12}). As for condition \ref{Econd2}, \begin{equation*}\begin{split}
N_{c'} & = 12 - c - d_\mathcal{B}(p)\\ & = 2 + n - d_\mathcal{B}(p),
\end{split}\end{equation*} hence, condition \ref{Econd2} holds if and only if \begin{equation}\label{EQ13}
d_\mathcal{B}(p)\geq 3.
\end{equation}
%%% use proof of advanced version for shorter argument for $d_\mathcal{B}(x)\geq 3$???
%\begin{clm}\label{Cdpgeq3}
%$d_\mathcal{B}(p)\geq 3$.
%\end{clm}
%\begin{prf}
%$\square$
%\end{prf}

Suppose to the contrary that (\ref{EQ13}) is false, and let $n_i$ denote the number of points $q'\in X$ such that $d_\mathcal{B}(q') = i$. Then, as $p$ was of maximum degree, $n_i > 0$ only if $0\leq i\leq 2$ or $i = n - c = 2n - 10$. Further, $n_{2n - 10}\leq 1$ for otherwise $\mathcal{H}(\mathcal{B})$ would be a tight star in contradiction with (\ref{EQ11a}). Now, counting the sum of degrees in $\mathcal{H}(\mathcal{B})$ in two ways, we obtain \begin{equation*}\begin{split}
3(2n - 10) & =\sum_{q'\in X} d_\mathcal{B}(q')\\ & = \sum_{i = 0}^{2n - 10} in_i\\ & = n_1 + 2 n_2 + (2n - 10) n_{2n - 10}\\ & \leq 2(n - 1) + (2n - 10),
\end{split}\end{equation*} which is equivalent to $n\leq 9$.

Since $n\geq 9$ by the theorem's assumption, the only problematic case is $n = 9$, so suppose that this is indeed the case.
Notice that $n_{2n - 10}\neq 0,$ otherwise, we obtain by the previous argument that $3(2n - 10)\leq n_1 + 2 n_2\leq 2n$ in contradiction with $n = 9$. Thus, $n_{2n - 10} = 1$, which further implies $n_1 = 0$, otherwise \begin{equation}\begin{split}
2(2n - 10)& \leq n_1 + 2 n_2\\
& \leq 1 + 2(n - 2)
\end{split}\end{equation} in contradiction with $n = 9$ again.
% We can also show in a similar way that $n_1 = 0$,
So, we can conclude that \begin{equation}\label{EQ14}
n_2 = n - 1 = 8.
\end{equation}

Let $y$ denote the unique point for which $d_\mathcal{B}(y) = 2n - 10 = 8$, and let $G_y$ be the (not necessarily connected) link graph of $y$ in $\mathcal{H}(\mathcal{B})$, i.e. $G_y=\{X\backslash\{y\}, E_y\}$ where $E_y =\{T\backslash\{y\} : y\in T\in\Delta(\mathcal{B})\}$. Now, we obtain from (\ref{EQ14}) that
%Notice that every triangle of $\mathcal{B}$ contains $y$ and 
$$G_y\emph{ is a disjoint union of cycles}.$$
%as a consequence of (\ref{EQ14}).
Further, note the following easy consequences of Lemma \ref{Lsmallgr}.

\begin{clm}\label{Cn6}
Let $\mathcal{A}\in B(6, 2)$ be a betweenness structure with triangles $R$ and $T$. Then $|R\cap T|\neq 1$.
\end{clm}

%This immediately implies the following claim.
\begin{clm}\label{C2indep}
There is no set of points $Y\subset X\backslash\{y\}$, $|Y| = 5$ that induces two independent edges in $G_y$.
\end{clm}
%\begin{prf}
%If $G_y[Y]$ would consist of two independent edges, then the substructure $\mathcal{B}\vert_{Y\cup\{y\}}$ would contradict Claim \ref{Cn6}.
%$\square$
%\end{prf}
%Claim \ref{C2indep} follows from Claim \ref{Cn6} which, in turn, is an easy consequence of Lemma \ref{Lsmallgr}.

Now, if $G_y$ is connected, then we can clearly find an $Y\subset X\backslash\{y\}$ that contradicts Claim \ref{C2indep}. If $G_y$ is not connected but contains a triangle, then we can again construct such a $Y$: take two points, $x_1$ and $x_2$, from the triangle, two endpoints, $z_1$ and $z_2$ of an edge that is not in the triangle and a fifth point $w$ that is not adjacent to any of the previously chosen ones. Such a $w$ exists because $G_y$ has $n - 1 = 8$ vertices, at most $7$ of which are adjacent to $x_1$, $x_2$, $z_1$ or $z_2$.
% $n = 9$, the (unified) neighborhood of $x_1$ and $x_2$ is of size $1$, and $z_1$ and $z_2$ have at most two neighbors on the circle that contains them.
 Hence, we can conclude that $$G_y\emph{ is the disjoint union of two 4-cycles}$$ and so
%\begin{obs}\label{O4indep}
\begin{equation}\label{EQ14a}
\begin{split}
&\alpha(G_y) = 4,\emph{ and any independent set of size $4$ consists of}\\
&\emph{one-one opposing pair of vertices from each components of }G_y.
\end{split}
\end{equation}
%\end{obs}

Observation (\ref{EQ14a}) has two important consequences.

\begin{clm}\label{C4indep}
Let $A$ and $B$ be two independent sets of size $4$ in $G_y$. Then $|A\cap B|$ is even.
\end{clm}

\begin{clm}\label{Cyleq5}
Let $Y\subseteq X$ be a set of points such that $y\in Y$ and $\mathcal{B}\vert_Y$ is linear. Then $|Y|\leq 5$.
\end{clm}
%\begin{prf}
%$Y\backslash\{y\}$ is an independent set of $G_y$, hence, by (\ref{EQ14a}), $|Y\backslash\{y\}|\leq 4$. $\square$
%\end{prf}

Observe that $\mathcal{B} - y$ is a linear betweenness structure of order $8$, hence, it is also ordered by Proposition \ref{Plin} and we can index its points as $x_1, x_2,\ldots, x_8$ such that $$\mathcal{B} - y = [x_1, x_2,\ldots, x_8].$$

\begin{clm}\label{Cposbetw}$ $
\begin{enumerate}
%[label={(\roman*)}, ref={\roman*}]
\item $(x_i\ y\ x_j)_\mathcal{B}\Rightarrow |j - i|\geq 5$;\label{Eposbetw1}
\item $(x_i\ x_j\ y)_\mathcal{B}\Rightarrow |j - i|\leq 3$.\label{Eposbetw2}
\end{enumerate}
\end{clm}
\begin{prf}
\ref{Eposbetw1}. Suppose that $i < j$ and set $Y =\{x_1,\ldots, x_i, y, x_j,\ldots, x_8\}$. Since $(x_1\ \ldots\ x_i\ x_j\ \ldots\ x_8)_\mathcal{B}$ and $(x_i\ y\ x_j)_\mathcal{B}$ hold, we obtain from Observation \ref{Oblowup} that $\mathcal{B}\vert_Y$ is ordered. Now, $|Y|\leq 5$ by Claim \ref{Cyleq5} and so $|j - i|\geq 5$.

\ref{Eposbetw2} Set $Y =\{x_i, x_{i + 1},\ldots, x_j, y\}$. Since $(x_i\ x_{i + 1}\ \ldots\ x_j)_\mathcal{B}$ and $(x_i\ x_j\ y)_\mathcal{B}$ hold, it follows from Observation \ref{Oblowup} that $\mathcal{B}\vert_Y$ is ordered. Now, $|Y|\leq 5$ by Claim \ref{Cyleq5} and so $|j - i|\leq 3$.
$\square$
\end{prf}

We call a pair $\{x_i, x_j\}$ an \emph{$\ell$-chord} if $|j - i| =\ell$ and the triple $\{x_i, x_j, y\}$ is collinear. Our next observation follows from Claim \ref{Cposbetw}.

\begin{crl}\label{CRposbetw}
Let $\{x_i, x_j\}$ be an $\ell$-chord and suppose that $i < j$. Then one of the following cases hold:
\begin{enumerate}
\item if $\ell\geq 5$, then $(x_1\ \ldots\ x_i\ y\ x_j\ \ldots\ x_8)_\mathcal{B}$ holds and so $\allowbreak\{x_1,\allowbreak\ldots,\allowbreak x_i,\allowbreak x_j,\allowbreak\ldots,\allowbreak x_8\}$ is an independent set in $G_y$;
\item if $\ell\leq 3$, then either $(y\ x_i\ x_{i + 1}\ \ldots\ x_j)_\mathcal{B}$ or $(x_i\ x_{i + 1}\ \ldots\ x_j\ y)_\mathcal{B}$ holds and so $\{x_i, x_{i + 1},\ldots, x_j\}$ is an independent set in $G_y$.
\end{enumerate}
\end{crl}

Let $h_\ell$ denote the number of $\ell$-chords for $1\leq\ell\leq 8$. Note that \begin{equation}\label{EQ15}
h_\ell\leq 8 -\ell
\end{equation} and \begin{equation}\label{EQ16}
h_4 = 0
\end{equation} by Claim \ref{Cposbetw}. Further, observe that $$\sum_{\ell = 1}^8 h_\ell = \binom{8}{2} - |\mathcal{B}|_\Delta = 20.$$ However, we will show below the contradiction $$\sum_{\ell = 1}^8 h_\ell\leq 19.$$

\begin{clm}\label{Cherd}
If $\{x_i, x_{i + 3}\}$ is a $3$-chord, then $\{x_i, x_{i + 2}\}$ and $\{x_{i + 1}, x_{i + 3}\}$ are $2$-chords.
\end{clm}
\begin{prf}
Because of Claim \ref{Cposbetw}, either $(y\ x_i\ x_{i + 3})_\mathcal{B}$ or $(x_i\ x_{i + 3}\ y)_\mathcal{B}$ holds, so $\{x_i,\allowbreak x_{i + 1},\allowbreak x_{i + 2},\allowbreak x_{i + 3},\allowbreak y\}$ induces a linear substructure by Observation \ref{Oblowup} and hence both $\{x_i, x_{i + 2}, y\}$ and $\{x_{i + 1}, x_{i + 3}, y\}$ are collinear.
$\square$
\end{prf}

\begin{clm}\label{Cchord}$ $
\begin{enumerate}
%[label={(\roman*)}, ref={\roman*}]
\item $h_3\leq 3$ and if $h_3 = 3$, then the $3$-chords are exactly $\{x_1, x_4\}$, $\{x_3, x_6\}$ and $\{x_5, x_8\}$;\label{Echord1}
\item $h_5\leq 2$ and if $h_5 = 2$, then the $5$-chords are exactly $\{x_1, x_6\}$ and $\{x_3, x_8\}$.\label{Echord2}
\end{enumerate}
\end{clm}
\begin{prf}

\ref{Echord1}. If $\{x_i, x_{i + 3}\}$ is a $3$-chord, then $Y_1 =\{x_i, x_{i + 1}, x_{i + 2}, x_{i + 3}\}$ is an independent set in $G_y$ by Corollary \ref{CRposbetw}.
%$(y\ x_i\ x_{i + 3})_\mathcal{B}$ or $(x_i\ x_{i + 3}\ y)_\mathcal{B}$ hold by Claim \ref{Cposbetw}, thus, $Y_1 =\{x_i, x_{i + 1}, x_{i + 2}, x_{i + 3}\}$ is an independent set by Observation \ref{Oblowup}.
Similarly, if $\{x_{i + 1}, x_{i + 4}\}$ is a $3$-chord, then $Y_2 =\{x_{i + 1}, x_{i + 2}, x_{i + 3}, x_{i + 4}\}$ is an independent set in $G_y$. However, $|Y_1\cap Y_2| = 3$ in contradiction with Claim \ref{C4indep}. Hence, we get the maximum possible number of $3$-chords if we take every second one starting with $\{x_1, x_4\}$.

\ref{Echord2}. Similarly to the previous case, if $\{x_i, x_{i + 5}\}$ and $\{x_{i + 1}, x_{i + 6}\}$ are both $5$-chords, then both $Y_1 =\{x_1,\ldots, x_i, x_{i + 5},\ldots, x_8\}$ and $Y_2 =\{x_1,\allowbreak\ldots,\allowbreak x_{i + 1},\allowbreak x_{i + 6},\allowbreak\ldots,\allowbreak x_8\}$ are independent sets by Corollary \ref{CRposbetw}, contradicting Claim \ref{C4indep}.
Hence, we get the maximum possible number of $5$-chords if we take every second one starting with $\{x_1, x_6\}$.
$\square$
\end{prf}

%The following statement is the corollary of Claim \ref{Cchord}.
\begin{clm}\label{Ch5geq0}
If $h_5 > 0$, then $h_2\leq 5$ and $h_3\leq 2$.
\end{clm}
\begin{prf}
Suppose that $h_5 > 0$ and let $\{x_i, x_{i + 5}\}$ be a $5$-chord. We already know from (\ref{EQ15}) that $h_2\leq 6$. Assume to the contrary that $h_2 = 6$. We obtain from Claim \ref{Cherd} that for all $1\leq k\leq 6$, $$\{x_k, x_{k + 2}\}\emph{ is a 2-chord}.$$

First, we show that for all $1\leq k\leq 6$ \begin{equation}\label{EQ2}
%\forall\,i\leq k\leq i + 2,\,
(y\ x_k\ x_{k + 2})_\mathcal{B}\Rightarrow (y\ x_{k + 1}\ x_{k + 3})_\mathcal{B}.
\end{equation}
Suppose to the contrary that $(y\ x_k\ x_{k + 2})_\mathcal{B}$ holds but $(y\ x_{k + 1}\ x_{k + 3})_\mathcal{B}$ is false. Then, because of Claim \ref{Cposbetw} and the fact that $\{x_{k + 1}, x_{k + 3}\}$ is a $2$-chord, $(x_{k + 1}\ x_{k + 3}\ y)_\mathcal{B}$ must be true. From this and $(x_{k + 1}\ x_{k + 2}\ x_{k + 3})_\mathcal{B}$, $(x_{k + 1}\ x_{k + 2}\ y)_\mathcal{B}$ follows by f.r.p. On the other hand, however, $(y\ x_k\ x_{k + 2})_\mathcal{B}$ and $(x_k\ x_{k + 1}\ x_{k + 2})_\mathcal{B}$ yield $(y\ x_{k + 1}\ x_{k + 2})_\mathcal{B}$, a contradiction.

Next, notice that \begin{equation}\label{EQ3}
(y\ x_i\ x_{i + 2})_\mathcal{B}
\end{equation} holds.
If not, then Claim \ref{Cposbetw} and the fact that $\{x_i, x_{i + 2}\}$ is a $2$-chord imply $(x_i\ x_{i + 2}\ y)_\mathcal{B}$. Since $\{x_i, x_{i + 5}\}$ is a $5$-chord, $(x_1\ \ldots\ x_i\ y\ x_{i + 5}\ \ldots\ x_8)_\mathcal{B}$ is true by Corollary \ref{CRposbetw}. Now, Observation \ref{Oblowup} implies $(x_1\ \ldots\ x_i\ x_{i + 2}\ y\ x_{i + 5}\ \ldots\ x_8)_\mathcal{B}$, giving an independent set $\{x_1,\ldots, x_i, x_{i + 2}, x_{i + 5},\ldots, x_8\}$ of size $5$ in $G_y$ in contradiction with (\ref{EQ14a}).
% As $\{x_i, x_{i + 2}\}$ was a $2$-chord, $(y\ x_i\ x_{i + 2})_\mathcal{B}$ must be true.

Similarly to (\ref{EQ3}), \begin{equation*}
(x_{i + 3}\ x_{i + 5}\ y)_\mathcal{B}
\end{equation*} holds, which leads to a contradiction as $(y\ x_{i + 3}\ x_{i + 5})$ follow from (\ref{EQ3}) by repeated application of (\ref{EQ2}).

For the second part of the claim, assume to the contrary that $h_3 > 2$. Then, because of Claim \ref{Cchord}, $h_3 = 3$ and the $3$-chords are exactly $\{x_1, x_4\}$, $\{x_3, x_6\}$ and $\{x_5, x_8\}$. Further, we obtain from Claim \ref{Cherd} that $h_2 = 6$, contradicting what we have just proved above.
$\square$
\end{prf}

Now, we can complete the proof of Part \ref{Edx5-c1} as follows.
\begin{itemize}
\item If $h_5 = 0$, then (\ref{EQ15}), (\ref{EQ16}) and Claim \ref{Cchord} yield
\begin{equation*}\begin{split}
\sum_{\ell = 1}^8 h_\ell & = h_1 + h_2 + h_3 + h_6 + h_7\\ & \leq 7 + 6 + 3 + 2 + 1 = 19;
\end{split}\end{equation*}
\item if $h_5 > 0$, then (\ref{EQ15}), (\ref{EQ16}), Claim \ref{Cchord} and Claim \ref{Ch5geq0} yield
\begin{equation*}\begin{split}
\sum_{\ell = 1}^8 h_\ell & = h_1 + h_2 + h_3 + h_5 + h_6 + h_7\\ & \leq 7 + 5 + 2 + 2 + 2 + 1 = 19.
\end{split}\end{equation*}
\end{itemize}

\ref{Edx5-c2}. Next, we prove that $d_\mathcal{B}(p) = n - 5$. Since we have just proved that $\mathcal{H}(\mathcal{B} - p)$ is a tight star, there are two points $y$ and $z$ different from $p$ of degree \begin{equation}\label{EQ17}
d_{\mathcal{B} - p}(y) = d_{\mathcal{B} - p}(z) = |\mathcal{B} - p|_\Delta = n - c - d_{\mathcal{B}}(p).
\end{equation} However, since $\mathcal{H}(\mathcal{B})$ is not a tight star, we can assume without loss of generality that $d_\mathcal{B}(z) < n - c$, thus, because of the maximality of $d_{\mathcal{B}}(p)$, $$d_{\mathcal{B}}(p)\geq d_{\mathcal{B}}(z)\geq d_{\mathcal{B} - p}(z),$$ from which $$d_{\mathcal{B}}(p)\geq (n - c) / 2 = n - 5$$ follows by (\ref{EQ17}). As we have already established $d_{\mathcal{B}}(p)\leq n - 5$ in (\ref{EQ12}), the proof is complete.

\ref{Edx5-c3}.
Let $q\in\ker(\mathcal{H}(\mathcal{B} - p))$ be a point such that $d_\mathcal{B}(q) < n - c$. Now, \begin{equation*}\begin{split}
d_{\mathcal{B} - p}(q) & = |\mathcal{B} - p|_\Delta\\
& = n - c - d_\mathcal{B}(p)\\
& = n - 5\\
& = d_\mathcal{B}(p),
\end{split}\end{equation*} hence, we obtain that $d_{\mathcal{B} - p}(q) = d_\mathcal{B}(q) = d_\mathcal{B}(p) = n - 5$.
%\begin{equation*}\begin{split}
%d_\mathcal{B}(p) & = d_{\mathcal{B} - p}(q)\\
%&\leq d_\mathcal{B}(q)\\
%&\leq d_\mathcal{B}(p)
%\end{split}\end{equation*} and so we have equality everywhere.
%In particular, $d_\mathcal{B}(q) = d_\mathcal{B}(p) = n - 5$.
%Further, $d_{\mathcal{B} - p}(q) = d_\mathcal{B}(q)$, which means
This also shows that any triangle that contains $p$ avoids $q$.
$\square$
\end{prf}

Since $\mathcal{B}\in B(n, n - c)$ and $n > 2c - 1$, we obtain from (\ref{EQ11a}) and Observation \ref{Otstar} that there exists a point $x\in X$ such that $$1 < d_\mathcal{B}(x) < n - c.$$ We can suppose that $x$ is such a point of maximum degree. Now, we can apply Claim \ref{Cdx5-c} to obtain that $d_\mathcal{B}(x) = n - 5$ and $\mathcal{H}(\mathcal{B} - x)$ is a tight star. Let $\{y, z\}$ be the kernel of $\mathcal{H}(\mathcal{B} - x)$, and suppose that $$d_\mathcal{B}(y)\geq d_{\mathcal{B}}(z).$$ we close the proof of Claim \ref{Cstruct} by considering the following two cases.

\begin{itemize}
\item \textsc{Case} 1: $d_\mathcal{B}(y) < n - c$;
\item \textsc{Case} 2: $d_\mathcal{B}(y) = n - c$.
\end{itemize}

\begin{clm}\label{Ccases}$ $
\begin{enumerate}
%[label={(\roman*)}, ref={\roman*}]
\item If $d_\mathcal{B}(y) < n - c$, then
there exists a point $w\in X\backslash\{x, y, z\}$ such that $d_\mathcal{B}(w) = d_\mathcal{B}(x) = d_\mathcal{B}(y) = d_\mathcal{B}(z) = n - 5$ and every triangle of $\mathcal{B}$ contains exactly one pair of points from the set $\{w, x, y, z\}$ that is either $\{w, x\}$ or $\{y, z\}$;\label{Ecases1}
\item if $d_\mathcal{B}(y) = n - c$, then
$d_\mathcal{B}(x) = d_\mathcal{B}(z) = n - 5$ and every triangle of $\mathcal{B}$ contains exactly one of the pairs $\{x, y\}$ and $\{y, z\}$.\label{Ecases2}
\end{enumerate}
\end{clm}
\begin{prf}

\ref{Ecases1}. We already know that $d_\mathcal{B}(x) = n - 5$. Let $T$ be an arbitrary triangle of $\mathcal{B}$.
Because $d_\mathcal{B}(z)\leq d_\mathcal{B}(y) < n - c$, we obtain from Part \ref{Edx5-c3} of Claim \ref{Cdx5-c} that $d_\mathcal{B}(y) = d_\mathcal{B}(z) = n - 5$ and \begin{equation}\label{EQ18b}
\{x, y\}, \{x, z\}\not\subset T.
\end{equation}

%%++++
%Because of the maximality of $d_\mathcal{B}(x)$, $d_\mathcal{B}(y) < n - c$ implies
%$$d_\mathcal{B}(x)\geq d_\mathcal{B}(y)\geq d_{\mathcal{B} - x}(y) = |\mathcal{B} - x|_\Delta = n - c - d_\mathcal{B}(x) = n - 5 = d_\mathcal{B}(x).$$ Hence, $d_\mathcal{B}(y) = n - 5$. Further, $d_\mathcal{B}(y)  = d_{\mathcal{B} - x}(y)$, thus, \begin{equation}\label{EQ18b}
%\{x, y\}\not\subset T.
%\end{equation}
%%----
%%++++
%Similarly, as $d_\mathcal{B}(z)\leq d_\mathcal{B}(y) < n - c$, we obtain that \begin{equation}\label{EQ18c}
%\{x, z\}\not\subset T.
%\end{equation}
%%----
Now, since every triangle that avoids $x$ contains both $y$ and $z$, we obtain from (\ref{EQ18b}) that \begin{equation}\label{EQ18d}
y\in T\Leftrightarrow z\in T.
\end{equation}

Since $d_\mathcal{B}(y) = d_\mathcal{B}(x)$, we obtain from Part \ref{Edx5-c1} of Claim \ref{Cdx5-c} that $\mathcal{B} - y$ is a tight star. Further, (\ref{EQ18b}) yields $$d_{\mathcal{B} - y}(x) = d_\mathcal{B}(x) = n - 5 = |\mathcal{B}|_\Delta - d_\mathcal{B}(y) = |\mathcal{B} - y|_\Delta,$$ hence, $x\in\ker(\mathcal{H}(\mathcal{B} - y))$. Let $w$ be the other point of that kernel. Point $w$ is different from $x$ and $y$ by definition. If $w = z$ would be true, then all $n - 5 > 1$ triangles of $\mathcal{B} - y$ would contain $x$ and $z$, which contradicts (\ref{EQ18b}). Therefore, $w, x, y$ and $z$ are distinct points.

Next, observe that $d_\mathcal{B}(w) < n - c$. Otherwise, all triangles would contain $w$ and so $w$ would be in the kernel of $\mathcal{B} - x$, which is impossible since $w\neq y$ and $w\neq z$. Now, we can apply Part \ref{Edx5-c3} of Claim \ref{Cdx5-c} with $p = y$ and $q = w$ to obtain that $d_\mathcal{B}(w) = n - 5$ and $\{w, y\}\not\subset T.$
%%++++
%Therefore, $$d_\mathcal{B}(x)\geq d_\mathcal{B}(w)\geq d_{\mathcal{B} - y}(w) = |\mathcal{B} - y|_\Delta = n - c - d_\mathcal{B}(y) = n - 5 = d_\mathcal{B}(x)$$ and hence $d_\mathcal{B}(w) = n - 5$.
%Further, $d_\mathcal{B}(w) = d_{\mathcal{B} - y}(w)$, hence, \begin{equation}\label{EQ18e}
%\{w, y\}\not\subset T.
%\end{equation}
%%----
Since every triangle that avoids $y$ contains both $w$ and $x$,
\begin{equation}\label{EQ18f}
w\in T\Leftrightarrow x\in T
\end{equation} follows.
%from (\ref{EQ18b}) and (\ref{EQ18e}).
%Finally, we show that triangle $T$ contains exactly one of the pairs $\{x, w\}$ or $\{y, z\}$, and it does not contain any of the pairs $\{x, y\}, \{x, z\}, \{w, y\}$ and $\{w, z\}$.
Now, triangle $T$ either contains $x$, or it is in $\Delta(\mathcal{B} - x)$ and hence contains $y$, so we obtain from (\ref{EQ18d}) and (\ref{EQ18f}) that $T$ contains exactly one of the pairs $\{w, x\}$ or $\{y, z\}$, and it does not contain any other pairs from $\{w, x, y, z\}$.

%\ref{Ecases1}. As $d_{\mathcal{B}'}(y) = d_{\mathcal{B}'}(z) = n - c - d_{\mathcal{B}}(x) = d_{\mathcal{B}}(x)$ and $d_{\mathcal{B}}(x)$ was chosen to be maximal, $$d_{\mathcal{B}}(x) = d_{\mathcal{B}}(y) = d_{\mathcal{B}}(z) = 5 - c$$ and \begin{equation}\label{EQ18}
%\emph{any triangle that contains }x\emph{ avoids }y\emph{ and }z.
%\end{equation}
%
%Now, we can swap the roles of $x$ and $y$ to obtain that $\mathcal{H}(\mathcal{B} - y)$ is a (non-empty) tight star and $x$ is in its kernel. Let $w$ denote the other point of said kernel. Notice that $w\neq z$ because of (\ref{EQ18}). Further, $d_{\mathcal{B}}(w)\neq n - c$, otherwise, since $d_{\mathcal{B} - y}(w) = 5 - c$, (\ref{EQ18}) implies $d_{\mathcal{B}'}(w) = d_{\mathcal{B}}(w) -  d_{\mathcal{B} - y}(w) = 5 - c$ as well, which would contradict $w\neq y, z$. Now, following the above argument, we obtain that $$d_{\mathcal{B}}(w) = 5 - c$$ and any triangle that contains $y$ avoids $x$ and $w$, hence, any triangle that contains either $x$ or $w$ is an edge of $\mathcal{H}(\mathcal{B} - y)$. In summary, $$d_{\mathcal{B}}(x) = d_{\mathcal{B}}(y) = d_{\mathcal{B}}(z) = d_{\mathcal{B}}(w) = 5 - c$$ and every triangle $T$ of $\mathcal{B}$ contains exactly one of the pairs $\{x, w\}$ and $\{y, z\}$.

\ref{Ecases2}. Notice that $d_{\mathcal{B}}(z) < n - c$ because $y$ is already contained in all triangles but $\mathcal{H}(\mathcal{B})$ is not a tight star by (\ref{EQ11a}). Now, we obtain from Part \ref{Edx5-c1} of Claim \ref{Cdx5-c} that $d_\mathcal{B}(z) = n - 5$ and no triangle of $\mathcal{B}$ contains both $x$ and $z$.
%%++++
%$$d_\mathcal{B}(x)\geq d_\mathcal{B}(z)\geq d_{\mathcal{B} - x}(z) = n - c - d_\mathcal{B}(x) = n - 5,$$ from which we obtain that $d_{\mathcal{B}}(z) = d_{\mathcal{B} - x}(z) = n - 5$
%and so no triangle contains $x$ and $z$ simultaneously.
%%----
Now, as $y$ is of full degree, it can be easily seen that every triangle of $\mathcal{B}$ contains exactly one of the pairs $\{x, y\}$ and $\{y, z\}$.
$\square$
\end{prf}

%Now, we close the proof of Claim \ref{Cstruct} by analyzing the two cases in Claim \ref{Ccases}.
Now, we analyze Cases 1-2 with the help of Claim \ref{Ccases}.

\begin{case}{1}
Claim \ref{Ccases} guarantees that the set of triangles of $\mathcal{B}$ can be divided into $wx$-triangles and $yz$-triangles, depending on whether they contain the pair $\{w, x\}$ or $\{y, z\}$.

Set $Y = X\backslash\{w, x, y, z\}$ and let $n_i$ denote the number of points of degree $i$ in $\mathcal{B}$.
Then, for every point $p\in X$, $d_\mathcal{B}(p)\leq 2\Leftrightarrow p\in Y$, hence, $n_0 + n_1 + n_2 = |Y| = n - 4$. On the other hand,
\begin{equation*}\begin{split}
n_1 + 2n_2 & =\sum_{p\in Y} d_\mathcal{B}(p)\\
& = |\mathcal{B}|_\Delta\\
& = 2n - 10,
\end{split}\end{equation*} so one of the following cases hold:
%Note that $n_i > 0\Rightarrow i\in\{0, 1, 2, n - 5\}$ and $n - 5 > 2$, so $n_{n - 5} = 4$ and $n_0 + n_1 + n_2 = n - 4$ (see Claim \ref{Ccases}). Now, for every point $p\in Y$, $d_\mathcal{B}(p)\leq 2$, hence, $$\sum_{p\in Y} d_\mathcal{B}(p)\leq 2(n - 4).$$ On the other hand, Claim \ref{Ccases} also yields $$\sum_{p\in Y} d_\mathcal{B}(p) = |\mathcal{B}|_\Delta = 2n - 10,$$ so one of the following cases hold:
\begin{itemize}
%[label={(\alph*)}, ref={\alph*}]
\item\textsc{Case} 1.1: $n_0 = 1, n_1 = 0$ and $n_2 = n - 5$;
\item\textsc{Case} 1.2: $n_0 = 0, n_1 = 2$ and $n_2 = n - 6$.
\end{itemize}
Case 1.1 is impossible because if $u$ denotes the point of degree $0$ and $q$ is a point of degree $2$, then $\mathcal{B}\vert_{\{q, u, w, x, y, z\}}$ would contradict Claim \ref{Cn6}.
As for Case 1.2, notice that there are an equal number of $wx$- and $yz$-triangles because $d_\mathcal{B}(x) = d_\mathcal{B}(y)$, hence, the two points of degree $1$ cannot be covered by the same tight star. This gives exactly Case \ref{Estruct1} of Claim \ref{Cstruct}.
\end{case}

\begin{case}{2}
Claim \ref{Ccases} guarantees that the set of triangles of $\mathcal{B}$ can be divided into $xy$-triangles and $yz$-triangles, depending on whether they contain $\{x, y\}$ or $\{y, z\}$. Set $Y = X\backslash\{x, y, z\}$.
Then, for every point $p\in X$, $d_\mathcal{B}(p)\leq 2\Leftrightarrow p\in Y$, hence, $n_0 + n_1 + n_2 = |Y| = n - 3$. On the other hand, $n_1 + 2n_2 = 2n - 10$ as we have shown in the previous case, so one of the following cases hold:
%\begin{equation*}\begin{split}
%n_1 + 2n_2 & =\sum_{p\in Y} d_\mathcal{B}(p)\\
%& = |\mathcal{B}|_\Delta\\
%& = 2n - 10,
%\end{split}\end{equation*} so one of the following cases hold:
%and $$\sum_{p\in Y} d_\mathcal{B}(p)\leq 2(n - 3).$$
%Note that $n_i > 0\Rightarrow i\in\{0, 1, 2, n - 5, n - c\}$ and $n - 5 > 2$, so $n_0 + n_1 + n_2 = n - 3$, $n_{n - 5} = 2$ and $n_{n - c} = 1$ (see Claim \ref{Ccases}). Now, for every point $p\in Y$, $d_\mathcal{B}(p)\leq 2$, hence, $$\sum_{p\in Y} d_\mathcal{B}(p)\leq 2(n - 3).$$
%On the other hand, $$\sum_{p\in Y} d_\mathcal{B}(p) = |\mathcal{B}|_\Delta = 2n - 10,$$ so one of the following cases hold:
\begin{itemize}
%[label={(\alph*)}, ref={\alph*}]
\item\textsc{Case} 2.1: $n_0 = 2, n_1 = 0$ and $n_2 = n - 5$;
\item\textsc{Case} 2.2: $n_0 = 1, n_1 = 2$ and $n_2 = n - 6$;
\item\textsc{Case} 2.3: $n_0 = 0, n_1 = 4$ and $n_2 = n - 7$.
\end{itemize}
Note that there are an equal number of $xy$- and $yz$-triangles because $d_\mathcal{B}(x) = d_\mathcal{B}(z)$, hence, in each case, exactly half of the points of degree $1$ are covered by $xy$-triangles. Now, we can see that Case 2.1 and Case 2.3 coincide with Case \ref{Estruct2} and Case \ref{Estruct3} of Claim \ref{Cstruct}, respectively.
As for Case 2.2, let $q$ denote the point of degree $0$, and $u$ and $v$ be the two points of degree $1$. Then $\mathcal{B}\vert_{\{q, u, v, x, y, z\}}$ contradicts Claim \ref{Cn6}.
\end{case}
$\square$
\end{prf}

Next, we analyze Cases \ref{Estruct1}, \ref{Estruct2} and \ref{Estruct3} of Claim \ref{Cstruct} in order to characterize $\mathcal{B}$.

\begin{figure}
\centering
\begin{subfigure}[b]{0.33\textwidth}
\centering
\includegraphics{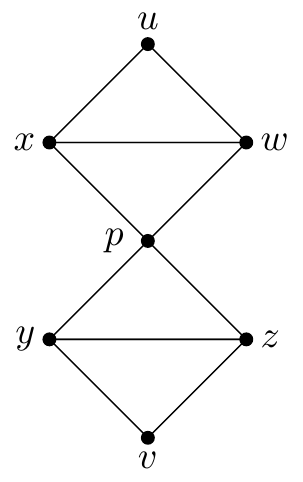}
%\begin{tikzpicture}[pont/.style={circle, fill=black, inner sep=0.5mm}]
%
%\begin{scope}[rotate=90]
%\node (A) at (0,0) [pont] {};
%\node (B) at (1,1) [pont] {};
%\node (C) at (1,-1) [pont] {};
%\node (D) at (2,0) [pont] {};
%\node (E) at (3,1) [pont] {};
%\node (F) at (3,-1) [pont] {};
%\node (G) at (4,0) [pont] {};
%
%\draw (A) node[anchor=north] {\small$v$};
%\draw (B) node[anchor=east] {\small$y$};
%\draw (C) node[anchor=west] {\small$z$};
%\draw (2,0.1) node[anchor=east] {\small$p$};
%\draw (E) node[anchor=east] {\small$x$};
%\draw (F) node[anchor=west] {\small$w$};
%\draw (G) node[anchor=south] {\small$u$};
%
%%\draw[thick, pattern=vertical lines] (0,0) -- (1,1) -- (1,-1) -- cycle;
%%\draw[thick, pattern=horizontal lines] (1,1) -- (1,-1) -- (2,0) -- cycle;
%%\draw[thick, pattern=vertical lines] (2,0) -- (3,1) -- (3,-1) -- cycle;
%%\draw[thick, pattern=horizontal lines] (3,1) -- (3,-1) -- (4,0) -- cycle;
%\draw (0,0) -- (1,1) -- (1,-1) -- cycle;
%\draw (1,1) -- (1,-1) -- (2,0) -- cycle;
%\draw (2,0) -- (3,1) -- (3,-1) -- cycle;
%\draw (3,1) -- (3,-1) -- (4,0) -- cycle;
%\end{scope}
%\end{tikzpicture}
\caption{$\tilde{\mathcal{H}}_1$}\label{SFmetr71}
\end{subfigure}%
\begin{subfigure}[b]{0.33\textwidth}
\centering
\includegraphics{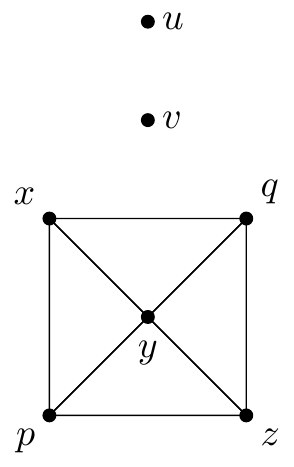}
%\begin{tikzpicture}[pont/.style={circle, fill=black, inner sep=0.5mm}]
%
%\begin{scope}[rotate=90]
%\node (A) at (0,0) [pont] {};
%\node (B) at (1,1) [pont] {};
%\node (C) at (-1,1) [pont] {};
%\node (D) at (-1,-1) [pont] {};
%\node (E) at (1,-1) [pont] {};
%\node (F) at (2,0) [pont] {};
%\node (G) at (3,0) [pont] {};
%
%%\draw[thick, pattern=vertical lines] (0,0) -- (1,1) -- (-1,1) -- cycle;
%%\draw[thick, pattern=horizontal lines] (0,0) -- (-1,1) -- (-1,-1) -- cycle;
%%\draw[thick, pattern=vertical lines] (0,0) -- (-1,-1) -- (1,-1) -- cycle;
%%\draw[thick, pattern=horizontal lines] (0,0) -- (1,-1) -- (1,1) -- cycle;
%\draw (0,0) -- (1,1) -- (-1,1) -- cycle;
%\draw (0,0) -- (-1,1) -- (-1,-1) -- cycle;
%\draw (0,0) -- (-1,-1) -- (1,-1) -- cycle;
%\draw (0,0) -- (1,-1) -- (1,1) -- cycle;
%
%\draw[fill=white,white] (-0.2,-0.14) rectangle (-0.52,0.14);
%\draw (-0.1,0) node[anchor=north] {\small$y$};
%\draw (B) node[anchor=south east] {\small$x$};
%\draw (C) node[anchor=north east] {\small$p$};
%\draw (D) node[anchor=north west] {\small$z$};
%\draw (E) node[anchor=south west] {\small$q$};
%\draw (F) node[anchor=west] {\small$v$};
%\draw (G) node[anchor=west] {\small$u$};
%\end{scope}
%\end{tikzpicture}
\caption{$\tilde{\mathcal{H}}_2$}\label{SFmetr72}
\end{subfigure}%
\begin{subfigure}[b]{0.33\textwidth}
\centering
\includegraphics{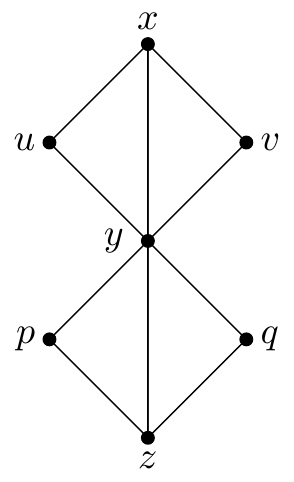}
%\begin{tikzpicture}[pont/.style={circle, fill=black, inner sep=0.5mm}]
%
%\begin{scope}[rotate=90]
%\node (A) at (0,0) [pont] {};
%\node (B) at (1,1) [pont] {};
%\node (C) at (1,-1) [pont] {};
%\node (D) at (2,0) [pont] {};
%\node (E) at (3,1) [pont] {};
%\node (F) at (3,-1) [pont] {};
%\node (G) at (4,0) [pont] {};
%
%\draw (A) node[anchor=north] {\small$z$};
%\draw (B) node[anchor=east] {\small$p$};
%\draw (C) node[anchor=west] {\small$q$};
%\draw (2,0.1) node[anchor=east] {\small$y$};
%\draw (E) node[anchor=east] {\small$u$};
%\draw (F) node[anchor=west] {\small$v$};
%\draw (G) node[anchor=south] {\small$x$};
%
%%\draw[thick, pattern=vertical lines] (0,0) -- (1,1) -- (2,0) -- cycle;
%%\draw[thick, pattern=horizontal lines] (0,0) -- (1,-1) -- (2,0) -- cycle;
%%\draw[thick, pattern=horizontal lines] (2,0) -- (3,1) -- (4,0) -- cycle;
%%\draw[thick, pattern=vertical lines] (2,0) -- (3,-1) -- (4,0) -- cycle;
%\draw (0,0) -- (1,1) -- (2,0) -- cycle;
%\draw (0,0) -- (1,-1) -- (2,0) -- cycle;
%\draw (2,0) -- (3,1) -- (4,0) -- cycle;
%\draw (2,0) -- (3,-1) -- (4,0) -- cycle;
%\end{scope}
%\end{tikzpicture}
\caption{$\tilde{\mathcal{H}}_3$}\label{SFmetr73}
\end{subfigure}%
\caption{Triangle hypergraphs in Claim \ref{Cmetr7}. The hyperedges are represented by triangles.}\label{Fmetr7}
\end{figure}

\begin{clm}\label{Cmetr7}
The following statements hold for the hypergraphs $\tilde{\mathcal{H}}_1$, $\tilde{\mathcal{H}}_2$ and $\tilde{\mathcal{H}}_3$ in Figure \ref{Fmetr7}:
\begin{enumerate}
\item $\tilde{\mathcal{H}}_1$ is metrizable, and for every betweenness structure $\mathcal{A}$ with triangle hypergraph $\tilde{\mathcal{H}}_1$, $\mathcal{A}\simeq\mathcal{T}_{7, 1}$;\label{Emetr71}
\item $\tilde{\mathcal{H}}_2$ is not metrizable;\label{Emetr72}
\item $\tilde{\mathcal{H}}_3$ is not metrizable.\label{Emetr73}
\end{enumerate}
\end{clm}

The proof of Claim \ref{Cmetr7} can be found in \cite{gnszabo2018qusapxprep} along with the other technical lemmas.

\begin{case}{\ref{Estruct1}}
Recall that there exist distinct points $u, v, w, x, y, z\in X$ such that with $X' = X\backslash\{u, v, w, x, y, z\}$, $$\Delta(\mathcal{B}) =\{\{p, w, x\} : p\in X'\}\cup\{\{p, y, z\} : p\in X'\}\cup\{\{u, w, x\}, \{v, y, z\}\}.$$
We show that $\mathcal{B}\simeq\mathcal{T}_{n, i}$ for some $1\leq i\leq\left\lceil\frac{n - 5}{2}\right\rceil$.

\begin{figure}[t]
\centering
\begin{subfigure}[b]{0.5\textwidth}
\centering
\includegraphics{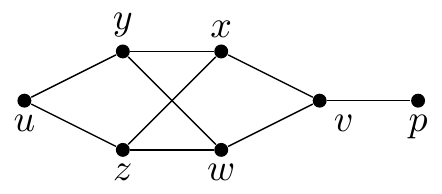}
%\begin{tikzpicture}[scale=1, pont/.style={circle, fill=black, inner sep=0.5mm},
%mybox/.style={rectangle, rounded corners, draw=black, fill=white, inner sep=0.5mm}]
%
%\node (U) at (0,0) [pont] {};
%\node (Z) at (1,-0.5) [pont] {};
%\node (Y) at (1,0.5) [pont] {};
%\node (W) at (2,-0.5) [pont] {};
%\node (X) at (2,0.5) [pont] {};
%\node (V) at (3,0) [pont] {};
%\node (P) at (4,0) [pont] {};
%
%\draw (U) node[anchor=north] {\small$u$};
%\draw (Y) node[anchor=south] {\small$y$};
%\draw (Z) node[anchor=north] {\small$z$};
%\draw (X) node[anchor=south] {\small$x$};
%\draw (W) node[anchor=north] {\small$w$};
%\draw (V) node[anchor=north west] {\small$v$};
%\draw (P) node[anchor=north] {\small$p$};
%
%\draw (U)--(Y) node {};
%\draw (U)--(Z) node {};
%\draw (Y)--(X) node {};
%\draw (Y)--(W) node {};
%\draw (Z)--(X) node {};
%\draw (Z)--(W) node {};
%\draw (X)--(V) node {};
%\draw (W)--(V) node {};
%\draw (V)--(P) node {};
%\end{tikzpicture}
\caption{}\label{SFtpi1}
\end{subfigure}%
\begin{subfigure}[b]{0.5\textwidth}
\centering
\includegraphics{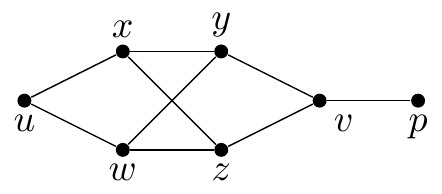}
%\begin{tikzpicture}[scale=1, pont/.style={circle, fill=black, inner sep=0.5mm},
%mybox/.style={rectangle, rounded corners, draw=black, fill=white, inner sep=0.5mm}]
%
%\node (U) at (0,0) [pont] {};
%\node (W) at (1,-0.5) [pont] {};
%\node (X) at (1,0.5) [pont] {};
%\node (Z) at (2,-0.5) [pont] {};
%\node (Y) at (2,0.5) [pont] {};
%\node (V) at (3,0) [pont] {};
%\node (P) at (4,0) [pont] {};
%
%\draw (U) node[anchor=north] {\small$u$};
%\draw (Y) node[anchor=south] {\small$y$};
%\draw (Z) node[anchor=north] {\small$z$};
%\draw (X) node[anchor=south] {\small$x$};
%\draw (W) node[anchor=north] {\small$w$};
%\draw (V) node[anchor=north west] {\small$v$};
%\draw (P) node[anchor=north] {\small$p$};
%
%\draw (U)--(X) node {};
%\draw (U)--(W) node {};
%\draw (Y)--(X) node {};
%\draw (Y)--(W) node {};
%\draw (Z)--(X) node {};
%\draw (Z)--(W) node {};
%\draw (Y)--(V) node {};
%\draw (Z)--(V) node {};
%\draw (V)--(P) node {};
%\end{tikzpicture}
\caption{}
\end{subfigure}%
\vspace{1em}
\begin{subfigure}[b]{0.5\textwidth}
\centering
\includegraphics{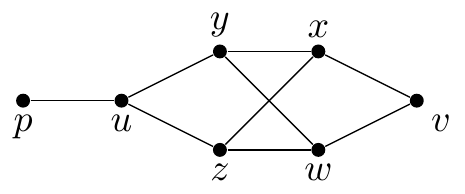}
%\begin{tikzpicture}[scale=1, pont/.style={circle, fill=black, inner sep=0.5mm},
%mybox/.style={rectangle, rounded corners, draw=black, fill=white, inner sep=0.5mm}]
%
%\node (P) at (0,0) [pont] {};
%\node (U) at (1,0) [pont] {};
%\node (Z) at (2,-0.5) [pont] {};
%\node (Y) at (2,0.5) [pont] {};
%\node (W) at (3,-0.5) [pont] {};
%\node (X) at (3,0.5) [pont] {};
%\node (V) at (4,0) [pont] {};
%
%\draw (P) node[anchor=north] {\small$p$};
%\draw (U) node[anchor=north] {\small$u$};
%\draw (Y) node[anchor=south] {\small$y$};
%\draw (Z) node[anchor=north] {\small$z$};
%\draw (X) node[anchor=south] {\small$x$};
%\draw (W) node[anchor=north] {\small$w$};
%\draw (V) node[anchor=north west] {\small$v$};
%
%\draw (P)--(U) node {};
%\draw (U)--(Y) node {};
%\draw (U)--(Z) node {};
%\draw (Y)--(X) node {};
%\draw (Y)--(W) node {};
%\draw (Z)--(X) node {};
%\draw (Z)--(W) node {};
%\draw (X)--(V) node {};
%\draw (W)--(V) node {};
%\end{tikzpicture}
\caption{}
\end{subfigure}%
\begin{subfigure}[b]{0.5\textwidth}
\centering
\includegraphics{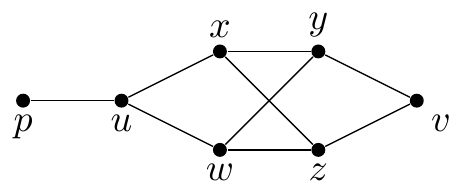}
%\begin{tikzpicture}[scale=1, pont/.style={circle, fill=black, inner sep=0.5mm},
%mybox/.style={rectangle, rounded corners, draw=black, fill=white, inner sep=0.5mm}]
%
%\node (P) at (0,0) [pont] {};
%\node (U) at (1,0) [pont] {};
%\node (W) at (2,-0.5) [pont] {};
%\node (X) at (2,0.5) [pont] {};
%\node (Z) at (3,-0.5) [pont] {};
%\node (Y) at (3,0.5) [pont] {};
%\node (V) at (4,0) [pont] {};
%
%\draw (P) node[anchor=north] {\small$p$};
%\draw (U) node[anchor=north] {\small$u$};
%\draw (Y) node[anchor=south] {\small$y$};
%\draw (Z) node[anchor=north] {\small$z$};
%\draw (X) node[anchor=south] {\small$x$};
%\draw (W) node[anchor=north] {\small$w$};
%\draw (V) node[anchor=north west] {\small$v$};
%
%\draw (P)--(U) node {};
%\draw (U)--(X) node {};
%\draw (U)--(W) node {};
%\draw (Y)--(X) node {};
%\draw (Y)--(W) node {};
%\draw (Z)--(X) node {};
%\draw (Z)--(W) node {};
%\draw (Y)--(V) node {};
%\draw (Z)--(V) node {};
%\end{tikzpicture}
\caption{}
\end{subfigure}%
\caption{Possible spanner graphs of betweenness structure $\mathcal{B}'$ in the proof of Case A of Theorem \ref{Tlinsize2}.}\label{Ftpi}
\end{figure}

Let $p$ be a point of $X'$ and consider the betweenness structure $\mathcal{B}' =\mathcal{B}\vert_{\{p, u, v, w, x, y, z\}}$. It is easy to see that $\mathcal{H}(\mathcal{B}') =\tilde{\mathcal{H}}_1$, so we obtain from Part \ref{Emetr71} of Claim \ref{Cmetr7} that $\mathcal{B}'\simeq\mathcal{T}_{7, 1}$. Further, it is easy to see that $\mathcal{B}'$ is induced by one of the graphs in Figure \ref{Ftpi},
%In summary, we can conclude that for any point $p\in X'$, $\mathcal{B}\vert_{\{u, v, w, x, y, z, p\}}$ is induced by one of the graphs in Figure \ref{Ftpi} (taking all of the symmetries into account), 
therefore, 
\begin{equation}\label{EQ18h}
(p\ u\ v)_\mathcal{B}\text{ or }(u\ v\ p)_\mathcal{B}
\end{equation} holds. We can assume without loss of generality that $G(\mathcal{B}')$ is the graph in Figure \ref{SFtpi1} and thus \begin{equation}\label{EQ18g}
(u\ y\ x\ v)_\mathcal{B}, (u\ z\ w\ v)_\mathcal{B}, (u\ y\ w\ v)_\mathcal{B}\emph{ and }(u\ z\ x\ v)_\mathcal{B}
\end{equation} hold.

We know from the triangle hypergraph that $\mathcal{B}\vert_{X\backslash\{w, x, y, z\}}$ is linear, moreover, it is ordered by Proposition  \ref{Plin} as $n - 4\geq 5$. This is also true for $\mathcal{B}\vert_{X'}$ as $X'\subseteq X\backslash\{w, x, y, z\}$. Let $x_1, x_2,\ldots, x_{n - 6}$ denote the points of $X'$ such that $$\mathcal{B}\vert_{X'} = [x_1, x_2,\ldots, x_{n - 6}],$$ and let $k$ and $\ell$ denote the position of $u$ and $v$ in $\mathcal{B}\vert_{X'\cup\{u\}}$ and $\mathcal{B}\vert_{X'\cup\{v\}}$, respectively. Since no $p\in X'$ is between $u$ and $v$ by (\ref{EQ18h}),
 $k = \ell$, and we can assume without loss of generality that $$(x_1\ \ldots\ x_{k - 1}\ u\ v\ x_{k}\ \ldots\ x_{n - 6})_\mathcal{B}.$$ We can also assume that $k\leq\left\lceil\frac{n - 5}{2}\right\rceil$. Keep in mind, however, that $k = 1$ is possible. Now, it follows from (\ref{EQ18g}) that $$(x_1\ \ldots\ x_{k - 1}\ u\ y\ x\ v\ x_{k}\ \ldots\ x_{n - 6})_\mathcal{B},$$ $$(x_1\ \ldots\ x_{k - 1}\ u\ z\ w\ v\ x_{k}\ \ldots\ x_{n - 6})_\mathcal{B},$$ $$(x_1\ \ldots\ x_{k - 1}\ u\ y\ w\ v\ x_{k}\ \ldots\ x_{n - 6})_\mathcal{B}$$ and $$(x_1\ \ldots\ x_{k - 1}\ u\ z\ x\ v\ x_{k}\ \ldots\ x_{n - 6})_\mathcal{B}.$$ With this, we covered all collinear triples of $\mathcal{B}$, and the resulting betweenness structure is metrizable and isomorphic to $\mathcal{T}_{n, k}$.
\end{case}

\begin{case}{\ref{Estruct2}}
Recall that there exist distinct points $u, v, x, y, z\in X$ such that with $X' = X\backslash\{u, v, x, y, z\}$, $$\Delta(\mathcal{B}) =\{\{p, x, y\} : p\in X'\}\cup\{\{p, y, z\} : p\in X'\}.$$
Set $\mathcal{B}' =\mathcal{B}\vert_{\{p, q, u, v, x, y, z\}}$ where $p, q\in X'$ are two arbitrary points.
Now, $\mathcal{H}(\mathcal{B}') =\tilde{\mathcal{H}}_2$, which is shown to be impossible by Part \ref{Emetr72} of Claim \ref{Cmetr7}.

%%Notice that both $\mathcal{B}' - x$ and $\mathcal{B}' - q$ are a quasilinear betweenness structures with $2$ triangles that form a tight star. The next observation follows from Theorem \ref{Tqus1}.

%%\begin{obs}\label{Okern}
%%Let $\mathcal{A}$ be a quasilinear betweenness structure of order $n\geq 6$ such that $\mathcal{H}(\mathcal{A})$ is a tight star. Then there is exactly one cyclic line in $\mathcal{A}$ and it contains $\ker(\mathcal{H}(\mathcal{A}))$.
%%\end{obs}

%%%\begin{prf}
%%%Let $Y$ be the ground set of $\mathcal{A}$. Because of Theorem \ref{Tqus1}, $\mathcal{A}$ is isomorphic to either $\mathcal{R}_{n, i}^4$ ($i\in I_n^4$) or $\mathcal{S}_n^4$. It is easy to see that there is exactly one subset $Y' =\{y_1, y_2, y_3, y_4\}\subset Y$ such that the betweennesses $(y_1\ y_2\ y_3)_\mathcal{A}$, $(y_2\ y_3\ y_4)_\mathcal{A}$, $(y_3\ y_4\ y_1)_\mathcal{A}$ and $(y_4\ y_1\ y_2)_\mathcal{A}$ hold. Further, $Y'$ contains $\ker(\mathcal{H}(\mathcal{A}))$. $\square$
%%%\end{prf}

%%We obtain from Observation \ref{Okern} that there is exactly one cyclic line $L_x$ in $\mathcal{B}' - x$, and it contains the kernel $\{y, z\}$. 
%%%and there is exactly one cyclic line $L_q$ in $\mathcal{B}' - q$ and it contains the star kernel $\{y, p\}$.
%%Further, $L_x$ does not contain $p$ or $q$ because $\{y, z, p\}$ and $\{y, z, q\}$ are triangles. Hence, $L_x$ is a cyclic line in $\mathcal{B}' - q$ as well that does not contain $p$, contradicting Observation \ref{Okern}.
\end{case}

\begin{case}{\ref{Estruct3}}
Recall that there exist distinct points $p, q, u, v, x, y, z\in X$ such that with $X' = X\backslash\{p, q, u, v, x, y, z\}$, $$\Delta(\mathcal{B}) =\{\{p', x, y\} : p'\in X'\}\cup\{\{p', y, z\} : p'\in X'\}\cup$$ $$\{\{u, x, y\}, \{v, x, y\}, \{p, y, z\}, \{q, y, z\}\}.$$ Set $\mathcal{B}' =\mathcal{B}\vert_{\{p, q, u, v, x, y, z\}}$. Now, $\mathcal{H}(\mathcal{B}') =\tilde{\mathcal{H}}_3$, which is again impossible by Part \ref{Emetr73} of Claim \ref{Cmetr7}.

%%Observe that $\mathcal{B}' - x$ and $\mathcal{B}' - z$ are quasilinear betweenness structures on $6$ points with $2$ triangles that form a tight star, hence, they are isomorphic to either $\mathcal{R}_{6, 1}^4$, $\mathcal{R}_{6, 2}^4$ or $\mathcal{S}_6^4$.

%%If $\mathcal{B}' - x\simeq\mathcal{R}_{6, 1}^4$, then we can assume by symmetry that $\mathcal{B}' - x =\mathcal{R}_{6, 1}^4(y, z; u, v, p, q)$. It is easy to see by Observation \ref{O6qlin} that $\mathcal{B}' - z =\mathcal{R}_{6, 1}^4(y,\allowbreak x;\allowbreak u,\allowbreak v,\allowbreak p,\allowbreak q)$, from which $(y\ u\ x)_\mathcal{B}$ follows. This is, however, impossible as $\{u, x, y\}$ was a triangle.
%%Next, if $\mathcal{B}' - x\simeq\mathcal{R}_{6, 2}^4$, then we can assume by symmetry that $\mathcal{B}' - x =\mathcal{R}_{6, 2}^4(y, z; p, u, v, q)$. It is easy to see that $\mathcal{B}' - z =\mathcal{R}_{6, 2}^4(y, x; p, u, v, q)$ by Observation \ref{O6qlin}, from which $(y\ u\ x)_\mathcal{B}$ follows, leading to a contradiction again.
%%Finally, if $\mathcal{B}' - x\simeq\mathcal{S}_6^4$, then we can assume by symmetry that $\mathcal{B}' - x =\mathcal{S}_6^4(y, z; u, p, q, v)$. Now, we obtain from Observation \ref{O6qlin} that $\mathcal{B}' - z =\mathcal{S}_6^4(y, x; u, p, q, v)$, from which $(u\ y\ x)_\mathcal{B}$ follows, a contradiction again. In summary, Case \ref{Estruct3} is impossible.
\end{case}

\ref{Tlinsize23}. Finally, the following examples prove the last point of Theorem \ref{Tlinsize2}: take $\mathcal{S}_6^4$ for $n = 6$, $\mathcal{S}_7^3$ for $n = 7$ and $\mathcal{S}_8^2$ for $n = 8$.
$\square$
\end{prfof}

We remark that Theorem \ref{Tqus2} and Theorem \ref{Tqus3} follows from Theorem \ref{Tlinsize1} if $n$ is large enough ($n\geq 8$ and $n\geq 9$, respectively). The reason why we have chosen another path to prove them is that Theorem \ref{Tlinsize1} does not say anything about small betweenness structures and complicated case analysis would still be required to characterize them.

%$$B(n, kn - c)\neq\emptyset\Rightarrow c\geq 2k.$$

%---------------------------------------------------------------------------------------------------

\section{Conclusion}\label{Sconc}

Extending the work of Richmond, Richmond \cite{richmond1997metric} and Beaudou \cite{beaudou2013lines} on hypergraph metrizability, we have characterized almost-metrizable betweenness structures of small linear co-size, including the largest non-linear (quasilinear) betweenness structures and the largest betweenness structures of co-size $2n - c$. We have observed that there are gaps in the size-spectrum of betweenness structures and we proposed interesting quantities that can be subject of future research.

We close the paper with our conjectures on the possible extension of Theorem \ref{Tlinsize1} and Theorem \ref{Tlinsize2} to betweenness structures of co-size $kn - c$.
We call a $3$-uniform hypergraph a \emph{tight $k$-star} if it is the (not necessarily edge-disjoint) union of $k$ tight stars on the same ground set.
% there exist $k$ pairs of vertices $\{x_1, y_1\}, \{x_2, y_2\},\ldots,\{x_k, y_k\}$ such that each edge of $\mathcal{H}$ contains one of these pairs.

\begin{cnj}\label{CNktight}
For all integers $k\geq 0$ and $c$, there exists a threshold $N_c^k > 0$ such that for all $n\geq N_c^k$ and betweenness structure $\mathcal{B}\in B(n, kn - c)$, $\mathcal{H}(\mathcal{B})$ is a tight $k$-star.
\end{cnj}

Further, the following inequality would be an easy consequence of Conjecture \ref{CNktight}.
\begin{cnj}
$$\vartheta_{\min}(k)\geq 2k.$$
\end{cnj}

%---------------------------------------------------------------------------------------------------

%\section{Appendix}\label{Sapp}

%\subsection{Characterization of Small Quasilinear Betweenness Structures}
%In this section, we characterize small quasilinear betweenness structures by proving Lemma \ref{Lsmallgr}.

\section*{Acknowledgment}
We are grateful to Pierre Aboulker for sharing his thoughts on some of the results presented here.

This work was supported by the National Research, Development and Innovation Office -- NKFIH, No. 108947.

\storecounter{thr}{L:thr}
\storecounter{dfn}{L:dfn}
\storecounter{lmm}{L:lmm}
\storecounter{prp}{L:prp}
\storecounter{crl}{L:crl}
\storecounter{clm}{L:clm}
\storecounter{fct}{L:fct}
\storecounter{obs}{L:obs}
\storecounter{cnj}{L:cnj}
\storecounter{prb}{L:prb}
\storecounter{figure}{L:figure}

\bibliographystyle{elsarticle-num}
\bibliography{fms}

\end{document}

\maketitle

\getcounter{thr}{L:thr}
\getcounter{dfn}{L:dfn}
\getcounter{lmm}{L:lmm}
\getcounter{prp}{L:prp}
\getcounter{crl}{L:crl}
\getcounter{clm}{L:clm}
\getcounter{fct}{L:fct}
\getcounter{obs}{L:obs}
\getcounter{cnj}{L:cnj}
\getcounter{prb}{L:prb}
\getcounter{figure}{L:figure}

\section{Results}

This is an appendix to \cite{gnszabo2018qusarxiv}, and contains the proofs of some technical results for small betweenness structures. We use definitions and references (e.g. theorem and figure numbering) from \cite{gnszabo2018qusarxiv} but we also repeat the results to be proved here along with the corresponding figures.

We prove the following results.

\setcounter{foo}{\thefigure}
\setcounterref{figure}{Fsmallgr}
\addtocounter{figure}{-1}
\begin{figure}[t]
\centering
\includegraphics{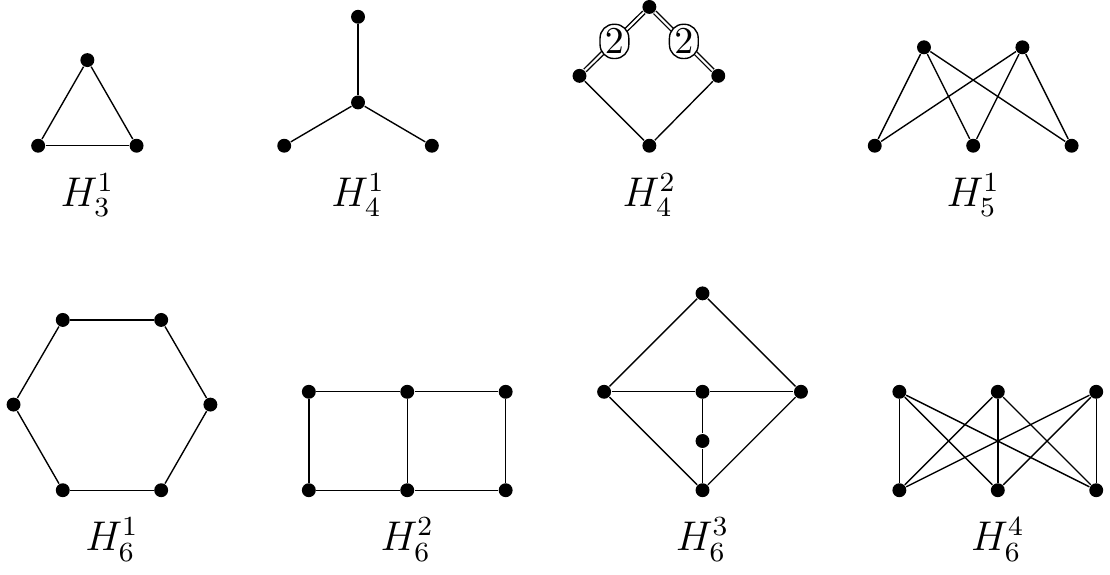}
\caption{Spanner graphs of small exceptional quasilinear betweenness structures. Edges of weight different from $1$ are indicated by double-lines.}
\end{figure}
\setcounter{figure}{\thefoo}

\setcounter{foo}{\theclm}
\setcounterref{lmm}{Lsmallgr}
\addtocounter{lmm}{-1}
\begin{lmm}
Up to isomorphism, the quasilinear betweenness structures of order at most $7$ are the following:
%A betweenness structure of order $1\leq n\leq 7$ is quasilinear if and only if it is isomorphic to one of the following betweenness structures:
\begin{itemize}
\item $\mathcal{R}_{n, i}^4$ for $5\leq n\leq 7$, $i\in I_n^4$;
\item $\mathcal{S}_n^4$ for $5\leq n\leq 7$;
\item $\mathcal{B}(G)$ where $G$ is one of the graphs in Figure \ref{Fsmallgr}.
\end{itemize}
\end{lmm}
\setcounter{lmm}{\thefoo}

\setcounter{foo}{\theclm}
\setcounterref{clm}{Cn5}
\addtocounter{clm}{-1}
\begin{clm}
Let $\mathcal{B}\in B(5, 2)$ be a regular betweenness structure. Then $\mathcal{H}(\mathcal{B})$ is a tight star.
\end{clm}
\setcounter{clm}{\thefoo}

\setcounter{foo}{\thefigure}
\setcounterref{figure}{Fmetr7}
\addtocounter{figure}{-1}
\begin{figure}
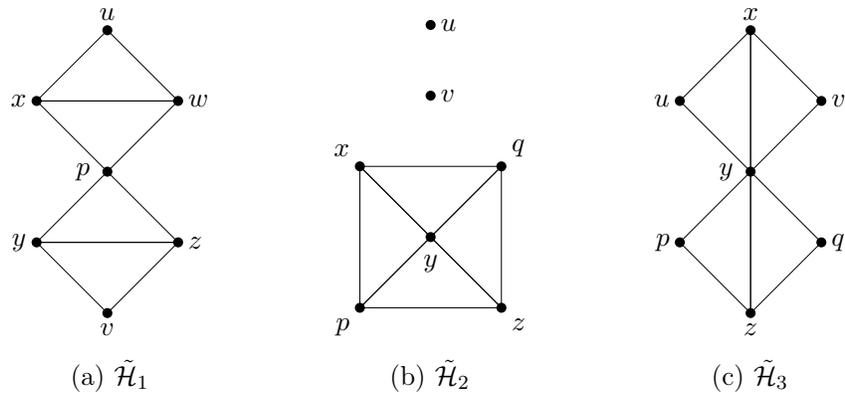

\centering
\begin{subfigure}[b]{0.33\textwidth}
\centering
\includegraphics{fms_qus-figure8}
%\begin{tikzpicture}[pont/.style={circle, fill=black, inner sep=0.5mm}]
%
%\begin{scope}[rotate=90]
%\node (A) at (0,0) [pont] {};
%\node (B) at (1,1) [pont] {};
%\node (C) at (1,-1) [pont] {};
%\node (D) at (2,0) [pont] {};
%\node (E) at (3,1) [pont] {};
%\node (F) at (3,-1) [pont] {};
%\node (G) at (4,0) [pont] {};
%
%\draw (A) node[anchor=north] {\small$v$};
%\draw (B) node[anchor=east] {\small$y$};
%\draw (C) node[anchor=west] {\small$z$};
%\draw (2,0.1) node[anchor=east] {\small$p$};
%\draw (E) node[anchor=east] {\small$x$};
%\draw (F) node[anchor=west] {\small$w$};
%\draw (G) node[anchor=south] {\small$u$};
%
%\draw[thick, pattern=vertical lines] (0,0) -- (1,1) -- (1,-1) -- cycle;
%\draw[thick, pattern=horizontal lines] (1,1) -- (1,-1) -- (2,0) -- cycle;
%\draw[thick, pattern=vertical lines] (2,0) -- (3,1) -- (3,-1) -- cycle;
%\draw[thick, pattern=horizontal lines] (3,1) -- (3,-1) -- (4,0) -- cycle;
%\end{scope}
%\end{tikzpicture}
\caption{$\tilde{\mathcal{H}}_1$}%\label{SFmetr71}
\end{subfigure}%
\begin{subfigure}[b]{0.33\textwidth}
\centering
\includegraphics{fms_qus-figure9}
%\begin{tikzpicture}[pont/.style={circle, fill=black, inner sep=0.5mm}]
%
%\begin{scope}[rotate=90]
%\node (A) at (0,0) [pont] {};
%\node (B) at (1,1) [pont] {};
%\node (C) at (-1,1) [pont] {};
%\node (D) at (-1,-1) [pont] {};
%\node (E) at (1,-1) [pont] {};
%\node (F) at (2,0) [pont] {};
%\node (G) at (3,0) [pont] {};
%
%\draw[thick, pattern=vertical lines] (0,0) -- (1,1) -- (-1,1) -- cycle;
%\draw[thick, pattern=horizontal lines] (0,0) -- (-1,1) -- (-1,-1) -- cycle;
%\draw[thick, pattern=vertical lines] (0,0) -- (-1,-1) -- (1,-1) -- cycle;
%\draw[thick, pattern=horizontal lines] (0,0) -- (1,-1) -- (1,1) -- cycle;
%
%\draw[fill=white,white] (-0.2,-0.14) rectangle (-0.52,0.14);
%\draw (-0.1,0) node[anchor=north] {\small$y$};
%\draw (B) node[anchor=south east] {\small$x$};
%\draw (C) node[anchor=north east] {\small$p$};
%\draw (D) node[anchor=north west] {\small$z$};
%\draw (E) node[anchor=south west] {\small$q$};
%\draw (F) node[anchor=west] {\small$v$};
%\draw (G) node[anchor=west] {\small$u$};
%\end{scope}
%\end{tikzpicture}
\caption{$\tilde{\mathcal{H}}_2$}%\label{SFmetr72}
\end{subfigure}%
\begin{subfigure}[b]{0.33\textwidth}
\centering
\includegraphics{fms_qus-figure10}
%\begin{tikzpicture}[pont/.style={circle, fill=black, inner sep=0.5mm}]
%
%\begin{scope}[rotate=90]
%\node (A) at (0,0) [pont] {};
%\node (B) at (1,1) [pont] {};
%\node (C) at (1,-1) [pont] {};
%\node (D) at (2,0) [pont] {};
%\node (E) at (3,1) [pont] {};
%\node (F) at (3,-1) [pont] {};
%\node (G) at (4,0) [pont] {};
%
%\draw (A) node[anchor=north] {\small$z$};
%\draw (B) node[anchor=east] {\small$p$};
%\draw (C) node[anchor=west] {\small$q$};
%\draw (2,0.1) node[anchor=east] {\small$y$};
%\draw (E) node[anchor=east] {\small$u$};
%\draw (F) node[anchor=west] {\small$v$};
%\draw (G) node[anchor=south] {\small$x$};
%
%\draw[thick, pattern=vertical lines] (0,0) -- (1,1) -- (2,0) -- cycle;
%\draw[thick, pattern=horizontal lines] (0,0) -- (1,-1) -- (2,0) -- cycle;
%\draw[thick, pattern=horizontal lines] (2,0) -- (3,1) -- (4,0) -- cycle;
%\draw[thick, pattern=vertical lines] (2,0) -- (3,-1) -- (4,0) -- cycle;
%\end{scope}
%\end{tikzpicture}
\caption{$\tilde{\mathcal{H}}_3$}%\label{SFmetr73}
\end{subfigure}%
\caption{Triangle hypergraphs in Claim \ref{Cmetr7}}%\label{Fmetr7}
\end{figure}
\setcounter{figure}{\thefoo}

\setcounter{foo}{\theclm}
\setcounterref{clm}{Cmetr7}
\addtocounter{clm}{-1}
\begin{clm}
The following statements hold for the hypergraphs $\tilde{\mathcal{H}}_1$, $\tilde{\mathcal{H}}_2$ and $\tilde{\mathcal{H}}_3$ in Figure \ref{Fmetr7}:
\begin{enumerate}
\item $\tilde{\mathcal{H}}_1$ is metrizable, and for every betweenness structure $\mathcal{A}$ with triangle hypergraph $\tilde{\mathcal{H}}_1$, $\mathcal{A}\simeq\mathcal{T}_{7, 1}$;
\item $\tilde{\mathcal{H}}_2$ is not metrizable;
\item $\tilde{\mathcal{H}}_3$ is not metrizable.
\end{enumerate}
\end{clm}
\setcounter{clm}{\thefoo}

\section{Proof of Lemma \ref{Lsmallgr}}

%\begin{prfof}{Lemma \ref{Lsmallgr}}
Let $\mathcal{B} = (X,\beta)$ be a quasilinear betweenness structure of order $n\leq 7$ and set $\mathcal{H} =\mathcal{H}(\mathcal{B})$. For the sake of this proof, we will denote the betweenness structure induced by the graph $H_n^i$ from Figure \ref{Fsmallgr} by $\mathcal{A}_n^i$.

%Let us first deal with the case $3\leq n\leq 4$.
There are no quasilinear betweenness structures of order $n < 3$, and if $n = 3$, then $\mathcal{B}$ is obviously induced by a triangle.
If $n = 4$, then $\mathcal{B}$ is of co-size $1$, so $\mathcal{H}$ is obviously a tight star. Thus, we can apply Lemma \ref{Lgen3} with $c = 3$ to obtain that $\mathcal{B}\simeq\mathcal{Q}_4^3\simeq\mathcal{A}_4^1$ or $\mathcal{B}\simeq\mathcal{R}_{4, 1}^3\simeq\mathcal{S}_4^3\simeq\mathcal{A}_4^2$.
%Let $X = \{x_1, x_2, x_3, x_4\}$ be the ground set of $\mathcal{B}$ and set $G = G(\mathcal{B})$. Notice that a for any betweenness $(x\ y\ z)_\mathcal{B}$, \begin{equation}\label{EQ18a}
%\{x, y\}, \{y, z\}\in E(G),
%\end{equation} otherwise $\mathcal{B}$ would be linear by f.r.p. This also means that each collinear triple induces a $2$-path in $G$. As $|\mathcal{B}|_\Delta = 1$, there are $3$ collinear triples in $\mathcal{B}$ and hence $3$ induced $2$-paths in $G$. It can be easily seen that among all of the (connected) simple graphs on $4$ points, only the cycle $C_4$ and the star $H_4^1$ meet these requirements. Now by (\ref{EQ18a}), there are three possible betweennesses with $G = H_4^1$ which form the betweenness structure $\mathcal{A}_4^1$, while chosing $3$ out of the $4$ possible betweennesses on $C_4$ gives a betweenness structure isomorphic to $\mathcal{A}_4^2$.
%induced by the second (weighted) graph in row $n = 4$ of Figure \ref{Fsmallgr}.

\begin{figure}
\centering
\begin{subfigure}[b]{0.5\textwidth}
\centering
\includegraphics{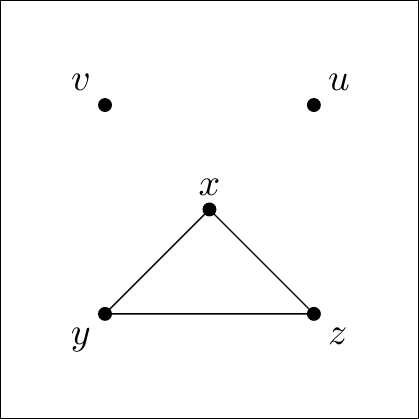}
%\begin{tikzpicture}[pont/.style={circle, fill=black, inner sep=0.5mm}]
%
%\draw (-2.121,-2.121) rectangle (2.121, 2.121);
%
%\node (U) at (1.061,1.061) [pont] {};
%\node (V) at (-1.061,1.061) [pont] {};
%\node (X) at (0,0) [pont] {};
%\node (Y) at (-1.061,-1.061) [pont] {};
%\node (Z) at (1.061,-1.061) [pont] {};
%
%%\draw[thick, pattern=vertical lines] (0,0) -- (-1.061,-1.061) -- (1.061,-1.061) -- cycle;
%\draw (0,0) -- (-1.061,-1.061) -- (1.061,-1.061) -- cycle;
%
%
%\draw (U) node[anchor=south west] {\small$u$};
%\draw (V) node[anchor=south east] {\small$v$};
%\draw (X) node[anchor=south] {\small$x$};
%\draw (Y) node[anchor=north east] {\small$y$};
%\draw (Z) node[anchor=north west] {\small$z$};
%\end{tikzpicture}
\caption{Case 1}\label{SFthyp1}
\end{subfigure}%
\begin{subfigure}[b]{0.5\textwidth}
\centering
\includegraphics{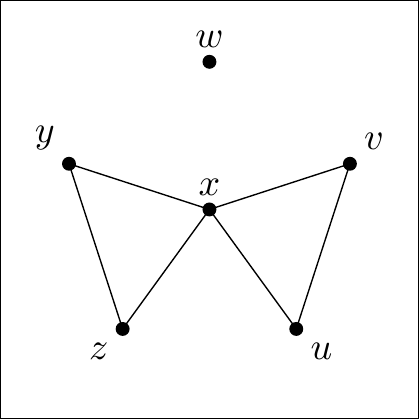}
%\begin{tikzpicture}[pont/.style={circle, fill=black, inner sep=0.5mm}]
%
%\draw (-2.121,-2.121) rectangle (2.121, 2.121);
%
%\node (X) at (0,0) [pont] {};
%\node (Y) at (-1.427,0.464) [pont] {};
%\node (Z) at (-0.882,-1.214) [pont] {};
%\node (U) at (0.882,-1.214) [pont] {};
%\node (V) at (1.427,0.464) [pont] {};
%\node (W) at (-0,1.5) [pont] {};
%
%%\draw[thick, pattern=north east lines] (0,0) -- (-1.427,0.464) -- (-0.882,-1.214) -- cycle;
%%\draw[thick, pattern=north west lines] (0,0) -- (0.882,-1.214) -- (1.427,0.464) -- cycle;
%\draw (0,0) -- (-1.427,0.464) -- (-0.882,-1.214) -- cycle;
%\draw (0,0) -- (0.882,-1.214) -- (1.427,0.464) -- cycle;
%
%\draw (X) node[anchor=south] {\small$x$};
%\draw (Y) node[anchor=south east] {\small$y$};
%\draw (Z) node[anchor=north east] {\small$z$};
%\draw (U) node[anchor=north west] {\small$u$};
%\draw (V) node[anchor=south west] {\small$v$};
%\draw (W) node[anchor=south] {\small$w$};
%\end{tikzpicture}
\caption{Case 2}\label{SFthyp2}
\end{subfigure}%
\vspace{1em}
\begin{subfigure}[b]{0.5\textwidth}
\centering
\includegraphics{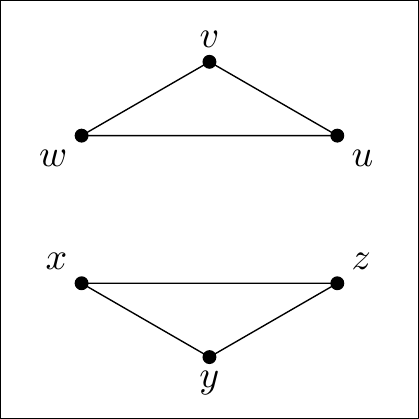}
%\begin{tikzpicture}[pont/.style={circle, fill=black, inner sep=0.5mm}]
%
%\draw (-2.121,-2.121) rectangle (2.121, 2.121);
%
%\node (X) at (-1.299,-0.75) [pont] {};
%\node (Y) at (0,-1.5) [pont] {};
%\node (Z) at (1.299,-0.75) [pont] {};
%\node (U) at (1.299,0.75) [pont] {};
%\node (V) at (0,1.5) [pont] {};
%\node (W) at (-1.299,0.75) [pont] {};
%
%%\draw[thick, pattern=vertical lines] (-1.299,-0.75) -- (0,-1.5) -- (1.299,-0.75) -- cycle;
%%\draw[thick, pattern=vertical lines] (1.299,0.75) -- (0,1.5) -- (-1.299,0.75) -- cycle;
%\draw (-1.299,-0.75) -- (0,-1.5) -- (1.299,-0.75) -- cycle;
%\draw (1.299,0.75) -- (0,1.5) -- (-1.299,0.75) -- cycle;
%
%\draw (X) node[anchor=south east] {\small$x$};
%\draw (Y) node[anchor=north] {\small$y$};
%\draw (Z) node[anchor=south west] {\small$z$};
%\draw (U) node[anchor=north west] {\small$u$};
%\draw (V) node[anchor=south] {\small$v$};
%\draw (W) node[anchor=north east] {\small$w$};
%\end{tikzpicture}
\caption{Case 3}\label{SFthyp3}
\end{subfigure}%
\begin{subfigure}[b]{0.5\textwidth}
\centering
\includegraphics{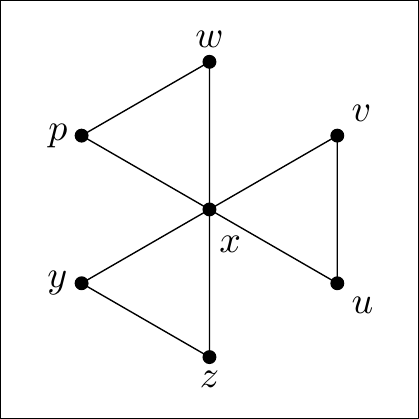}
%\begin{tikzpicture}[pont/.style={circle, fill=black, inner sep=0.5mm}]
%
%\draw (-2.121,-2.121) rectangle (2.121, 2.121);
%
%\node (Y) at (-1.299,-0.75) [pont] {};
%\node (Z) at (0,-1.5) [pont] {};
%\node (U) at (1.299,-0.75) [pont] {};
%\node (V) at (1.299,0.75) [pont] {};
%\node (W) at (0,1.5) [pont] {};
%\node (P) at (-1.299,0.75) [pont] {};
%\node (X) at (0,0) [pont] {};
%
%%\draw[thick, pattern=north west lines] (0,0) -- (-1.299,-0.75) -- (0,-1.5) -- cycle;
%%\draw[thick, pattern=vertical lines] (0,0) -- (1.299,-0.75) -- (1.299,0.75) -- cycle;
%%\draw[thick, pattern=north east lines] (0,0) -- (0,1.5) -- (-1.299,0.75) -- cycle;
%\draw (0,0) -- (-1.299,-0.75) -- (0,-1.5) -- cycle;
%\draw (0,0) -- (1.299,-0.75) -- (1.299,0.75) -- cycle;
%\draw (0,0) -- (0,1.5) -- (-1.299,0.75) -- cycle;
%
%\draw (Z) node[anchor=north] {\small$z$};
%\draw (Y) node[anchor=east] {\small$y$};
%\draw (0.21, -0.35) node {\small$x$};
%\draw (U) node[anchor=north west] {\small$u$};
%\draw (V) node[anchor=south west] {\small$v$};
%\draw (W) node[anchor=south] {\small$w$};
%\draw (P) node[anchor=east] {\small$p$};
%\end{tikzpicture}
\caption{Case 4}\label{SFthyp4}
\end{subfigure}%
\caption{Possible triangle hypergraphs of small exceptional quasilinear betweenness structures. The hyperedges are represented by triangles.}\label{Fthyp}
\end{figure}

Next, suppose that $5\leq n\leq 7$. Since for all $n\geq 3$, $$\tau(n, 0)= n - 4$$ by Part \ref{Equs11} Theorem \ref{Tqus1}, we obtain from Lemma \ref{Lgen2a} that $\mathcal{H}$ is a $\Delta$-star.
If $n = 5$ and $\mathcal{B}$ is regular or $n\geq 6$ and $\mathcal{H}$ is a tight star, then $\mathcal{B}$ is isomorphic to either $\mathcal{R}_{n, i}^4$ for some $i\in I_n^4$ or to $\mathcal{S}_n^4$, as shown by Lemma \ref{Lgen3} with parameter $c = 4$.
Therefore, either $n = 5$ and $\mathcal{B}$ is irregular or $n\geq 6$ and $\mathcal{H}$ is a non-tight $\Delta$-star. The four possible triangle hypergraphs are shown in Figure \ref{Fthyp}.
%%% regular -> acyclic, irregular -> cyclic???

\begin{figure}
\centering
\begin{subfigure}[b]{0.5\textwidth}
\centering
\includegraphics{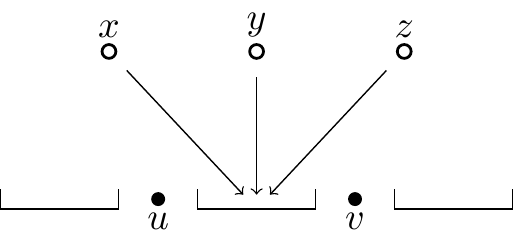}
%\begin{tikzpicture}[pont/.style={circle, fill=black, inner sep=0.5mm},
 %ppont/.style={circle, draw=black, thick, fill=white, minimum size = 4pt, inner sep=0pt, outer sep=5.5pt}]
%
%\node (U) at (0,0) [pont] {};
%\node (V) at (2,0) [pont] {};
%\node (X) at (-0.5,1.5) [ppont] {};
%\node (Y) at (1,1.5) [ppont] {};
%\node (Z) at (2.5,1.5) [ppont] {};
%\node (DEST) at (1,-0.1) {};
%
%\draw (-1.6,0.1)--(-1.6,-0.1)--(-0.4,-0.1)--(-0.4,0.1);
%\draw (0.4,0.1)--(0.4,-0.1)--(1.6,-0.1)--(1.6,0.1);
%\draw (2.4,0.1)--(2.4,-0.1)--(3.6,-0.1)--(3.6,0.1);
%
%\draw[->] (X)--(DEST);
%\draw[->] (Y)--(DEST);
%\draw[->] (Z)--(DEST);
%
%\draw (U) node[anchor=north] {\small$u$};
%\draw (V) node[anchor=north] {\small$v$};
%\draw (X) node[anchor=south] {\small$x$};
%\draw (Y) node[anchor=south] {\small$y$};
%\draw (Z) node[anchor=south] {\small$z$};
%
%\end{tikzpicture}
\caption{Case 1.1}
\end{subfigure}%
\begin{subfigure}[b]{0.5\textwidth}
\centering
\includegraphics{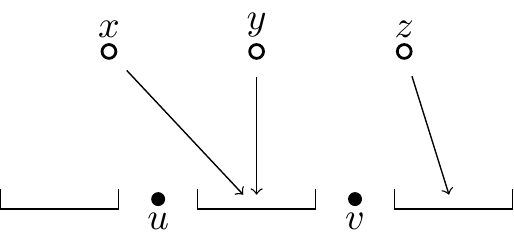}
%\begin{tikzpicture}[pont/.style={circle, fill=black, inner sep=0.5mm},
 %ppont/.style={circle, draw=black, thick, fill=white, minimum size = 4pt, inner sep=0pt, outer sep=5.5pt}]
%
%\node (U) at (0,0) [pont] {};
%\node (V) at (2,0) [pont] {};
%\node (X) at (-0.5,1.5) [ppont] {};
%\node (Y) at (1,1.5) [ppont] {};
%\node (Z) at (2.5,1.5) [ppont] {};
%\node (DEST1) at (1,-0.1) {};
%\node (DEST2) at (3,-0.1) {};
%
%\draw (-1.6,0.1)--(-1.6,-0.1)--(-0.4,-0.1)--(-0.4,0.1);
%\draw (0.4,0.1)--(0.4,-0.1)--(1.6,-0.1)--(1.6,0.1);
%\draw (2.4,0.1)--(2.4,-0.1)--(3.6,-0.1)--(3.6,0.1);
%
%\draw[->] (X)--(DEST1);
%\draw[->] (Y)--(DEST1);
%\draw[->] (Z)--(DEST2);
%
%\draw (U) node[anchor=north] {\small$u$};
%\draw (V) node[anchor=north] {\small$v$};
%\draw (X) node[anchor=south] {\small$x$};
%\draw (Y) node[anchor=south] {\small$y$};
%\draw (Z) node[anchor=south] {\small$z$};
%
%\end{tikzpicture}
\caption{Case 1.2}
\end{subfigure}%
\vspace{1em}
\begin{subfigure}[b]{0.5\textwidth}
\centering
\includegraphics{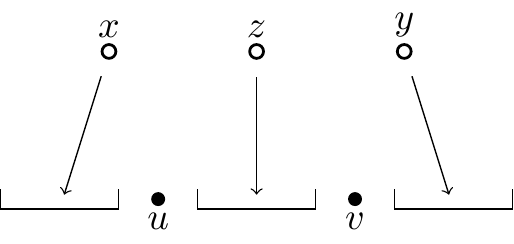}
%\begin{tikzpicture}[pont/.style={circle, fill=black, inner sep=0.5mm},
 %ppont/.style={circle, draw=black, thick, fill=white, minimum size = 4pt, inner sep=0pt, outer sep=5.5pt}]
%
%\node (U) at (0,0) [pont] {};
%\node (V) at (2,0) [pont] {};
%\node (X) at (-0.5,1.5) [ppont] {};
%\node (Z) at (1,1.5) [ppont] {};
%\node (Y) at (2.5,1.5) [ppont] {};
%\node (DEST1) at (-1,-0.1) {};
%\node (DEST2) at (1,-0.1) {};
%\node (DEST3) at (3,-0.1) {};
%
%\draw (-1.6,0.1)--(-1.6,-0.1)--(-0.4,-0.1)--(-0.4,0.1);
%\draw (0.4,0.1)--(0.4,-0.1)--(1.6,-0.1)--(1.6,0.1);
%\draw (2.4,0.1)--(2.4,-0.1)--(3.6,-0.1)--(3.6,0.1);
%
%\draw[->] (X)--(DEST1);
%\draw[->] (Y)--(DEST3);
%\draw[->] (Z)--(DEST2);
%
%\draw (U) node[anchor=north] {\small$u$};
%\draw (V) node[anchor=north] {\small$v$};
%\draw (X) node[anchor=south] {\small$x$};
%\draw (Y) node[anchor=south] {\small$y$};
%\draw (Z) node[anchor=south] {\small$z$};
%
%\end{tikzpicture}
\caption{Case 1.3}
\end{subfigure}%
\begin{subfigure}[b]{0.5\textwidth}
\centering
\includegraphics{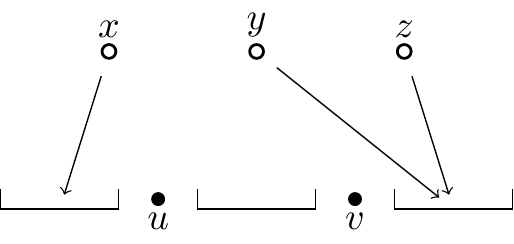}
%\begin{tikzpicture}[pont/.style={circle, fill=black, inner sep=0.5mm},
 %ppont/.style={circle, draw=black, thick, fill=white, minimum size = 4pt, inner sep=0pt, outer sep=5.5pt}]
%
%\node (U) at (0,0) [pont] {};
%\node (V) at (2,0) [pont] {};
%\node (X) at (-0.5,1.5) [ppont] {};
%\node (Y) at (1,1.5) [ppont] {};
%\node (Z) at (2.5,1.5) [ppont] {};
%\node (DEST1) at (-1,-0.1) {};
%\node (DEST2) at (3,-0.1) {};
%
%\draw (-1.6,0.1)--(-1.6,-0.1)--(-0.4,-0.1)--(-0.4,0.1);
%\draw (0.4,0.1)--(0.4,-0.1)--(1.6,-0.1)--(1.6,0.1);
%\draw (2.4,0.1)--(2.4,-0.1)--(3.6,-0.1)--(3.6,0.1);
%
%\draw[->] (X)--(DEST1);
%\draw[->] (Y)--(DEST2);
%\draw[->] (Z)--(DEST2);
%
%\draw (U) node[anchor=north] {\small$u$};
%\draw (V) node[anchor=north] {\small$v$};
%\draw (X) node[anchor=south] {\small$x$};
%\draw (Y) node[anchor=south] {\small$y$};
%\draw (Z) node[anchor=south] {\small$z$};
%
%\end{tikzpicture}
\caption{Case 1.4}
\end{subfigure}%
\caption{Cases 1.1-1.4. in the proof of Lemma \ref{Lsmallgr}}\label{Fc1ca-d}
\end{figure}

\begin{case}{1}
%I
In this case, $n = 5$, $X =\{u, v, x, y, z\}$ and $T =\{x, y, z\}$ is the sole triangle of $\mathcal{B}$ (see Figure \ref{SFthyp1}).
Because of the case's assumption, there exists a subset $Y\subseteq X$ such that \begin{equation}\label{EQ19}
\mathcal{B}\vert_Y\simeq\mathcal{C}_4. 
\end{equation} Points $u$ and $v$ partition $T$ into three parts: $(\cdot\ u\ v)_\mathcal{B} =\{p\in X : (p\ u\ v)_\mathcal{B}\}$, $(u\cdot v)_\mathcal{B} =\{p\in X : (u\ p\ v)_\mathcal{B}\}$ and $(u\ v\ \cdot)_\mathcal{B}  =\{p\in X : (u\ v\ p)_\mathcal{B}\}$. Notice that $T\not\subseteq Y$, therefore, $u, v\in Y$. We obtain from (\ref{EQ19}) that \begin{equation}\label{EQ20}
\emph{either }|(u\cdot v)_{\mathcal{B}\vert_Y}|\geq 2\emph{, or }|(\cdot\ u\ v)_{\mathcal{B}\vert_Y}|\geq 1\emph{ and }|(u\ v\ \cdot)_{\mathcal{B}\vert_Y}|\geq 1.
\end{equation} Note that we can replace $\mathcal{B}\vert_Y$ with $\mathcal{B}$ in (\ref{EQ20}), so the following cases are possible up to symmetry (see Figure \ref{Fc1ca-d}): %%% replace y with z in case C???
\begin{itemize}
%[label={\Alph*.}, ref={\Alph*}]
\item\textsc{Case} 1.1: $(u\ x\ v)_\mathcal{B}\wedge (u\ y\ v)_\mathcal{B}\wedge (u\ z\ v)_\mathcal{B}$;%\label{Ec1ca}
\item\textsc{Case} 1.2: $(u\ x\ v)_\mathcal{B}\wedge (u\ y\ v)_\mathcal{B}\wedge (u\ v\ z)_\mathcal{B}$;%\label{Ec1cb}
\item\textsc{Case} 1.3: $(x\ u\ v)_\mathcal{B}\wedge (u\ v\ y)_\mathcal{B}\wedge (u\ z\ v)_\mathcal{B}$;%\label{Ec1cc}
\item\textsc{Case} 1.4: $(x\ u\ v)_\mathcal{B}\wedge (u\ v\ y)_\mathcal{B}\wedge (u\ v\ z)_\mathcal{B}$.%\label{Ec1cd}
\end{itemize}

\begin{case}{1.1}
%A
We can suppose without loss of generality that $Y =\{u, v, x, y\}$. Now, because of (\ref{EQ19}), $(x\ u\ y)_\mathcal{B}$ and $(x\ v\ y)_\mathcal{B}$ hold. The betweenness structures $\mathcal{B} - x$ and $\mathcal{B} - y$ are linear, so they are isomorphic to either $\mathcal{P}_4$ or $\mathcal{C}_4$. We show that \begin{equation}\label{EQ5}
\mathcal{B} - x\simeq\mathcal{B} - y\simeq\mathcal{C}_4.
\end{equation}

Relation $(x\ z\ u)_\mathcal{B}$ cannot hold, otherwise $(x\ z\ u)_\mathcal{B}$ and $(x\ u\ y)_\mathcal{B}$ would imply $(x\ z\ y)_\mathcal{B}$ by f.r.p., in contradiction with $T$ being a triangle.
Similarly, we obtain that none of the betweennesses $(x\ z\ v)_\mathcal{B}$, $(y\ z\ u)_\mathcal{B}$ and $(y\ z\ v)_\mathcal{B}$ hold.
As a consequence, neither $\mathcal{B} - x$ nor $\mathcal{B} - y$ is isomorphic to $\mathcal{P}_4$. For example, if $\mathcal{B} - x\simeq\mathcal{P}_4$, then because of $(u\ y\ v)_\mathcal{B}$ and $(u\ z\ v)_\mathcal{B}$, $(u\ z\ y)_\mathcal{B}$ or $(y\ z\ v)_\mathcal{B}$ hold in contradiction with the previous assertion.

Now, (\ref{EQ5}) yields $(x\ u\ z)_\mathcal{B}$, $(x\ v\ z)_\mathcal{B}$, $(y\ u\ z)_\mathcal{B}$ and $(y\ v\ z)_\mathcal{B}$, therefore, we have $9$ nontrivial betweennesses: $$(u\ x\ v)_\mathcal{B}, (u\ y\ v)_\mathcal{B}, (u\ z\ v)_\mathcal{B}, (x\ u\ y)_\mathcal{B}, (x\ v\ y)_\mathcal{B}$$ $$(x\ u\ z)_\mathcal{B}, (x\ v\ z)_\mathcal{B}, (y\ u\ z)_\mathcal{B}\text{ and } (y\ v\ z)_\mathcal{B}.$$ These betweennesses cover all collinear triples of $\mathcal{B}$ and show that $\mathcal{B}\simeq\mathcal{B}(K_{2, 3}) =\mathcal{A}_5^1$.
\end{case}

\begin{case}{1.2}
%B
It follows from (\ref{EQ20}) that $Y =\{u, v, x, y\}$ and so $(x\ u\ y)_\mathcal{B}$ and $(x\ v\ y)_\mathcal{B}$ hold by (\ref{EQ19}). Relations $(u\ x\ z)_\mathcal{B}$ and $(x\ v\ z)_\mathcal{B}$ follows from $(u\ x\ v)_\mathcal{B}$ and $(u\ v\ z)_\mathcal{B}$ by f.r.p. Similarly, $(u\ y\ z)_\mathcal{B}$ and $(y\ v\ z)_\mathcal{B}$ follows from $(u\ y\ v)_\mathcal{B}$ and $(u\ v\ z)_\mathcal{B}$, hence, considering (\ref{EQ19}) and the case's initial assumptions, we have $9$ non-trivial betweennesses: $$(x\ u\ y)_\mathcal{B}, (u\ y\ v)_\mathcal{B}, (y\ v\ x)_\mathcal{B}, (v\ x\ u)_\mathcal{B}, (u\ v\ z)_\mathcal{B},$$ $$(u\ x\ z)_\mathcal{B}, (x\ v\ z)_\mathcal{B}, (u\ x\ v)_\mathcal{B}\text{ and } (u\ v\ z)_\mathcal{B}.$$ These betweennesses cover all collinear triples of $\mathcal{B}$ and show that $\mathcal{B}\simeq\mathcal{R}_{5, 1}^4$.
\end{case}

\begin{case}{1.3}
%C
It follows from (\ref{EQ20}) that $Y =\{u, v, x, y\}$ and so $(v\ y\ x)_\mathcal{B}$ and $(y\ x\ u)_\mathcal{B}$ hold by (\ref{EQ19}).
Further, $(x\ u\ v)_\mathcal{B}$ and $(u\ z\ v)_\mathcal{B}$ imply $(x\ u\ z)_\mathcal{B}$ by f.r.p. Similarly, $(u\ v\ y)_\mathcal{B}$ and $(u\ z\ v)_\mathcal{B}$ imply $(u\ z\ y)_\mathcal{B}$ and $(z\ v\ y)_\mathcal{B}$.

Now, considering the case's assumptions, we have $9$ non-trivial betweennesses: $$(x\ u\ v)_\mathcal{B}, (u\ v\ y)_\mathcal{B}, (u\ z\ v)_\mathcal{B}, (v\ y\ x)_\mathcal{B}, (y\ x\ u)_\mathcal{B}$$ $$(x\ u\ z)_\mathcal{B}, (x\ z\ v)_\mathcal{B}, (u\ z\ y)_\mathcal{B}\text{ and } (z\ v\ y)_\mathcal{B}.$$  These betweennesses cover all collinear triples of $\mathcal{B}$ and show that $\mathcal{B}\simeq\mathcal{S}_5^4$.
\end{case}

\begin{case}{1.4}
It follows from (\ref{EQ20}) that $\mathcal{B} - x\not\simeq\mathcal{C}_4$, so we can assume without loss of generality that $Y =\{u, v, x, y\}$ and hence $(y\ x\ u)_\mathcal{B}$ and $(v\ y\ x)_\mathcal{B}$ hold by (\ref{EQ19}). Moreover, since $\mathcal{B} - x\simeq\mathcal{P}_4$, $(u\ v\ y\ z)_\mathcal{B}$ or $(u\ v\ z\ y)_\mathcal{B}$ is true. However, in case of $(u\ v\ z\ y)_\mathcal{B}$, $(v\ z\ y)_\mathcal{B}$ and $(v\ y\ x)_\mathcal{B}$ would lead to the contradiction $(z\ y\ x)_\mathcal{B}$ by f.r.p., hence, $(u\ v\ y\ z)_\mathcal{B}$ holds and in particular, \begin{equation}\label{EQ5a}
(v\ y\ z)_\mathcal{B}
\end{equation} must be true. Now, there are two cases to consider.

\begin{wcase}{1.4.1}{$\mathcal{B} - y\simeq\mathcal{C}_4$}
%i
In this case, relations $(z\ x\ u)_\mathcal{B}$ and $(v\ z\ x)_\mathcal{B}$ hold. However, $(v\ z\ x)_\mathcal{B}$ and (\ref{EQ5a}) imply $(y\ z\ x)_\mathcal{B}$ by f.r.p., which is impossible for $T$ is a triangle.
\end{wcase}

\begin{wcase}{1.4.2}{$\mathcal{B} - y\simeq\mathcal{P}_4$}
%ii
Now, $(x\ u\ v\ z)_\mathcal{B}$ holds. However, $(x\ v\ z)_\mathcal{B}$ and (\ref{EQ5a}) leads to a contradiction again.
\end{wcase}
\end{case}

\begin{figure}
\centering
\includegraphics{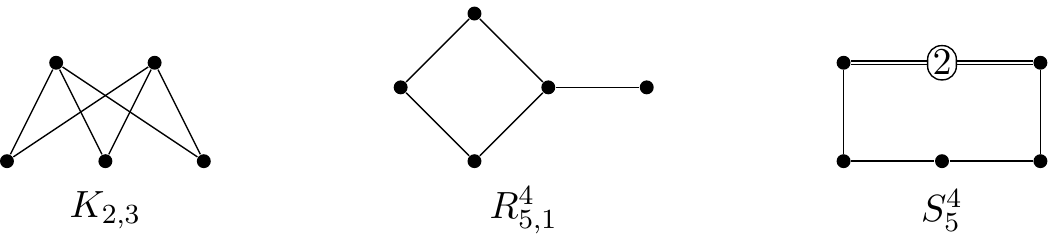}
%\begin{tikzpicture}[pont/.style={circle, fill=black, inner sep=0.5mm},
%mybox/.style={rectangle, rounded corners, draw=black, fill=white, inner sep=0.5mm}]
%
%\begin{scope}[shift={(0,0)}]
%\node (A) at (0,0) [pont] {};
%\node (B) at (1,0) [pont] {};
%\node (C) at (2,0) [pont] {};
%\node (D) at (0.5,1) [pont] {};
%\node (E) at (1.5,1) [pont] {};
%
%\draw (A)--(D) node {};
%\draw (B)--(D) node {};
%\draw (C)--(D) node {};
%\draw (A)--(E) node {};
%\draw (B)--(E) node {};
%\draw (C)--(E) node {};
%
%\draw (1,-0.5) node {\small$K_{2, 3}$};
%\end{scope}
%
%\begin{scope}[shift={(4,0.75)}]
%\node (A) at (0,0) [pont] {};
%\node (B) at (1.5,0) [pont] {};
%\node (C) at (2.5,0) [pont] {};
%\node (D) at (0.75,0.75) [pont] {};
%\node (E) at (0.75,-0.75) [pont] {};
%
%\draw (A)--(D) node {};
%\draw (A)--(E) node {};
%\draw (D)--(B) node {};
%\draw (E)--(B) node {};
%\draw (B)--(C) node {};
%
%\draw (1.25,-1.25) node {\small$R_{5, 1}^4$};
%\end{scope}
%
%\begin{scope}[shift={(8.5,0)}]
%\node (A) at (0,1) [pont] {};
%\node (B) at (0,0) [pont] {};
%\node (C) at (1,0) [pont] {};
%\node (D) at (2,0) [pont] {};
%\node (E) at (2,1) [pont] {};
%
%\draw (A)--(B) node {};
%\draw (B)--(C) node {};
%\draw (C)--(D) node {};
%\draw (D)--(E) node {};
%\draw[double] (E)--(A) node[midway, mybox] {\small$2$};
%%\draw (E)--(A) node {};
%
%\draw (1,-0.5) node {\small$S_5^4$};
%\end{scope}
%\end{tikzpicture}
\caption{Spanner graphs of $5$-point quasilinear betweenness structures}\label{Fqlinord5}
\end{figure}

In summary, we obtain that the quasilinear betweenness structures on $5$ points are exactly the ones induced by the graphs listed in Figure \ref{Fqlinord5}.
\end{case}

\begin{case}{2}
%II
In this case, $n = 6$, $X =\{u, v, w, x, y, z\}$ and there are two triangles in $\mathcal{B}$ that intersect in one vertex: $T =\{x\ y\ z\}$ and $R =\{x\ u\ v\}$ (see Figure \ref{SFthyp2}). Now, $\mathcal{B} - x$ is a linear betweenness structure of order $5$, therefore, it is ordered by Proposition \ref{Plin}. Further, $\mathcal{B} - y$, $\mathcal{B} - z$, $\mathcal{B} - u$ and $\mathcal{B} - v$ are quasilinear betweenness structures of order $5$, hence, they are induced by one of the $5$-vertex graphs in Figure \ref{Fqlinord5}. However, the spanner graph cannot be isomorphic to $K_{2, 3}$ because each substructure in question intersect the ordered substructure $\mathcal{B} - x$ in $4$ points and thus has a $4$-point ordered substructure itself.

\begin{figure}
\centering
\begin{subfigure}[b]{0.5\textwidth}
\centering
\includegraphics{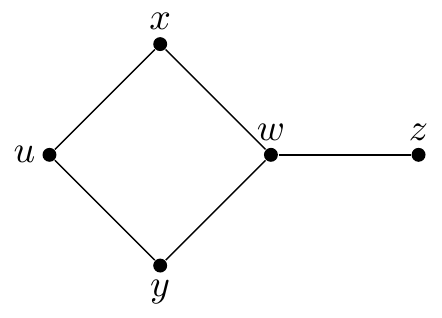}
%\begin{tikzpicture}[scale=1.5, pont/.style={circle, fill=black, inner sep=0.5mm}]
%\node (A) at (0,0) [pont] {};
%\node (B) at (1.5,0) [pont] {};
%\node (C) at (2.5,0) [pont] {};
%\node (D) at (0.75,0.75) [pont] {};
%\node (E) at (0.75,-0.75) [pont] {};
%
%\draw (A)--(D) node {};
%\draw (A)--(E) node {};
%\draw (D)--(B) node {};
%\draw (E)--(B) node {};
%\draw (B)--(C) node {};
%
%\draw (A) node[anchor=east] {\small$u$};
%\draw (B) node[anchor=south] {\small$w$};
%\draw (C) node[anchor=south] {\small$z$};
%\draw (D) node[anchor=south] {\small$x$};
%\draw (E) node[anchor=north] {\small$y$};
%\end{tikzpicture}
\caption{Case 2.1}
\end{subfigure}%
\begin{subfigure}[b]{0.5\textwidth}
\centering
\includegraphics{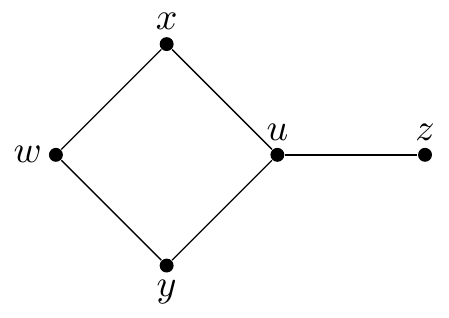}
%\begin{tikzpicture}[scale=1.5, pont/.style={circle, fill=black, inner sep=0.5mm}]
%\node (A) at (0,0) [pont] {};
%\node (B) at (1.5,0) [pont] {};
%\node (C) at (2.5,0) [pont] {};
%\node (D) at (0.75,0.75) [pont] {};
%\node (E) at (0.75,-0.75) [pont] {};
%
%\draw (A)--(D) node {};
%\draw (A)--(E) node {};
%\draw (D)--(B) node {};
%\draw (E)--(B) node {};
%\draw (B)--(C) node {};
%
%\draw (A) node[anchor=east] {\small$w$};
%\draw (B) node[anchor=south] {\small$u$};
%\draw (C) node[anchor=south] {\small$z$};
%\draw (D) node[anchor=south] {\small$x$};
%\draw (E) node[anchor=north] {\small$y$};
%\end{tikzpicture}
\caption{Case 2.2}
\end{subfigure}%
\vspace{1em}
\begin{subfigure}[b]{0.5\textwidth}
\centering
\includegraphics{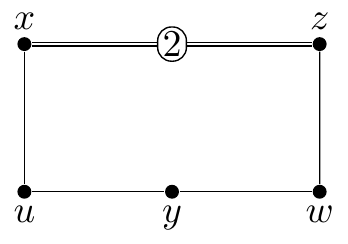}
%\begin{tikzpicture}[scale=1.5, pont/.style={circle, fill=black, inner sep=0.5mm},
%mybox/.style={rectangle, rounded corners, draw=black, fill=white, inner sep=0.5mm}]
%\node (A) at (0,1) [pont] {};
%\node (B) at (0,0) [pont] {};
%\node (C) at (1,0) [pont] {};
%\node (D) at (2,0) [pont] {};
%\node (E) at (2,1) [pont] {};
%
%\draw (A)--(B) node {};
%\draw (B)--(C) node {};
%\draw (C)--(D) node {};
%\draw (D)--(E) node {};
%\draw[double] (E)--(A) node[midway, mybox] {\small$2$};
%
%\draw (A) node[anchor=south] {\small$x$};
%\draw (B) node[anchor=north] {\small$u$};
%\draw (C) node[anchor=north] {\small$y$};
%\draw (D) node[anchor=north] {\small$w$};
%\draw (E) node[anchor=south] {\small$z$};
%\end{tikzpicture}
\caption{Case 2.3}
\end{subfigure}%
\begin{subfigure}[b]{0.5\textwidth}
\centering
\includegraphics{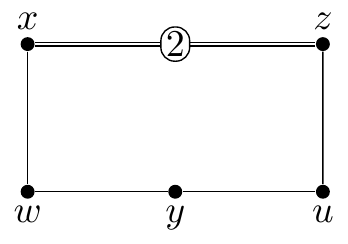}
%\begin{tikzpicture}[scale=1.5, pont/.style={circle, fill=black, inner sep=0.5mm},
%mybox/.style={rectangle, rounded corners, draw=black, fill=white, inner sep=0.5mm}]
%\node (A) at (0,1) [pont] {};
%\node (B) at (0,0) [pont] {};
%\node (C) at (1,0) [pont] {};
%\node (D) at (2,0) [pont] {};
%\node (E) at (2,1) [pont] {};
%
%\draw (A)--(B) node {};
%\draw (B)--(C) node {};
%\draw (C)--(D) node {};
%\draw (D)--(E) node {};
%\draw[double] (E)--(A) node[midway, mybox] {\small$2$};
%
%\draw (A) node[anchor=south] {\small$x$};
%\draw (B) node[anchor=north] {\small$w$};
%\draw (C) node[anchor=north] {\small$y$};
%\draw (D) node[anchor=north] {\small$u$};
%\draw (E) node[anchor=south] {\small$z$};
%\end{tikzpicture}
\caption{Case 2.4}
\end{subfigure}%
\caption{Cases 2.1-2.4. in the proof of Lemma \ref{Lsmallgr}}\label{Fc2}
\end{figure}

Consider now $\mathcal{B} - v$. Because of the argument above, it is induced by either $R_{5, 1}^4$ or $S_5^4$. As $(\mathcal{B} - v) - x$ must be ordered, we can determine the position of $x$ in the spanner graph of $\mathcal{B} - v$ up to symmetry. After that, $y$ and $z$ must be the two uniquely determined vertices for which $\{x, y, z\}$ forms a triangle ($y$ and $z$ are in symmetric position in $\mathcal{B}$, so they are interchangeable). Finally, we can place $u$ and $w$ in two possible ways. There are four cases altogether shown in Figure \ref{Fc2}.

\begin{case}{2.1}
%A
The ordering of $\mathcal{B} - x$ is an extension of the ordering of $(\mathcal{B} - x) - v = [u, y, w, z]$. Let $i$ denote the position of $v$ in this extended ordering.

\begin{wcase}{2.1.1}{$i = 1$}
$(v\ u\ x)_{\mathcal{B} - x}\wedge (u\ x\ w)_{\mathcal{B} - v}\overset{\text{f.r.p.}}{\Rightarrow} (v\ u\ x)_\mathcal{B}$, which is impossible since $R =\{x, u, v\}$ was a triangle.
%Relations $(v\ u\ w)_{\mathcal{B} - x} $ and $(u\ x\ w)_{\mathcal{B} - v}$ imply $(v\ u\ x)_\mathcal{B}$ by f.r.p., which is impossible since $R =\{x, u, v\}$ was a triangle.
\end{wcase}

\begin{wcase}{2.1.2}{$i = 2$}
$(x\ u\ y)_{\mathcal{B} - v}\wedge (u\ v\ y)_{\mathcal{B} - x}\overset{\text{f.r.p.}}{\Rightarrow} (x\ u\ v)_\mathcal{B}\ \lightning$
%Relations $(x\ u\ y)_{\mathcal{B} - v}$ and $(u\ v\ y)_{\mathcal{B} - x}$ imply $(x\ u\ v)_\mathcal{B}$ by f.r.p., a contradiction again.
\end{wcase}

\begin{wcase}{2.1.3}{$i = 3$}
As the triple $\{v, x, z\}$ is collinear, so one of the following cases hold:
\begin{itemize}
\item\textsc{Case} 2.1.3.1: $(v\ x\ z)_\mathcal{B}$;
\item\textsc{Case} 2.1.3.2: $(x\ v\ z)_\mathcal{B}$;
\item\textsc{Case} 2.1.3.3: $(x\ z\ v)_\mathcal{B}$.
\end{itemize}

\begin{case}{2.1.3.1}
%a
$(v\ x\ z)_\mathcal{B}\wedge (u\ v\ z)_{\mathcal{B} - x}\overset{\text{f.r.p.}}{\Rightarrow} (u\ v\ x)_\mathcal{B}\ \lightning$
%Relations $(v\ x\ z)_\mathcal{B}$ and $(u\ v\ z)_{\mathcal{B} - x}$ imply $(u\ v\ x)_\mathcal{B}$ by f.r.p., a contradiction.
\end{case}

\begin{case}{2.1.3.2}
%b
$(x\ v\ z)_\mathcal{B}\wedge (u\ x\ z)_{\mathcal{B} - v}\overset{\text{f.r.p.}}{\Rightarrow} (u\ x\ v)_\mathcal{B}\ \lightning$
%Relations $(x\ v\ z)_\mathcal{B}$ and $(u\ x\ z)_{\mathcal{B} - v}$ imply $(u\ x\ v)_\mathcal{B}$ by f.r.p., a contradiction.
%%contradicting $(x\ y\ w)_\mathcal{B}$.
\end{case}

\begin{case}{2.1.3.3}
%c
$(x\ z\ v)_\mathcal{B}\wedge (x\ w\ z)_{\mathcal{B} - v}\overset{\text{f.r.p.}}{\Rightarrow} (w\ z\ v)_\mathcal{B}$, contradicting $(v\ w\ z)_{\mathcal{B} - x}$.
%Relations $(x\ z\ v)_\mathcal{B}$ and $(x\ w\ z)_{\mathcal{B} - v}$ imply $(w\ z\ v)_\mathcal{B}$ by f.r.p., contradicting $(v\ w\ z)_{\mathcal{B} - x}$.
\end{case}
\end{wcase}
%Relations $(x\ u\ w)_\mathcal{B}$ and $(u\ v\ w)_\mathcal{B}$ imply $(x\ u\ v)_\mathcal{B}$ by f.r.p., a contradiction again.
%\end{wcase}

\begin{wcase}{2.1.4}{$i = 4$}
$(u\ w\ v)_{\mathcal{B} - x}\wedge (u\ x\ w)_{\mathcal{B} - v}\overset{\text{f.r.p.}}{\Rightarrow} (u\ x\ v)_\mathcal{B}\ \lightning$
%Relations $(u\ w\ v)_{\mathcal{B} - x}$ and $(u\ x\ w)_{\mathcal{B} - v}$ imply $(u\ x\ v)_\mathcal{B}$ by f.r.p., a contradiction.
\end{wcase}
%As $v, x$ and $y$ are collinear, one of the following cases hold:
%\begin{itemize}
%\item\textsc{Case} 2.1.4.1: $(v\ x\ y)_\mathcal{B}$;
%\item\textsc{Case} 2.1.4.2: $(x\ v\ y)_\mathcal{B}$;
%\item\textsc{Case} 2.1.4.2: $(x\ y\ v)_\mathcal{B}$.
%\end{itemize}
%
%\begin{case}{2.1.4.1}%a
%Relations $(v\ x\ y)_\mathcal{B}$ and $(x\ u\ y)_\mathcal{B}$ imply $(v\ x\ u)_\mathcal{B}$ by f.r.p., a contradiction.
%\end{case}
%
%\begin{case}{2.1.4.2}%b
%Relations $(x\ v\ y)_\mathcal{B}$ and $(y\ w\ v)_\mathcal{B}$ imply $(x\ w\ y)_\mathcal{B}$ by f.r.p., contradicting $(x\ y\ w)_\mathcal{B}$.
%\end{case}
%
%\begin{case}{2.1.4.3}%c
%Relations $(x\ y\ v)_\mathcal{B}$ and $(x\ u\ y)_\mathcal{B}$ imply $(x\ u\ v)_\mathcal{B}$ by f.r.p., a contradiction.
%\end{case}
%\end{wcase}

\begin{wcase}{2.1.5}{$i = 5$}
$(u\ w\ v)_{\mathcal{B} - x}\wedge (u\ x\ w)_{\mathcal{B} - v}\overset{\text{f.r.p.}}{\Rightarrow} (u\ x\ v)_\mathcal{B}\ \lightning$
%Relations $(u\ w\ v)_{\mathcal{B} - x}$ and $(u\ x\ w)_{\mathcal{B} - v}$ imply $(u\ x\ v)_\mathcal{B}$ by f.r.p., a contradiction.
\end{wcase}
\end{case}

\begin{case}{2.2}
%B
Similarly to the previous case, $\mathcal{B} - x$ is an extension of the ordering of $(\mathcal{B} - x) - v = [w, y, u, z]$. Let $i$ denote the position of $v$ in this extended order.

\begin{wcase}{2.2.1}{$i = 1$}
$(v\ w\ u)_{\mathcal{B} - x}\wedge (w\ x\ u)_{\mathcal{B} - v}\overset{\text{f.r.p.}}{\Rightarrow} (v\ x\ u)_\mathcal{B}\ \lightning$
%Relations $(v\ w\ u)_{\mathcal{B} - x}$ and $(w\ x\ u)_{\mathcal{B} -v}$ imply $(v\ x\ u)_\mathcal{B}$ by f.r.p., a contradiction.
%As $v, x$ and $y$ are collinear, one of the following relations hold:
%\begin{itemize}%[label={\alph*.}, ref={\alph*}]
%\item\textsc{Case} 2.2.1.1: $(v\ x\ y)_\mathcal{B}$;
%\item\textsc{Case} 2.2.1.2: $(x\ v\ y)_\mathcal{B}$;
%\item\textsc{Case} 2.2.1.3: $(x\ y\ v)_\mathcal{B}$.
%\end{itemize}
%
%\begin{case}{2.2.1.1}%a
%Relations $(v\ x\ y)_\mathcal{B}$ and $(v\ y\ u)_\mathcal{B}$ imply $(v\ x\ u)_\mathcal{B}$ by f.r.p., a contradiction.
%\end{case}
%
%\begin{case}{2.2.1.2}%b
%Relations $(x\ v\ y)_\mathcal{B}$ and $(x\ y\ u)_\mathcal{B}$ imply $(x\ v\ u)_\mathcal{B}$ by f.r.p. a contradiction.
%\end{case}
%
%\begin{case}{2.2.1.3}%c
%Relations $(x\ y\ v)_\mathcal{B}$ and $(x\ w\ y)_\mathcal{B}$ imply $(w\ y\ v)_\mathcal{B}$ by f.r.p., contradicting $(v\ w\ y)_\mathcal{B}$.
%\end{case}
\end{wcase}

\begin{wcase}{2.2.2}{$i = 2$}
As the triple $\{v, x, z\}$ is collinear, one of the following cases hold:
\begin{itemize}
%[label={\alph*.}, ref={\alph*}]
\item\textsc{Case} 2.2.2.1: $(v\ x\ z)_\mathcal{B}$;
\item\textsc{Case} 2.2.2.2: $(x\ v\ z)_\mathcal{B}$;
\item\textsc{Case} 2.2.2.3: $(x\ z\ v)_\mathcal{B}$.
\end{itemize}

\begin{case}{2.2.2.1}
%a
$(v\ x\ z)_\mathcal{B}\wedge (x\ u\ z)_{\mathcal{B} - v}\overset{\text{f.r.p.}}{\Rightarrow} (v\ x\ u)_\mathcal{B}\ \lightning$
%Relations $(v\ x\ z)_\mathcal{B}$ and $(x\ u\ z)_{\mathcal{B} - v}$ imply $(v\ x\ u)_\mathcal{B}$ by f.r.p., a contradiction.
\end{case}

\begin{case}{2.2.2.2}
%b
$(x\ v\ z)_\mathcal{B}\wedge (v\ u\ z)_{\mathcal{B} - x}\overset{\text{f.r.p.}}{\Rightarrow} (x\ v\ u)_\mathcal{B}\ \lightning$
%Relations $(x\ v\ z)_\mathcal{B}$ and $(v\ u\ z)_{\mathcal{B} - x}$ imply $(x\ v\ u)_\mathcal{B}$ by f.r.p., a contradiction.
\end{case}

\begin{case}{2.2.2.3}
%c
$(x\ z\ v)_\mathcal{B}\wedge (x\ u\ z)_{\mathcal{B} - v}\overset{\text{f.r.p.}}{\Rightarrow} (x\ u\ v)_\mathcal{B}\ \lightning$
%Relations $(x\ z\ v)_\mathcal{B}$ and $(x\ u\ z)_{\mathcal{B} - v}$ imply $(x\ u\ v)_\mathcal{B}$ by f.r.p., a contradiction.
\end{case}
%Relations $(w\ v\ u)_\mathcal{B}$ and $(x\ w\ u)_\mathcal{B}$ imply $(x\ v\ u)_\mathcal{B}$ by f.r.p., a contradiction.
\end{wcase}

\begin{wcase}{2.2.3}{$i = 3$}
$(y\ u\ x)_{\mathcal{B} - v}\wedge (y\ v\ u)_{\mathcal{B} - x}\overset{\text{f.r.p.}}{\Rightarrow} (v\ u\ x)_\mathcal{B}\ \lightning$
%Relations $(y\ u\ x)_{\mathcal{B} - v}$ and $(y\ v\ u)_{\mathcal{B} - x}$ imply $(v\ u\ x)_\mathcal{B}$ by f.r.p., a contradiction.
\end{wcase}

\begin{wcase}{2.2.4}{$i = 4$}
$(x\ u\ z)_{\mathcal{B} - v}\wedge (u\ v\ z)_{\mathcal{B} - x}\overset{\text{f.r.p.}}{\Rightarrow} (x\ u\ v)_\mathcal{B}\ \lightning$
%Relations $(x\ u\ z)_{\mathcal{B} - v}$ and $(u\ v\ z)_{\mathcal{B} - x}$ imply $(x\ u\ v)_\mathcal{B}$ by f.r.p., a contradiction.
\end{wcase}

\begin{wcase}{2.2.5}{$i = 5$}
$(w\ u\ v)_{\mathcal{B} - x}\wedge (w\ x\ u)_{\mathcal{B} - v}\overset{\text{f.r.p.}}{\Rightarrow} (x\ u\ v)_\mathcal{B}\ \lightning$
%Relations $(w\ u\ v)_{\mathcal{B} - x}$ and $(w\ x\ u)_{\mathcal{B} - v}$ imply $(x\ u\ v)_\mathcal{B}$ by f.r.p., a contradiction.
%As $v, x$ and $y$ are collinear, one of the following relations hold:
%\begin{itemize}%[label={\alph*.}, ref={\alph*}]
%\item\textsc{Case} 2.2.5.1: $(v\ x\ y)_\mathcal{B}$;
%\item\textsc{Case} 2.2.5.2: $(x\ v\ y)_\mathcal{B}$;
%\item\textsc{Case} 2.2.5.3: $(x\ y\ v)_\mathcal{B}$.
%\end{itemize}
%
%\begin{case}{2.2.5.1}%a
%Relations $(v\ x\ y)_\mathcal{B}$ and $(x\ w\ y)_\mathcal{B}$ imply $(v\ w\ y)_\mathcal{B}$ by f.r.p., contradicting $(w\ y\ v)_\mathcal{B}$.
%\end{case}
%
%\begin{case}{2.2.5.2}%b
%Relations $(x\ v\ y)_\mathcal{B}$ and $(v\ u\ y)_\mathcal{B}$ imply $(x\ v\ u)_\mathcal{B}$ by f.r.p., a contradiction.
%\end{case}
%
%\begin{case}{2.2.5.3}%c
%Relations $(x\ y\ v)_\mathcal{B}$ and $(y\ u\ v)_\mathcal{B}$ imply $(x\ u\ v)_\mathcal{B}$ by f.r.p., a contradiction.
%\end{case}
\end{wcase}
\end{case}

\begin{case}{2.3}
%C
Similarly to the previous cases, $\mathcal{B} - x$ is an extension of the ordering of $(\mathcal{B} - x) - v = [u, y, w, z]$. Let $i$ denote the position of $v$ in this extended order.

\begin{wcase}{2.3.1}{$i = 1$}
$(v\ u\ z)_{\mathcal{B} - x}\wedge (u\ x\ z)_{\mathcal{B} - v}\overset{\text{f.r.p.}}{\Rightarrow} (v\ u\ x)_\mathcal{B}\ \lightning$
%Relations $(v\ u\ z)_{\mathcal{B} - x}$ and $(u\ x\ z)_{\mathcal{B} - v}$ imply $(v\ u\ x)_\mathcal{B}$ by f.r.p., a contradiction.
\end{wcase}

\begin{wcase}{2.3.2}{$i = 2$}
$(x\ u\ y)_{\mathcal{B} - v}\wedge (u\ v\ y)_{\mathcal{B} - x}\overset{\text{f.r.p.}}{\Rightarrow} (x\ u\ v)_\mathcal{B}\ \lightning$
%Relations $(x\ u\ y)_{\mathcal{B} - v}$ and $(u\ v\ y)_{\mathcal{B} - x}$  imply $(x\ u\ v)_\mathcal{B}$ by f.r.p., a contradiction.
\end{wcase}

\begin{wcase}{2.3.3}{$i = 3$}
$(x\ u\ w)_{\mathcal{B} - v}\wedge (u\ v\ w)_{\mathcal{B} - x}\overset{\text{f.r.p.}}{\Rightarrow} (x\ u\ v)_\mathcal{B}\ \lightning$
%Relations $(x\ u\ w)_{\mathcal{B} - v}$ and $(u\ v\ w)_{\mathcal{B} - x}$ imply $(x\ u\ v)_\mathcal{B}$ by f.r.p., a contradiction.
%As $v, x$ and $z$ are collinear, one of the following relations hold:
%\begin{itemize}%[label={\alph*.}, ref={\alph*}]
%\item\textsc{Case} 2.3.3.1: $(v\ x\ z)_\mathcal{B}$;
%\item\textsc{Case} 2.3.3.2: $(x\ v\ z)_\mathcal{B}$;
%\item\textsc{Case} 2.3.3.3: $(x\ z\ v)_\mathcal{B}$.
%\end{itemize}
%
%\begin{case}{2.3.3.1}%a
%Relations $(v\ x\ z)_\mathcal{B}$ and $(y\ v\ z)_\mathcal{B}$ imply $(y\ x\ z)_\mathcal{B}$ by f.r.p., a contradiction.
%\end{case}
%
%\begin{case}{2.3.3.2}%b
%Relations $(x\ v\ z)_\mathcal{B}$ and $(u\ x\ z)_\mathcal{B}$ imply $(u\ x\ v)_\mathcal{B}$ by f.r.p., a contradiction.
%\end{case}
%
%\begin{case}{2.3.3.3}%c
%Relations $(x\ z\ v)_\mathcal{B}$ and $(x\ w\ z)_\mathcal{B}$ imply $(w\ z\ v)_\mathcal{B}$ by f.r.p., contradicting $(v\ w\ z)_\mathcal{B}$.
%\end{case}
\end{wcase}

\begin{wcase}{2.3.4}{$i = 4$}
As the triple $\{v, x, y\}$ is collinear, one of the following cases hold:
\begin{itemize}
%[label={\alph*.}, ref={\alph*}]
\item\textsc{Case} 2.3.4.1: $(v\ x\ y)_\mathcal{B}$;
\item\textsc{Case} 2.3.4.2: $(x\ v\ y)_\mathcal{B}$;
\item\textsc{Case} 2.3.4.3: $(x\ y\ v)_\mathcal{B}$.
\end{itemize}

\begin{case}{2.3.4.1}
%a
$(v\ x\ y)_\mathcal{B}\wedge (x\ u\ y)_{\mathcal{B} - v}\overset{\text{f.r.p.}}{\Rightarrow} (v\ x\ u)_\mathcal{B}\ \lightning$
%Relations $(v\ x\ y)_\mathcal{B}$ and $(x\ u\ y)_{\mathcal{B} - v}$ imply $(v\ x\ u)_\mathcal{B}$ by f.r.p., a contradiction.
\end{case}

\begin{case}{2.3.4.2}
%b
$(x\ v\ y)_\mathcal{B}\wedge (x\ y\ w)_{\mathcal{B} - v}\overset{\text{f.r.p.}}{\Rightarrow} (v\ y\ w)_\mathcal{B}$, contradicting $(y\ w\ v)_{\mathcal{B} - x}$.
%Relations $(x\ v\ y)_\mathcal{B}$ and $(x\ y\ w)_{\mathcal{B} - v}$ imply $(v\ y\ w)_\mathcal{B}$ by f.r.p., contradicting $(y\ w\ v)_{\mathcal{B} - x}$.
\end{case}

\begin{case}{2.3.4.3}
%c
$(x\ y\ v)_\mathcal{B}\wedge (x\ u\ y)_{\mathcal{B} - v}\overset{\text{f.r.p.}}{\Rightarrow} (x\ u\ v)_\mathcal{B}\ \lightning$
%Relations $(x\ y\ v)_\mathcal{B}$ and $(x\ u\ y)_{\mathcal{B} - v}$ imply $(x\ u\ v)_\mathcal{B}$ by f.r.p., a contradiction.
%%contradicting $(v\ w\ z)_\mathcal{B}$.
\end{case}
%Relations $(u\ x\ w)_\mathcal{B}$ and $(u\ w\ v)_\mathcal{B}$ imply $(u\ x\ v)_\mathcal{B}$ by f.r.p., a contradiction.
\end{wcase}

\begin{wcase}{2.3.5}{$i = 5$}
$(u\ z\ v)_{\mathcal{B} - x}\wedge (u\ x\ z)_{\mathcal{B} - v}\overset{\text{f.r.p.}}{\Rightarrow} (u\ x\ v)_\mathcal{B}\ \lightning$
%Relations $(u\ z\ v)_{\mathcal{B} - x}$ and $(u\ x\ z)_{\mathcal{B} - v}$ imply $(u\ x\ v)_\mathcal{B}$ by f.r.p., a contradiction.
\end{wcase}
\end{case}

\begin{case}{2.4}
%D
Similarly to the previous cases, $\mathcal{B} - x$ is an extension of the ordering of $(\mathcal{B} - x) - v = [w, y, u, z]$. Let $i$ denote the position of $v$ in this extended order.

\begin{wcase}{2.4.1}{$i = 1$}
As the triple $\{v, x, y\}$ is collinear, one of the following cases hold:
\begin{itemize}
%[label={\alph*.}, ref={\alph*}]
\item\textsc{Case} 2.4.1.1: $(v\ x\ y)_\mathcal{B}$;
\item\textsc{Case} 2.4.1.2: $(x\ v\ y)_\mathcal{B}$;
\item\textsc{Case} 2.4.1.3: $(x\ y\ v)_\mathcal{B}$.
\end{itemize}

\begin{case}{2.4.1.1}
%a
$(v\ x\ y)_\mathcal{B}\wedge (v\ y\ u)_{\mathcal{B} - x}\overset{\text{f.r.p.}}{\Rightarrow} (v\ x\ u)_\mathcal{B}\ \lightning$
%Relations $(v\ x\ y)_\mathcal{B}$ and $(v\ y\ u)_{\mathcal{B} - x}$ imply $(v\ x\ u)_\mathcal{B}$ by f.r.p., a contradiction.
\end{case}

\begin{case}{2.4.1.2}
%b
$(x\ v\ y)_\mathcal{B}\wedge (x\ y\ u)_{\mathcal{B} - v}\overset{\text{f.r.p.}}{\Rightarrow} (x\ v\ u)_\mathcal{B}\ \lightning$
%Relations $(x\ v\ y)_\mathcal{B}$ and $(x\ y\ u)_{\mathcal{B} - v}$ imply $(x\ v\ u)_\mathcal{B}$ by f.r.p., a contradiction.
\end{case}

\begin{case}{2.4.1.3}
%c
$(x\ y\ v)_\mathcal{B}\wedge (x\ w\ y)_{\mathcal{B} - v}\overset{\text{f.r.p.}}{\Rightarrow} (w\ y\ v)_\mathcal{B}$, contradicting $(v\ w\ y)_{\mathcal{B} - x}$.
%Relations $(x\ y\ v)_\mathcal{B}$ and $(x\ w\ y)_{\mathcal{B} - v}$ imply $(w\ y\ v)_\mathcal{B}$ by f.r.p., contradicting $(v\ w\ y)_{\mathcal{B} - x}$.
\end{case}
%Relations $(w\ x\ u)_\mathcal{B}$ and $(v\ w\ u)_\mathcal{B}$ imply $(v\ x\ u)_\mathcal{B}$ by f.r.p., a contradiction.
\end{wcase}

\begin{wcase}{2.4.2}{$i = 2$}
$(x\ w\ u)_{\mathcal{B} - v}\wedge (w\ v\ u)_{\mathcal{B} - x}\overset{\text{f.r.p.}}{\Rightarrow} (x\ v\ u)_\mathcal{B}\ \lightning$
%Relations $(x\ w\ u)_{\mathcal{B} - v}$ and $(w\ v\ u)_{\mathcal{B} - x}$ imply $(x\ v\ u)_\mathcal{B}$ by f.r.p., a contradiction.
%As $v, x$ and $z$ are collinear, one of the following relations hold:
%\begin{itemize}%[label={\alph*.}, ref={\alph*}]
%\item\textsc{Case} 2.4.2.1 $(v\ x\ z)_\mathcal{B}$;
%\item\textsc{Case} 2.4.2.2 $(x\ v\ z)_\mathcal{B}$;
%\item\textsc{Case} 2.4.2.3 $(x\ z\ v)_\mathcal{B}$.
%\end{itemize}
%
%\begin{case}{2.4.2.1}%a
%Relations $(v\ x\ z)_\mathcal{B}$ and $(x\ u\ z)_\mathcal{B}$ imply $(v\ x\ u)_\mathcal{B}$ by f.r.p., a contradiction.
%\end{case}
%
%\begin{case}{2.4.2.2}%b
%Relations $(x\ v\ z)_\mathcal{B}$ and $(v\ y\ z)_\mathcal{B}$ imply $(x\ y\ z)_\mathcal{B}$ by f.r.p., a contradiction.
%\end{case}
%
%\begin{case}{2.4.2.3}%c
%Relations $(x\ z\ v)_\mathcal{B}$ and $(x\ u\ z)_\mathcal{B}$ imply $(x\ u\ v)_\mathcal{B}$ by f.r.p., a contradiction.
%\end{case}
\end{wcase}

\begin{wcase}{2.4.3}{$i = 3$}
$(x\ y\ u)_{\mathcal{B} - v}\wedge (y\ v\ u)_{\mathcal{B} - x}\overset{\text{f.r.p.}}{\Rightarrow} (x\ v\ u)_\mathcal{B}\ \lightning$
%Relations $(x\ y\ u)_{\mathcal{B} - v}$ and $(y\ v\ u)_{\mathcal{B} - x}$ imply $(x\ v\ u)_\mathcal{B}$ by f.r.p., a contradiction.
%%Relations $(y\ v\ u)_\mathcal{B}$ and $(y\ u\ x)_\mathcal{B}$ imply $(v\ u\ x)_\mathcal{B}$ by f.r.p., a contradiction.
\end{wcase}

\begin{wcase}{2.4.4}{$i = 4$}
$(u\ z\ x)_{\mathcal{B} - v}\wedge (u\ v\ z)_{\mathcal{B} - x}\overset{\text{f.r.p.}}{\Rightarrow} (u\ v\ x)_\mathcal{B}\ \lightning$
%Relations $(u\ z\ x)_{\mathcal{B} - v}$ and $(u\ v\ z)_{\mathcal{B} - x}$ imply $(u\ v\ x)_\mathcal{B}$ by f.r.p., a contradiction.
%%Relations $(u\ v\ z)_\mathcal{B}$ and $(x\ u\ z)_\mathcal{B}$ imply $(x\ u\ v)_\mathcal{B}$ by f.r.p., a contradiction.
\end{wcase}

\begin{wcase}{2.4.5}{$i = 5$}
As the triple $\{v, x, y\}$ is collinear, one of the following cases hold:
\begin{itemize}
%[label={\alph*.}, ref={\alph*}]
\item\textsc{Case} 2.4.5.1: $(v\ x\ y)_\mathcal{B}$;
\item\textsc{Case} 2.2.5.2: $(x\ v\ y)_\mathcal{B}$;
\item\textsc{Case} 2.2.5.3: $(x\ y\ v)_\mathcal{B}$.
\end{itemize}

\begin{case}{2.4.5.1}
%a
$(v\ x\ y)_\mathcal{B}\wedge (x\ w\ y)_{\mathcal{B} - v}\overset{\text{f.r.p.}}{\Rightarrow} (v\ w\ y)_\mathcal{B}$, contradicting $(w\ y\ v)_{\mathcal{B} - x}$.
%Relations $(v\ x\ y)_\mathcal{B}$ and $(x\ w\ y)_{\mathcal{B} - v}$ imply $(v\ w\ y)_\mathcal{B}$ by f.r.p., contradicting $(w\ y\ v)_{\mathcal{B} - x}$.
\end{case}

\begin{case}{2.4.5.2}
%b
$(x\ v\ y)_\mathcal{B}\wedge (y\ u\ v)_{\mathcal{B} - x}\overset{\text{f.r.p.}}{\Rightarrow} (x\ v\ u)_\mathcal{B}\ \lightning$
%Relations $(x\ v\ y)_\mathcal{B}$ and $(y\ u\ v)_{\mathcal{B} - x}$ imply $(x\ v\ u)_\mathcal{B}$ by f.r.p., a contradiction.
\end{case}

\begin{case}{2.4.5.3}
%c
$(x\ y\ v)_\mathcal{B}\wedge (y\ u\ v)_{\mathcal{B} - x}\overset{\text{f.r.p.}}{\Rightarrow} (x\ u\ v)_\mathcal{B}\ \lightning$
%Relations $(x\ y\ v)_\mathcal{B}$ and $(y\ u\ v)_{\mathcal{B} - x}$ imply $(x\ u\ v)_\mathcal{B}$ by f.r.p., a contradiction.
\end{case}
%Relations $(w\ x\ u)_\mathcal{B}$ and $(w\ u\ v)_\mathcal{B}$ imply $(x\ u\ v)_\mathcal{B}$ by f.r.p., a contradiction.
\end{wcase}
\end{case}
\end{case}

\begin{case}{3}
%III
In this case, $n = 6$, $X =\{u, v, w, x, y, z\}$ and there are two, disjoint triangles in $\mathcal{B}$: $T =\{x, y, z\}$ and $R =\{u, v, w\}$ (see Figure \ref{SFthyp3}). Now, for any point $p\in X$, $\mathcal{B} - p$ is a quasilinear betweenness structure of order $5$.

\begin{wcase}{3.1}{$\mathcal{B}$ contains an ordered substructure on $4$ points}
%A
We can suppose without loss of generality that $(\mathcal{B} - x) - u$ is an ordered substructure. We know from Case 1 that there are two betweenness structures, $\mathcal{R}_{5, 1}^4$ and $\mathcal{S}_5^4$, of order $5$ that contain an ordered substructure on $4$ points. Thus, considering $\mathcal{B} - x$ and $\mathcal{B} - u$, there are three cases up to symmetry:
\begin{itemize}
%[label={\roman*.}, ref={\roman*}]
\item\textsc{Case} 3.1.1 $\mathcal{B} - x\simeq\mathcal{B} - u\simeq\mathcal{R}_{5, 1}^4$;
\item\textsc{Case} 3.1.2 $\mathcal{B} - x\simeq\mathcal{R}_{5, 1}^4$, $\mathcal{B} - u\simeq\mathcal{S}_5^4$;
\item\textsc{Case} 3.1.3 $\mathcal{B} - x\simeq\mathcal{B} - u\simeq \mathcal{S}_5^4$.
\end{itemize}

\begin{figure}
\centering
\begin{subfigure}[b]{\textwidth}
\centering
\includegraphics{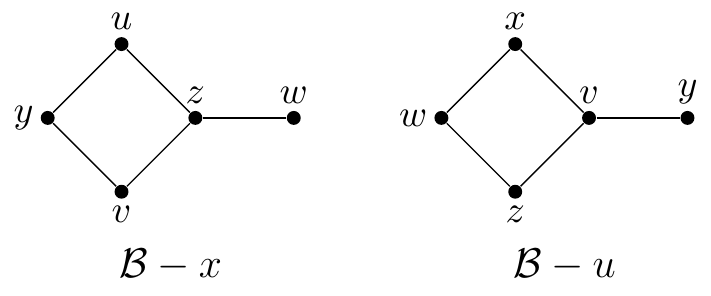}
%\begin{tikzpicture}[scale=1, pont/.style={circle, fill=black, inner sep=0.5mm},
%mybox/.style={rectangle, rounded corners, draw=black, fill=white, inner sep=0.5mm}]
%\begin{scope}[shift={(0,0)}]
%\node (A) at (0,0) [pont] {};
%\node (B) at (1.5,0) [pont] {};
%\node (C) at (2.5,0) [pont] {};
%\node (D) at (0.75,0.75) [pont] {};
%\node (E) at (0.75,-0.75) [pont] {};
%
%\draw (A)--(D) node {};
%\draw (A)--(E) node {};
%\draw (D)--(B) node {};
%\draw (E)--(B) node {};
%\draw (B)--(C) node {};
%
%\draw (A) node[anchor=east] {\small$y$};
%\draw (B) node[anchor=south] {\small$z$};
%\draw (C) node[anchor=south] {\small$w$};
%\draw (D) node[anchor=south] {\small$u$};
%\draw (E) node[anchor=north] {\small$v$};
%
%\draw (1.25,-1.5) node {$\mathcal{B} - x$};
%\end{scope}
%
%\begin{scope}[shift={(4,0)}]
%\node (A) at (0,0) [pont] {};
%\node (B) at (1.5,0) [pont] {};
%\node (C) at (2.5,0) [pont] {};
%\node (D) at (0.75,0.75) [pont] {};
%\node (E) at (0.75,-0.75) [pont] {};
%
%\draw (A)--(D) node {};
%\draw (A)--(E) node {};
%\draw (D)--(B) node {};
%\draw (E)--(B) node {};
%\draw (B)--(C) node {};
%
%\draw (A) node[anchor=east] {\small$w$};
%\draw (B) node[anchor=south] {\small$v$};
%\draw (C) node[anchor=south] {\small$y$};
%\draw (D) node[anchor=south] {\small$x$};
%\draw (E) node[anchor=north] {\small$z$};
%
%\draw (1.25,-1.5) node {$\mathcal{B} - u$};
%\end{scope}
%\end{tikzpicture}
\caption{Case 3.1.1}\label{SFc31}
\end{subfigure}%
\vspace{1em}
\begin{subfigure}[b]{\textwidth}
\centering
\includegraphics{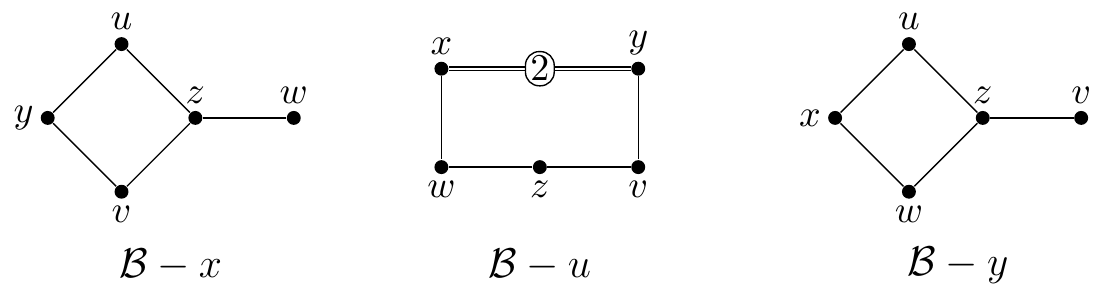}
\caption{Case 3.1.2}\label{SFc32}
\end{subfigure}%
\vspace{1em}
\begin{subfigure}[b]{\textwidth}
\centering
\includegraphics{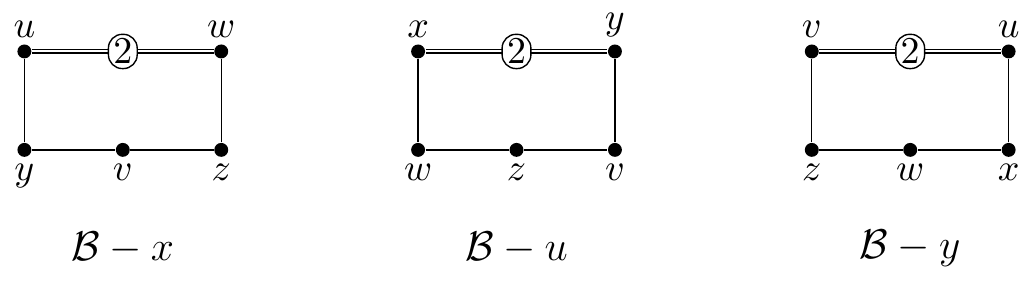}
\caption{Case 3.1.3}\label{SFc33}
\end{subfigure}%
\caption{Cases 3.1.1-3.1.3 in the proof of Lemma \ref{Lsmallgr}}%\label{Fc3}
\end{figure}

\begin{case}{3.1.1}
%i
Since $R$ is a triangle and $(\mathcal{B} - x) - u$ is ordered, we can suppose without loss of generality that $\mathcal{B} - x$ is induced by the corresponding graph shown in Figure \ref{SFc31}.
Similarly, as $(\mathcal{B} - u) - x$ is ordered, $T$ is a triangle and $(y\ v\ z\ w)_{\mathcal{B} - x}$ holds by the previous argument, $\mathcal{B} - u$ must be induced by the corresponding graph shown in Figure \ref{SFc31}.
Next, observe the following.

\begin{clm}\label{C4equ}$ $
\begin{enumerate}
%[label={(\alph*)}, ref={\emph{\alph*}}]
\item $(x\ y\ u)_\mathcal{B}\Leftrightarrow (x\ v\ u)_\mathcal{B}\Leftrightarrow (x\ z\ u)_\mathcal{B}\Leftrightarrow (x\ w\ u)_\mathcal{B}$;\label{E4equ1}
\item $\neg (y\ x\ u)_\mathcal{B}\wedge\neg (x\ u\ v)_\mathcal{B}\wedge\neg (z\ x\ u)_\mathcal{B}\wedge\neg (x\ u\ w)_\mathcal{B}$.\label{E4equ2}
\end{enumerate}
\end{clm}
\begin{prf}

\ref{E4equ1}. From Part \ref{E4equ1} we only show that $(x\ y\ u)_\mathcal{B}\Leftrightarrow (x\ v\ u)_\mathcal{B}$. Equivalences $(x\ v\ u)_\mathcal{B}\Leftrightarrow (x\ z\ u)_\mathcal{B}$ and $(x\ z\ u)_\mathcal{B}\Leftrightarrow (x\ w\ u)_\mathcal{B}$ can be proved in a similar way.

Relations $(x\ y\ u)_{\mathcal{B}}$ and $(x\ v\ y)_{\mathcal{B} - u}$ imply $(x\ v\ u)_\mathcal{B}$ by f.r.p. Conversely, $(x\ v\ u)_{\mathcal{B}}$ and $(v\ y\ u)_{\mathcal{B} - x}$ imply $(x\ y\ u)_\mathcal{B}$.

\ref{E4equ2}. From Part \ref{E4equ2} we prove $\neg (y\ x\ u)_\mathcal{B}$ as an example. Relation $(y\ u\ z)_{\mathcal{B} - x}$ holds, so $(y\ x\ u)_\mathcal{B}$ would imply $(y\ x\ z)_\mathcal{B}$ by f.r.p., a contradiction.
$\square$
\end{prf}

Now, if all of the four betweennesses in Part \ref{E4equ1} of Claim \ref{C4equ} hold, then with the $14$ betweennesses from $\mathcal{B} - x$ and $\mathcal{B} - u$, we have altogether $18$ non-trivial betweennesses that form a betweenness structure isomorphic to $\mathcal{A}_6^2$.
%induced by the second six point graph in Figure \ref{Fsmallgr}.

Otherwise, Part \ref{E4equ2} of Claim \ref{C4equ} yields $(x\ u\ y)_\mathcal{B}$, $(u\ x\ v)_\mathcal{B}$, $(x\ u\ z)_\mathcal{B}$ and $(u\ x\ w)_\mathcal{B}$ because the underlying triples must be collinear. Now, with the 14 betweennesses from $\mathcal{B} - x$ and $\mathcal{B} - u$, we have altogether $18$ non-trivial betweennesses that form a betweenness structure isomorphic to $\mathcal{A}_6^3$
% induced by the third six point graph in Figure \ref{Fsmallgr}.
\end{case}

\begin{case}{3.1.2}
%ii
As $R =\{u, v, w\}$ is a triangle and $(\mathcal{B} - x) - u$ is ordered, we can assume without loss of generality that $\mathcal{B} - x$ is induced by the corresponding graph in Figure \ref{SFc32}.
Similarly, as $(\mathcal{B} - u) - x$ is ordered and $(y\ v\ z\ w)_{\mathcal{B} - x}$ holds by the previous argument, $\mathcal{B} - u$ is induced by the corresponding graph in Figure \ref{SFc32}.

%Consider $\mathcal{B} - y$, which is a quasilinear betweenness structure of order $5$.
Observe that $(\mathcal{B} - x) - y$ is induced by a star and $(\mathcal{B} - u) - y$ is induced by a path, hence, we obtain from $(x\ w\ z\ v)_{\mathcal{B} - u}$ that $$\mathcal{B} - y\simeq\mathcal{R}_{5, 1}^4$$ and its adjacency graph is the one shown in Figure \ref{SFc32}.
%From this graph, we immediately obtain $$(x\ u\ z)_\mathcal{B}, (x\ u\ v)_\mathcal{B}\text{ and }(u\ x\ w)_\mathcal{B}.$$
Finally, $(u\ x\ w)_{\mathcal{B} - y}$ and $(y\ u\ w)_{\mathcal{B} - x}$ imply $(y\ u\ x)_\mathcal{B}$ by f.r.p. Together with the $17$ betweennesses from $\mathcal{B} - x$, $\mathcal{B} - u$ and $\mathcal{B} - y$, we have altogether $18$ non-trivial betweennesses that form a betweenness structure isomorphic to $\mathcal{A}_6^2$.
% induced by the second $6$ point graph shown in Figure \ref{Fsmallgr}.
\end{case}

\begin{case}{3.1.3}
%iii
As $R$ is a triangle and $(\mathcal{B} - x) - u$ is ordered, we can suppose without loss of generality that the first graph in Figure \ref{SFc33} induces $\mathcal{B} - x$.
Similarly, as $(\mathcal{B} - u) - x$ is ordered and $(y\ v\ z\ w)_{\mathcal{B} - x}$ holds by the previous argument, $\mathcal{B} - u$ is induced by the corresponding graph in Figure \ref{SFc33}.
It is also true that $$\mathcal{B} - y\simeq\mathcal{S}_5^4,$$ otherwise we would be back to Case 3.1.2. Now, because of $(x\ w\ z\ v)_{\mathcal{B} - u}$, $\mathcal{B} - y$ is induced by the third graph in Figure \ref{SFc33}.
Finally, $(y\ u\ w)_{\mathcal{B} - x}$ and $(u\ x\ w)_{\mathcal{B} - y}$ imply $(y\ u\ x)_\mathcal{B}$ by f.r.p. Together with the $17$ betweennesses from $\mathcal{B} - x$, $\mathcal{B} - u$ and $\mathcal{B} - y$, we have $18$ non-trivial betweennesses that form a betweenness structure isomorphic to $\mathcal{A}_6^1\simeq\mathcal{C}_6$.
% induced by a cycle of lenght $6$ (Figure \ref{Fsmallgr}).
\end{case}
\end{wcase}

\begin{wcase}{3.2}{$\mathcal{B}$ does not contain any ordered substructures on $4$ points}
%B
In this case, for all $p\in X$, $$\mathcal{B} - p\simeq\mathcal{B}(K_{2, 3}).$$ and the unique triangle must be the class of size $3$ of $K_{2, 3}$. Thus, we can determine the spanner graphs of all the substructures $\mathcal{B} - p$, $p\in X$, from which we obtain all betweennesses of $\mathcal{B}$. It follows that $\mathcal{B}$ is isomorphic to $\mathcal{A}_6^4\simeq\mathcal{B}(K_{3, 3})$.
\end{wcase}

To summarize Case 3, $\mathcal{B}$ is induced by a graph isomorphic to one of the $6$-vertex graphs in Figure \ref{Fsmallgr}.
\end{case}

\begin{case}{4}
%IV
In this case, $n = 7$, $X =\{p, u, v, w, x, y, z\}$ and there are three triangles in $\mathcal{B}$ that intersect in one point: $\{x, y, z\},\{x, u, v\}$ and $\{x, p, w\}$ (see Figure \ref{SFthyp4}). Now, $\mathcal{B} - y$ is a quasilinear betweenness structure of order $6$ with two triangles intersecting in $1$ point, hence, it belongs to Case 2. However, we have already shown that Case 2 is impossible. $\square$
\end{case}

\begin{figure}[t]
\centering
\begin{subfigure}[b]{0.3\textwidth}
\centering
\includegraphics{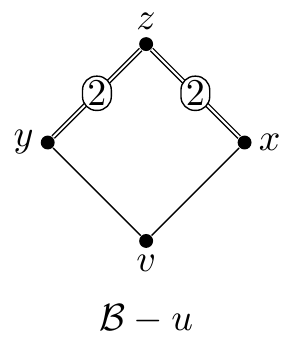}
%\begin{tikzpicture}[pont/.style={circle, fill=black, inner sep=0.5mm},
%mybox/.style={rectangle, rounded corners, draw=black, fill=white, inner sep=0.5mm}]
%
%\node (A) at (0,0) [pont] {};
%\node (B) at (1,1) [pont] {};
%\node (C) at (0,2) [pont] {};
%\node (D) at (-1,1) [pont] {};
%
%\draw (A)--(B) node {};
%\draw[double] (B)--(C) node[midway, mybox] {\small$2$};
%\draw[double] (C)--(D) node[midway, mybox] {\small$2$};
%\draw (D)--(A) node {};
%
%\draw (A) node[anchor=north] {\small$v$};
%\draw (B) node[anchor=west] {\small$x$};
%\draw (C) node[anchor=south] {\small$z$};
%\draw (D) node[anchor=east] {\small$y$};
%
%\draw (0,-0.5) node[anchor=north] {\small$\mathcal{B} - u$};
%\end{tikzpicture}
\caption{Case 1}\label{SFn5a}
\end{subfigure}
\begin{subfigure}[b]{0.3\textwidth}
\centering
\includegraphics{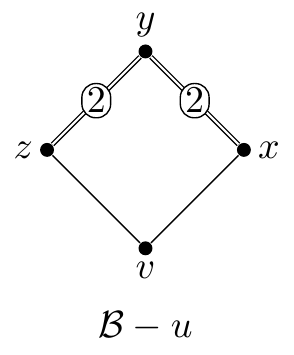}
%\begin{tikzpicture}[pont/.style={circle, fill=black, inner sep=0.5mm},
%mybox/.style={rectangle, rounded corners, draw=black, fill=white, inner sep=0.5mm}]
%
%\node (A) at (0,0) [pont] {};
%\node (B) at (1,1) [pont] {};
%\node (C) at (0,2) [pont] {};
%\node (D) at (-1,1) [pont] {};
%
%\draw (A)--(B) node {};
%\draw[double] (B)--(C) node[midway, mybox] {\small$2$};
%\draw[double] (C)--(D) node[midway, mybox] {\small$2$};
%\draw (D)--(A) node {};
%
%\draw (A) node[anchor=north] {\small$v$};
%\draw (B) node[anchor=west] {\small$x$};
%\draw (C) node[anchor=south] {\small$y$};
%\draw (D) node[anchor=east] {\small$z$};
%
%\draw (0,-0.5) node[anchor=north] {\small$\mathcal{B} - u$};
%\end{tikzpicture}
\caption{Case 2}\label{SFn5b}
\end{subfigure}%

\vspace{1em}
\begin{subfigure}[b]{1.0\textwidth}
\centering
\includegraphics{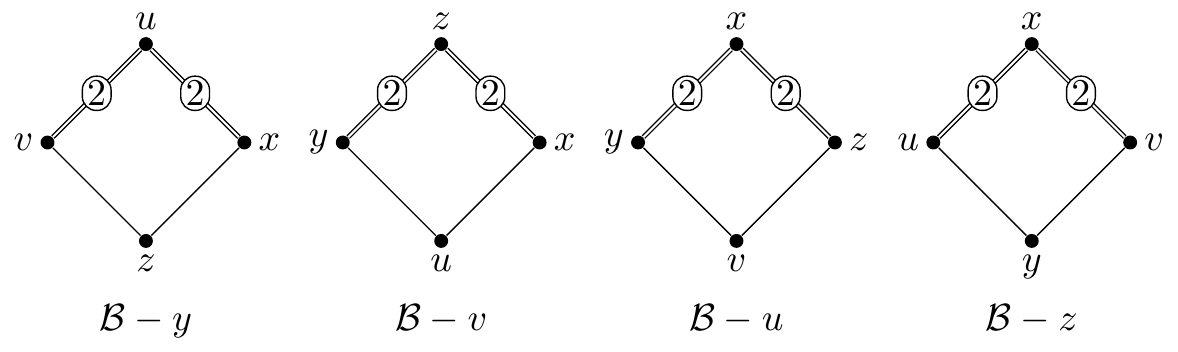}
\caption{Case 3}\label{SFn5c}
\end{subfigure}%
\caption{Cases 1-3 in the proof of Claim \ref{Cn5}}%\label{Fn5}
\end{figure}

%\begin{clm}\label{Cn5}
%Let $\mathcal{B}\in B(5, 2)$ be a regular betweenness structure. Then $\mathcal{H}(\mathcal{B})$ is a tight star.
%\end{clm}

\section{Proof of Claim \ref{Cn5}}

%\begin{prfof}{Claim \ref{Cn5}}
%\begin{prf}
Suppose to the contrary that $\mathcal{H}(\mathcal{B})$ is not a tight star. Then clearly $|R\cap T| = 1$. Let $x, y, z, u$ and $v$ be the points of $\mathcal{B}$ such that $R =\{x, u, v\}$ and $T =\{x, y, z\}$.
Since $\mathcal{B}$ is regular by assumption, $\mathcal{B} - x$ must be ordered, so one of the following cases hold up to symmetry:
\begin{itemize}
\item\textsc{Case} 1: $(u\ v\ y\ z)_\mathcal{B}$;
\item\textsc{Case} 3: $(u\ y\ z\ v)_\mathcal{B}$;
\item\textsc{Case} 2: $(u\ y\ v\ z)_\mathcal{B}$.

\end{itemize}

Before we analyze these cases, observe that for any point $p\neq x$, $\mathcal{B} - p$ is a quasilinear betweenness structure on $4$ points, thus, it is induced by one of the two $4$-vertex graphs in Figure \ref{Fsmallgr}. Further, notice that each betweenness of $\mathcal{A}_4^1$ has its middle point outside of the unique triangle.

\begin{case}{1}
Since $T$ is the only triangle of $\mathcal{B} - u$ and $(v\ y\ z)_\mathcal{B}$ holds by the case's assumption, $\mathcal{B} - u$ must be induced by the graph in Figure \ref{SFn5a} and consequently, $(v\ x\ z)_{\mathcal{B} - u}$ holds. However, this and $(u\ v\ z)_\mathcal{B}$ imply $(u\ v\ x)_\mathcal{B}$ by f.r.p., in contradiction with $R$ being a triangle.
% $(x\ v\ y)_\mathcal{B}$ holds.
%We obtain by symmetry that $\mathcal{B} - v$ is also induced by the corresponding graph in Figure \ref{SFn5a}, which implies $(x\ y\ v)_\mathcal{B}$, in contradiction with our previous observation.
\end{case}

\begin{case}{2}
Similarly to the previous case, $\mathcal{B} - u$ is induced by the corresponding graph in Figure \ref{SFn5b}, so $(y\ x\ v)_{\mathcal{B} - u}$ holds. However, this and $(u\ y\ v)_\mathcal{B}$ imply $(u\ x\ v)_\mathcal{B}$ by f.r.p. in contradiction with $R$ being a triangle.
\end{case}

\begin{case}{3}
As $(u\ v\ z)_\mathcal{B}$ holds by the case's assumption, $\mathcal{B} - y$ and $\mathcal{B} - v$ must be induced by the corresponding graphs in Figure \ref{SFn5c} and thus $(x\ z\ v)_{\mathcal{B} - y}$ and $(x\ u\ y)_{\mathcal{B} - v}$ hold.
From these betweennesses and the case's assumption, we obtain that $\mathcal{B} - u$ and $\mathcal{B} - z$ are induced by the corresponding graphs in Figure \ref{SFn5c}. But now, $(x\ y\ v)_{\mathcal{B} - u}$ contradicts $(x\ v\ y)_{\mathcal{B} - z}$. $\square$
\end{case}
%\end{prfof}

\begin{figure}[t]
\centering
\includegraphics{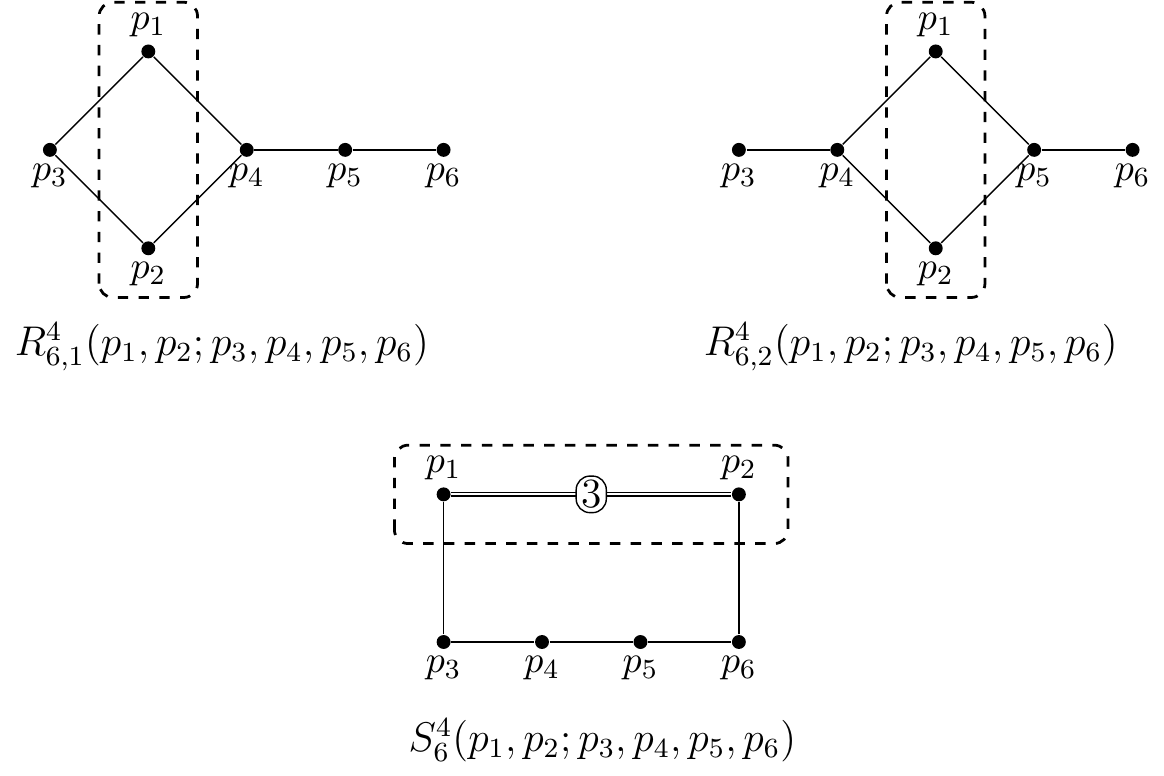}
\caption{Vertex-labeled spanner graphs of quasilinear betweenness structures on $6$ points. The kernel of each triangle hypergraph is encircled by dashed line.}\label{Fapi}
\end{figure}

\section{Proof of Claim \ref{Cmetr7}}

%\begin{prfof}{Claim \ref{Cmetr7}}
For six distinct points $p_1, p_2, p_3, p_4, p_5, p_6$, let $\mathcal{R}_{6, 1}^4(p_1, p_2; p_3, p_4, p_5, p_6)$, $\mathcal{R}_{6, 2}^4(p_1,\allowbreak p_2;\allowbreak p_3,\allowbreak p_4,\allowbreak p_5,\allowbreak p_6)$ and $\mathcal{S}_6^4(p_1, p_2; p_3, p_4, p_5, p_6)$ denote the betweenness structures induced by the graphs shown in Figure \ref{Fapi}.

\begin{case}{\ref{Emetr71}}
It is clear that $\tilde{\mathcal{H}}_1$ is metrizable as $\mathcal{H}(\mathcal{T}_{7, 1})\simeq\tilde{\mathcal{H}}_1$. Next, let $\mathcal{A}$ be a betweenness structure such that $\mathcal{H}(\mathcal{A}) =\tilde{\mathcal{H}}_1$ and label the points of $\tilde{\mathcal{H}}_1$ as in Figure \ref{SFmetr71}. We show that $\mathcal{A}\simeq\mathcal{T}_{7, 1}$.

It can be easily seen that $\mathcal{A} - x$, $\mathcal{A} - w$, $\mathcal{A} - y$ and $\mathcal{A} - z$ are all quasilinear betweenness structures of order $6$ and their triangle hypergraphs are tight starts, hence, they are isomorphic to either $\mathcal{R}_{6, 1}^4$, $\mathcal{R}_{6, 2}^4$ or $\mathcal{S}_6^4$ by Lemma \ref{Lgen3}. Observe the following.

\begin{obs}\label{O6qlin}
Let $\mathcal{B}\in B(6, 2)$ be a betweenness structure on ground set $X =\{p_1, p_2, p_3, p_4, p_5, q\}$ such that $(p_1\ p_2\ p_3\ p_4\ p_5)_\mathcal{B}$ holds and $\mathcal{H}(\mathcal{B})$ is a tight star with kernel $K =\{p_i, q\}$. Then
\begin{enumerate}
%[label={(\roman*)}, ref={\roman*}]
\item $i = 1\Rightarrow\mathcal{B} =\mathcal{S}_6^4(p_1, q; p_2, p_3, p_4, p_5)$;
\item $i = 2\Rightarrow\mathcal{B} =\mathcal{R}_{6, 1}^4(q, p_2; p_1, p_3, p_4, p_5)$;
\item $i = 3\Rightarrow\mathcal{B} =\mathcal{R}_{6, 2}^4(q, p_3; p_1, p_2, p_4, p_5)$.
\end{enumerate}
(Cases $i = 4$ and $i = 5$ can be obtained by symmetry.)
\end{obs}

Suppose first that one of the betweenness structures $\mathcal{A} - q'$, $q'\in\{x, y, z, w\}$ is isomorphic to $\mathcal{R}_{6, 1}^4$. We may assume that $q' = x$. Then $$\mathcal{A} - x =\mathcal{R}_{6, 1}^4(y, z; p_3, p_4, p_5, p_6)$$ where $$\{p_3, p_4\} =\{u, w\}\emph{ and }\{p_5, p_6\} =\{p, v\},$$ so there are four possibilities.

If $\mathcal{A} - x =\mathcal{R}_{6, 1}^4(y, z; w, u, p, v)$, then $(w\ z\ u\ p\ v)_\mathcal{A}$ holds, hence, $\mathcal{A} - y =\mathcal{S}_6^4(w, x; z, u, p, v)$ by Observation \ref{O6qlin}.
Further, since $\{y, z\} =\ker(\mathcal{H}(\mathcal{A} - w))$ and $(y\ u\ p\ v)_{\mathcal{A} - x}$ and $(z\ u\ p\ v)_{\mathcal{A} - x}$ hold, $\mathcal{A} - w$ must be $\mathcal{R}_{6, 1}^4(y, z; x, u, p, v).$ Now, $(x\ u\ p)_{\mathcal{A} - w}$ holds in contradiction with $(u\ p\ x)_{\mathcal{A} - y}$.

If $\mathcal{A} - x =\mathcal{R}_{6, 1}^4(y, z; w, u, v, p)$, then similarly to the previous case, $\mathcal{A} - y =\mathcal{S}_6^4(w, x; z, u, v, p)$ and $\mathcal{A} - w =\mathcal{R}_{6, 1}^4(y, z; x, u, v, p)$, which contradict on the triple $\{u, p, x\}$.

If $\mathcal{A} - x =\mathcal{R}_{6, 1}^4(y, z; u, w, p, v)$, then $\mathcal{A} - y$ must be $\mathcal{R}_{6, 2}^4(x, w; u, z, p, w)$, which contradicts the fact that $\{x, w, p\}$ is a triangle.

Lastly, we show that if $\mathcal{A} - x =\mathcal{R}_{6, 1}^4(y, z; u, w, v, p)$, then $\mathcal{A}\simeq\mathcal{T}_{7, 1}$. It is easy to see by Observation \ref{O6qlin} that $\mathcal{A} - y =\mathcal{R}_{6, 2}^4(x, w; u, z, v, p)$ and $\mathcal{A} - z =\mathcal{R}_{6, 2}^4(x, w; u, y, v, p)$. Both implies that $\mathcal{A} - w =\mathcal{R}_{6, 1}^4(y, z; u, x, v, p)$. These substructures are consistent with one another and cover all collinear triples of $\mathcal{A}$. It is easy to see now that $\mathcal{A}$ is induced by the vertex-labeled $T_{6, 1}$ in Figure \ref{SFtpi1}.

Next, suppose that none of the betweenness structures $\mathcal{A} - q'$, $q'\in\{x, y, z, w\}$ is isomorphic to $\mathcal{R}_{6, 1}^4$, but one of them, for example $\mathcal{A} - x$, is isomorphic to $\mathcal{R}_{6, 2}^4$. Then $$\mathcal{A} - x =\mathcal{R}_{6, 2}^4(y, z; p_3, p_4, p_5, p_6)$$ where $$\{p_3, p_6\} =\{p, v\}\emph{ and }\{p_4, p_5\} =\{u, w\},$$ so, there are two possibilities up to symmetry.
If $\mathcal{A} - x =\mathcal{R}_{6, 2}^4(y, z; p, u, w, v)$, then $\mathcal{A} - y =\mathcal{R}_{6, 1}^4(x, w; v, z, u, p)$, while if $\mathcal{A} - x =\mathcal{R}_{6, 2}^4(y, z; p, w, u, v)$, then $\mathcal{A} - y =\mathcal{R}_{6, 1}^4(x, w; p, z, u, v)$ by Observation \ref{O6qlin}. In both cases, $\mathcal{A} - y$ contradicts our previous assumption on $\mathcal{A} - q'$.

Finally, suppose that all of the betweenness structures $\mathcal{A} - q'$, $q'\in\{x, y, z, w\}$ are isomorphic to $\mathcal{S}_6^4$.
Then $$\mathcal{A} - x =\mathcal{S}_6^4(y, z; p_3, p_4, p_5, p_6)$$ where $$\{p_3, p_6\} =\{u, w\}\emph{ and }\{p_4, p_5\} =\{p, v\},$$ so there are two possibilities up to symmetry.
If $\mathcal{A} - x =\mathcal{S}_6^4(y, z; u, p, v, w)$, then $\mathcal{A} - y =\mathcal{R}_{6, 1}^4(x, w; z, v, p, u)$, while if $\mathcal{A} - x =\mathcal{S}_6^4(y, z; u, v, p, w)$, then $\mathcal{A} - y =\mathcal{R}_{6, 1}^4(x, w; z, p, v, u)$ by Observation \ref{O6qlin}.  In both cases, $\mathcal{A} - y$ contradicts our previous assumption on $\mathcal{A} - q'$.

In summary, we can conclude that $\mathcal{A}\simeq\mathcal{T}_{7, 1}$.
\end{case}

\begin{case}{\ref{Emetr72}}
Suppose to the contrary that $\tilde{\mathcal{H}}_2$ is metrizable and let $\mathcal{A}$ be a betweenness structure such that $\mathcal{H}(\mathcal{A}) =\tilde{\mathcal{H}}_2$. Label the points of $\tilde{\mathcal{H}}_2$ as in Figure \ref{SFmetr72}.
Notice that both $\mathcal{A} - x$ and $\mathcal{A} - q$ are a quasilinear betweenness structures with $2$ triangles that form a tight star. The next observation follows from Theorem \ref{Tqus1}.

\begin{obs}\label{Okern}
Let $\mathcal{A}$ be a quasilinear betweenness structure of order $n\geq 6$ such that $\mathcal{H}(\mathcal{A})$ is a tight star. Then there is exactly one cyclic line in $\mathcal{A}$ and it contains $\ker(\mathcal{H}(\mathcal{A}))$.
\end{obs}
%\begin{prf}
%Let $Y$ be the ground set of $\mathcal{A}$. Because of Theorem \ref{Tqus1}, $\mathcal{A}$ is isomorphic to either $\mathcal{R}_{n, i}^4$ ($i\in I_n^4$) or $\mathcal{S}_n^4$. It is easy to see that there is exactly one subset $Y' =\{y_1, y_2, y_3, y_4\}\subset Y$ such that the betweennesses $(y_1\ y_2\ y_3)_\mathcal{A}$, $(y_2\ y_3\ y_4)_\mathcal{A}$, $(y_3\ y_4\ y_1)_\mathcal{A}$ and $(y_4\ y_1\ y_2)_\mathcal{A}$ hold. Further, $Y'$ contains $\ker(\mathcal{H}(\mathcal{A}))$. $\square$
%\end{prf}

We obtain from Observation \ref{Okern} that there is exactly one cyclic line $L_x$ in $\mathcal{A} - x$, and it contains the kernel $\{y, z\}$. 
%and there is exactly one cyclic line $L_q$ in $\mathcal{A} - q$ and it contains the star kernel $\{y, p\}$.
Further, $L_x$ does not contain $p$ or $q$ because $\{y, z, p\}$ and $\{y, z, q\}$ are triangles. Hence, $L_x$ is a cyclic line in $\mathcal{A} - q$ as well that does not contain $p$, contradicting Observation \ref{Okern}.
\end{case}

\begin{case}{\ref{Emetr73}}
Suppose to the contrary that $\tilde{\mathcal{H}}_3$ is metrizable and let $\mathcal{A}$ be a betweenness structure such that $\mathcal{H}(\mathcal{A}) =\tilde{\mathcal{H}}_3$. Label the points of $\tilde{\mathcal{H}}_3$ as in Figure \ref{SFmetr73}.
Observe that $\mathcal{A} - x$ and $\mathcal{A} - z$ are quasilinear betweenness structures on $6$ points with $2$ triangles that form a tight star, hence, they are isomorphic to either $\mathcal{R}_{6, 1}^4$, $\mathcal{R}_{6, 2}^4$ or $\mathcal{S}_6^4$.

If $\mathcal{A} - x\simeq\mathcal{R}_{6, 1}^4$, then we can assume by symmetry that $\mathcal{A} - x =\mathcal{R}_{6, 1}^4(y, z; u, v, p, q)$. It is also easy to see by Observation \ref{O6qlin} that $\mathcal{A} - z =\mathcal{R}_{6, 1}^4(y,\allowbreak x;\allowbreak u,\allowbreak v,\allowbreak p,\allowbreak q)$, from which $(y\ u\ x)_\mathcal{A}$ follows. This is, however, impossible as $\{u, x, y\}$ was a triangle.

Next, if $\mathcal{A} - x\simeq\mathcal{R}_{6, 2}^4$, then we can assume by symmetry that $\mathcal{A} - x =\mathcal{R}_{6, 2}^4(y, z; p, u, v, q)$. It is easy to see that $\mathcal{A} - z =\mathcal{R}_{6, 2}^4(y, x; p, u, v, q)$ by Observation \ref{O6qlin}, from which the contradiction $(y\ u\ x)_{\mathcal{A}}$ follows again.

Finally, if $\mathcal{A} - x\simeq\mathcal{S}_6^4$, then we can assume by symmetry that $\mathcal{A} - x =\mathcal{S}_6^4(y, z; u, p, q, v)$. Now, we obtain from Observation \ref{O6qlin} that $\mathcal{A} - z =\mathcal{S}_6^4(y, x; u, p, q, v)$, from which $(u\ y\ x)_{\mathcal{A}}$ follows, leading to a contradiction again. $\square$
\end{case}
%\end{prfof}

\bibliographystyle{elsarticle-num}
\bibliography{fms}